\documentclass[preprint,12pt]{elsarticle}

\usepackage{fullpage}
\usepackage{graphicx}
\usepackage{subcaption}
\usepackage{xcolor}
\usepackage{multirow}
\graphicspath{ {./pics/} }
\usepackage{tikz}
\usetikzlibrary{arrows,positioning,calc}
\usepackage{amsmath,amsfonts,bm}
\DeclareMathOperator{\tr}{tr}
\usepackage{color}
\usepackage{float}
\usepackage{newfloat}
\usepackage{hyperref}
\usepackage{lineno}
\hypersetup{colorlinks,bookmarksopen,linkcolor=blue,citecolor=blue,urlcolor=blue}
\pgfdeclarelayer{bg} 
\pgfsetlayers{bg,main} 
\definecolor{myblue}{RGB}{0,115,189}

\DeclareFloatingEnvironment{algorithm}
\definecolor{mygreen}{RGB}{49,156,54}

\journal{International Journal of Solids and Structures}

\biboptions{authoryear}

\begin{document}
\begin{frontmatter}

\title{Bayesian Approach to Micromechanical Parameter Identification Using Integrated Digital Image Correlation}

\author[CTU_math]{L. Gaynutdinova}
\ead{liya.gaynutdinova@fsv.cvut.cz}
\author[TUe]{O. Rokoš}
\ead{o.rokos@tue.nl}
\author[CTU_mech]{J. Havelka}
\ead{jan.havelka@fsv.cvut.cz}
\author[CTU_math]{I. Pultarová}
\ead{ivana.pultarova@cvut.cz}
\author[CTU_mech]{J. Zeman\corref{correspondingauthor}}
\ead{jan.zeman@cvut.cz}

\address[CTU_math]{Department of Mathematics, Faculty of Civil Engineering, Czech Technical University in Prague, 166 29 Prague 6, Thákurova 2077/7, Czech Republic}
\address[CTU_mech]{Department of Mechanics, Faculty of Civil Engineeing, Czech Technical University in Prague, 166 29 Prague 6, Thákurova 2077/7, Czech Republic}
\address[TUe]{Mechanics of Materials, Department of Mechanical Engineering, Eindhoven University of Technology, P.O. Box 513, 5600 MB Eindhoven, The Netherlands}
\cortext[correspondingauthor]{Corresponding author.}

\begin{abstract}
Micromechanical parameters are essential in understanding the behaviour of materials with a heterogeneous structure, which helps to predict complex physical processes such as delamination, cracks, and plasticity. However, identifying these parameters is challenging due to micro-macro length scale differences, required high resolution, and ambiguity in boundary conditions, among others. The Integrated Digital Image Correlation (IDIC) method, a state-of-the-art full-field deterministic approach to parameter identification, is widely used but suffers from high sensitivity to boundary data errors and is limited to identification of parameters within well-posed problems. This article employs Bayesian approach to estimate micromechanical shear and bulk moduli of fibre-reinforced
composite samples under plane strain assumption, and to improve handling of boundary noise. The main purpose of this article is to quantify the effect of uncertainty in the boundary conditions in the stochastic setting. To this end, the Metropolis--Hastings Algorithm (MHA) is employed to estimate probability distributions of bulk and shear moduli and boundary condition parameters using IDIC, considering a fibre-reinforced composite sample under plane strain assumption.  The performance and robustness of the MHA are compared to two versions of deterministic IDIC method, under artificially introduced random and systematic errors in kinematic boundary conditions. Although MHA is shown to be computationally more expensive and in certain cases less accurate than the recently introduced Boundary-Enriched IDIC, it offers significant advantages, in particular being able to optimize a large number of parameters while obtaining statistical characterization as well as insights into individual parameter relationships. The paper furthermore highlights the benefits of the non-normalised approach to parameter identification with MHA (leading, within deterministic IDIC, to an ill-posed formulation), which significantly improves the robustness in handling the boundary noise.
\end{abstract}


\end{frontmatter}

\section{Introduction}
Often the only way of correctly interpreting the behavior of real materials with a heterogeneous structure is by observing them on a microscopic level, so we can account for strain localization, plasticity, delamination and cracks. The spatially heterogeneous character of these phenomena prompts the need for non-intrusive full-field measurement techniques in experimental mechanics. Digital imaging enabled the development of a highly accurate method called Digital Image Correlation (DIC) \citep{keating}, which is used to  assess the spatial transformation between two images. In practice, DIC is implemented as a computer program that automatically tracks regions of an object from one configuration to another, from which displacements can be inferred. Affordability and availability of the equipment and computer programs contributed to the popularity of this method~\citep{avril, viggiani}. 

DIC technique is often used as an input for identifying model parameters, typically material, leading to methods such as Integrated DIC (IDIC) \citep{leclerc}. The method relies on deterministic optimization of the least square difference between two or more images of a sample captured during an experiment, i.e., in the reference and a deformed configuration. This approach minimizes information losses and provides highly accurate results.

\begin{figure}[!htbp]
    \centering
	\scalebox{0.9}
	{
	\begin{tikzpicture}[>=stealth]
		\linespread{1}
		\tikzset{
			mynode/.style={inner sep=0,outer sep=0},
			myarrow/.style={myblue,thick},
		}	
		
		\begin{pgfonlayer}{bg}
		\node[mynode] (specimen) {
			\includegraphics[scale=0.23]{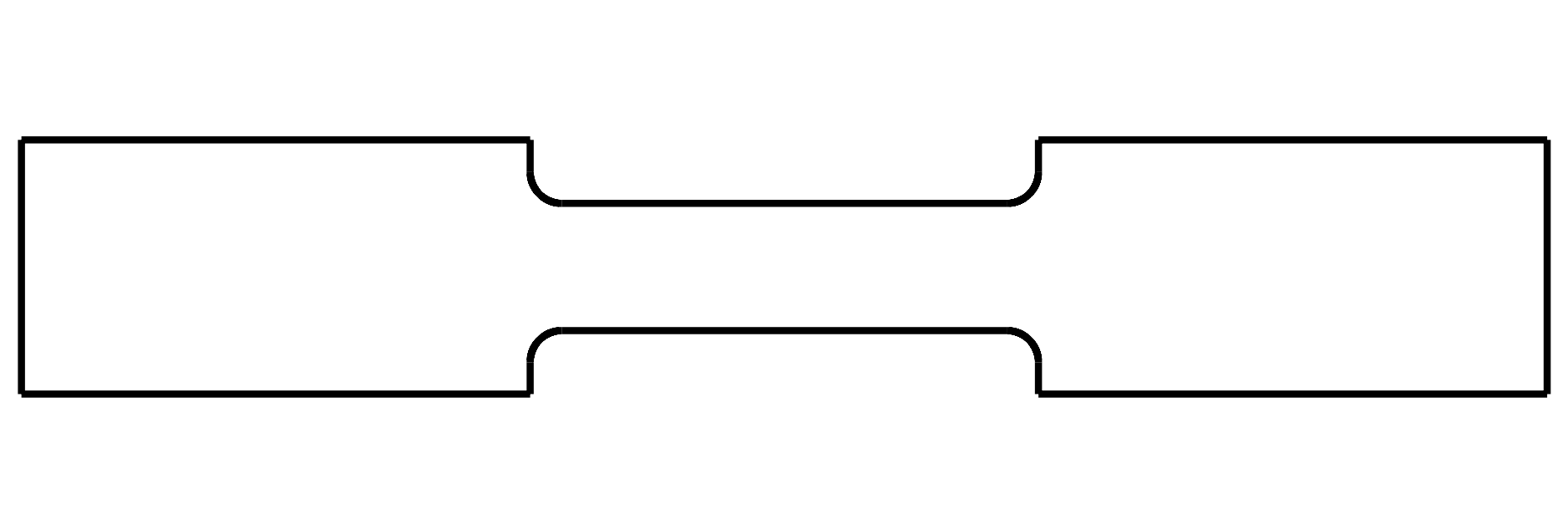}
		};
		\end{pgfonlayer}
		\node[mynode,above=1.0em of specimen] (mve) {
			\includegraphics[scale=1]{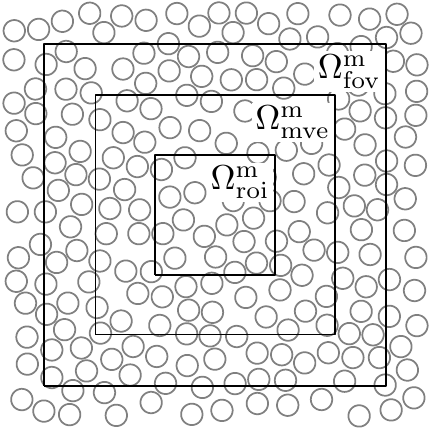}
		};
		\draw (-0.1,-0.1) rectangle (0.1,0.1);
		
		\coordinate []  (Fa) at (-2.9,0);
		\coordinate []  (Fb) at (-4.9,0);
		\node[below=0.1em of Fb,anchor=north] {$\color{myblue}{\boldsymbol{F}_\mathrm{exp},\ \boldsymbol{u}_\mathrm{D}}$};
		\coordinate []  (Fc) at (2.9,0);
		\coordinate []  (Fd) at (4.9,0);
		\node[below=0.1em of Fd,anchor=north] {$\color{myblue}{\boldsymbol{F}_\mathrm{exp},\ \boldsymbol{u}_\mathrm{D}}$};
		\draw[->,myblue,line width=0.5mm] (Fa) to (Fb);
		\draw[->,myblue,line width=0.5mm] (Fc) to (Fd);
		
		\node[mynode,left=1.0em of mve, shift={(0.0,0.5)}] (matrix) {\footnotesize matrix~};
		\node[mynode,left=1.0em of mve, shift={(0.0,-0.5)}] (inclusions) {\footnotesize inclusions~};
		\draw[->] (matrix.east) to (-2.0,3.65);
		\draw[->] (inclusions.east) to (-1.9,2.6);
		
		\begin{pgfonlayer}{bg}
		\draw[gray,dashed] (0.1,-0.1) to (mve.south east);
		\draw[gray,dashed] (-0.1,-0.1) to (mve.south west);
		\draw[gray,dashed] (0.1,0.1) to (mve.north east);
		\draw[gray,dashed] (-0.1,0.1) to (mve.north west);
		\end{pgfonlayer}
		
		
	\end{tikzpicture}}
	\caption{Scheme of the virtual experiment, microscale (top) and macroscale (bottom). Here, a macroscopic specimen is subjected to a tensile load, by either prescribing Neumann (external force $\bm{F}_{\text{exp}}$) or Dirichlet (displacement vector $\boldsymbol{u}_\mathrm{D}$) boundary conditions. Deformations of the microstructure $\Omega_{\text{mve}}^{\text{m}}$ (captured within the microscopic FOV $\Omega_{\text{fov}}^{\text{m}}$) are observed at the microscale with optical or scanning electron microscopy. The micro-images are correlated on the microscopic ROI, $\Omega_{\text{roi}}^{\text{m}}$.}
	\label{fig:experiment}
\end{figure}
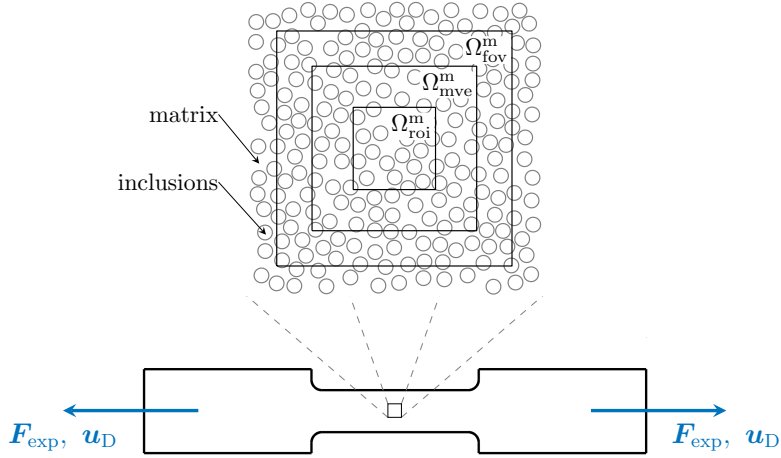

IDIC, however, comes with a set of challenges, particularly when used in a multi-scale setting. First of all, it requires a mechanical model with suitable constitutive relations and boundary conditions. Because producing a test specimen on a sufficiently small scale is involved, and the manufacturing process itself may influence the microscopic parameters, it is highly preferential to use the actual (macro-scale) product in the measuring process. To capture the microstructural displacements, the size of a digital pixel associated with DIC must be sufficiently small, as well as of the applied speckle pattern~\citep[see, e.g.,][]{Johan:2019}. While modern commercially available optical microscopes are able to provide high resolution, the problem lies in the fact that it is too time-consuming to scan the entire specimen with sufficient detail and then simulate it. For these reasons, only a subdomain of the specimen is typically scanned within a microscopic Field Of View (FOV), $\Omega_\mathrm{fov}^\mathrm{m}$, see Fig.~\ref{fig:experiment}, within which a Microstructural Volume Element~(MVE), $\Omega_{\text{mve}}^{\text{m}}$, is considered and a microstructural IDIC
model is constructed and correlated inside a microscopic Region Of Interest~(ROI), $\Omega_{\text{roi}}^{\text{m}}$. Although beneficial, this reduction brings a two-fold complication: (i) only material parameter ratios can be identified, because the Dirichlet boundary conditions are applied along the whole boundary, and hence the problem is inherently ill-posed, and (ii) those boundary conditions are not known and have to be identified (from local deformations at the boundary of $\Omega_{\text{mve}}^{\text{m}}$), since macroscopically applied boundary conditions (which are known) fall outside of the microscopic FOV, $\Omega_\mathrm{fov}^\mathrm{m}$. A potential solution are ``virtual boundaries'' \citep{kremmer}, which, however, may not be suitable for highly heterogeneous microstructures. On the other hand, high accuracy in boundary conditions prescribed to the MVE model is crucial \citep{rokos,rokos:2023}, as even small errors may significantly deteriorate accuracy of the identified parameters.

The most accurate way to establish MVE boundary conditions, according to \cite{shakoor}, is to employ Global DIC (GDIC) \citep{besnard}. In this method, the displacements are identified on the entire specimen and are subsequently interpolated as boundary conditions for the microstructural IDIC. In general, GDIC introduces: (i) kinematic smoothing effects when large elements or globally supported interpolation functions are used, and (ii) random errors when relatively small elements or locally supported interpolation functions are used. Because boundary conditions of the microstructure are kept fixed during the IDIC parameter identification procedure, the microstructural model has to compensate by adjusting its material parameters, causing the inaccurate identification of these parameters \citep{ruybalid}. The approach proposed by \cite{rokos}, referred to as the Boundary-Enriched IDIC (BE-IDIC), incorporates all Degrees Of Freedom (DOFs) associated with the virtual boundaries as IDIC DOFs. The method significantly improves accuracy of the identified parameters while maintaining robustness with respect to the image noise, as demonstrated in a case study on identifying matrix and fiber shear and bulk moduli for a fibre-reinforced composite sample under plane strain. Although the improved accuracy comes with a price of higher computational and memory requirements, the main weak point of microstructural IDIC/BE-IDIC is in dealing with the possible ill-posedness of the identification problem. For example, the choice of Dirichlet boundary conditions leads to infinite linearly dependent solutions for the material parameters. Because the aforementioned techniques rely on deterministic optimization methods, such as the Gauss--Newton algorithm, only material parameter ratios can be identified. This can typically be done by fixing one of the parameters to a predetermined value, which subsequently needs to be normalized~\citep[see][for an example and more details]{rokos:2023}. 

This article proposes a stochastic method for the parameter identification in the setting identical to \citep{rokos, rokos:2023}. Here, the Metropolis--Hastings Algorithm~(MHA), is used for the minimization process, while the IDIC is used as a full-field measurement technique. In contrast to the described deterministic methods, stochastic inversion allows to infer probability distributions of the unknown model parameters instead of single values, treating each iteration as an experimental measurement. While Bayesian inference is commonly used for mechanical parameter identification \citep{rappel,kucerova,Yue2022, Kuhn2021, THOMAS2022109630}, the effect of the uncertain boundary conditions has only been quantified in the deterministic setting. For this purpose, we choose not to employ any surrogate models for alleviating high computational costs. Additionally, we show that the Markov chain based method is able to overcome the ill-posedness of the inverse problem because the sampling depends only on the prior distribution and the previous state. The main novelties of this manuscript are therefore threefold: (i)~combination of the Bayesian approach with IDIC and MHA, (ii)~alleviation of potential problems with ill-conditioning or even ill-posedness due to many (kinematic) parameters in the deterministic version of BE-IDIC, and (iii)~the possibility to solve for ill-posed problems without normalizing micromechanical parameters, leading potentially to more accurate identifications.

The article is structured as follows: Section~\ref{section:idic} first introduces the relevant background on DIC and IDIC within the adopted two-dimensional setting. Readers already familiar with these topics or interested mostly in the stochastic approach may directly proceed to Section~\ref{section:stochastic}, where probability densities are derived first in Section~\ref{sec:bayes}, using Bayesian theory, which are converted into probability integrals for obtaining relevant statistical quantities in Section~\ref{sec:mha}, computed subsequently through a suitable sampling method (i.e., the Metropolis--Hastings algorithm, discussed at the end of Section~\ref{sec:mha}). The particular virtual experiment adopted throughout this manuscript, i.e., mechanical problem, loading and boundary conditions, macro- and micro-structural geometry, constitutive model, and identified material parameters are detailed in Section~\ref{model}. The following sections contain systematic overview of the results of the virtual numerical experiment. In particular, in Section~\ref{section:experiments} the MHA's sensitivity with respect to random and systematic errors in the boundary conditions are quantified and compared to deterministic IDIC in terms of mean values and modes. Section~\ref{section:relaxed_BC} then introduces parametrization of applied boundary conditions, and the resulting MHA with boundary DOFs is compared to BE-IDIC. Section~\ref{section:ill-posed} studies the behavior of the non-normalized MHA and compares it to the normalized version. The effect of the normalization choice in the BE-IDIC and the post-processing of the non-normalized MHA are also examined. The takeaways from the numerical experiments are finally summarized in Section~\ref{sec:summary} with the outlook on further research.

Throughout this article, scalar variables are denoted using italic font, $a$, array variables using sans serif font, $\mathsf{u}$, vectors and tensors are rendered in a boldface font, $\boldsymbol{v}$ or $\boldsymbol{A}$, single contraction is denoted~$\boldsymbol{A}\cdot\boldsymbol{v} = A_{ij}v_j$, for a second-order tensor $\boldsymbol{A}$ and vector $\boldsymbol{v}$, while ${\bm \nabla}_0$ denotes the gradient operator with respect to the reference configuration, ${\bm \nabla}_0 \boldsymbol{v} = \frac{\partial v_j}{\partial X_i} \boldsymbol{e}_i\boldsymbol{e}_j$, where $\boldsymbol{e}_i$ is a set of coordinate basis vectors. The hat ($\widehat{\bullet}$) denotes arbitrary admissible values, whereas the absence of hats indicates corresponding minimizers. 

\section{Deterministic Approach to Parameter Identification}\label{section:idic}
This section gives a brief exposition on the DIC full-field identification technique used in the numerical experiments, as well as the deterministic methods IDIC and BE-IDIC used in lieu of the benchmark.

A mechanical test is considered, as outlined in Fig.~\ref{fig:experiment}. DIC is used to assess spatial transformations before and after the specimen is deformed. A region of a photographed domain is tracked between the images, which allows to infer the displacement field upon proper regularization \citep{roux}. A camera has a static FOV, $\Omega_{\text{fov}}^{\text{m}}$, that contains the MVE, $\Omega_{\text{mve}}^{\text{m}}$,---the sub-domain modelled with FEM---which in turn contains the Region of Interest (ROI), $\Omega_{\text{roi}}^{\text{m}}$, which is used for correlating the images before and after deformation. Although the MVE and ROI may coincide, the camera's FOV is chosen such that the ROI (or MVE, if necessary) remains within it even after deformation. The images are stored as integer-valued arrays for both initial and deformed configuration, where each integer is associated with a pixel and denotes its brightness. 

DIC measurements can be used to identify a set of model parameters, i.e., to find a vector $\bm{\lambda}\in\mathbb{R}^{n_{\lambda}}$ that minimizes the difference between the values in the reference image and in the corresponding material points in the deformed image in the least squares sense, i.e.,
\begin{equation}
\bm{\lambda} \in \underset{\widehat{\bm \lambda} \in \mathbb{R}^{n_{\lambda}}}{\arg \min}\;{\cal R}_{\mathrm{dic}}(\widehat{\bm{\lambda}}).
\label{lambda}
\end{equation}
In Eq.~\eqref{lambda}, $\widehat{\bm \lambda}$ is a column matrix that stores the sought material parameters or kinematic DOFs at the boundary of MVE, $\mathbb{R}^{n_{\lambda}}$ denotes an $n_{\lambda}$-dimensional real space, and $\cal R_{\mathrm{dic}}$ is a non-convex cost function given as
\begin{equation}
{\cal R}_{\mathrm{dic}}(\widehat{\bm{\lambda}}) = \frac{1}{2} \int_{\Omega_{\mathrm{roi}}} \left[f(\bm{X}) - g(\bm{X}+\bm{u}(\bm{X},\widehat{\bm \lambda}))\right]^2 
\mathrm{d}\bm{X}, 
\label{cal_R}
\end{equation}
assuming that the brightness is conserved under this transformation. Here, $\bm{u}(\bm{X},\widehat{\bm \lambda}) = [u_1(\bm{X},\widehat{\bm \lambda}),u_2(\bm{X},\widehat{\bm \lambda})]^{\sf T}$ is an approximate displacement field depending on the set of model parameters $\widehat{\bm \lambda}$, $\bm{X}=[X_1,X_2]^{\sf T} \allowbreak \in~\Omega_{\mathrm{mve}} \subset \mathbb{R}^2$ stores the material coordinates in the reference configuration, $f(\bm{X})$ represents the initial image, whereas $g(\bm{X}+\bm{u}(\bm{X},\widehat{\bm \lambda}))$ the deformed image mapped onto the initial configuration. In a particular case of Global DIC~(GDIC), the displacement field~$\boldsymbol{u}(\boldsymbol{x},\boldsymbol{\lambda})$ is directly approximated with a continuous locally or globally supported polynomials, where parameters~$\boldsymbol{\lambda}$ store their corresponding coefficients of linear combination, cf., e.g., \citep{gdic:2012,besnard} for more details. This formulation seeks for kinematic full-field information, unlike integrated approach seeking directly for material parameters detailed below.

\subsection{Integrated Digital Image Correlation}\label{section:IDIC}

IDIC is a method initially proposed by \cite{roux} for experimental identification of typically material parameters based on Eq.~\eqref{lambda}. The displacement field ${\bm{u}}(\bm{X},\widehat{\bm \lambda})$  is obtained by solving the response of the underlying system, in our case mechanical system, which is governed by
\begin{equation}
\label{divP}
\begin{aligned}
    {\bm \nabla}_0 \cdot {\bm P}^{\sf T}({\bm u}({\bm X},\widehat{\bm \lambda}),\widehat{\bm \lambda}) &= {\bm 0}, && \bm{X}\in \Omega_{\mathrm{mve}},\\
    \bm{u}(\bm{X}) &= \bm{u}_{\partial\Omega_{\mathrm{mve}}}(\bm{X}), && \bm{X}\in \partial\Omega_{\mathrm{mve}},
\end{aligned}
\end{equation}
where ${\bm \nabla}_0$ denotes the gradient operator in the reference configuration, ${\bm P}$ is the first Piola--Kirchhoff stress tensor described by an underlying constitutive law discussed in more detail below in Section~\ref{sec:constitutive_model}, and ${\bm u}_{\partial \Omega_{\mathrm{mve}}}$ is a prescribed displacement on the boundary $\partial \Omega_{\mathrm{mve}}$, see, e.g., \cite{Tadmor} for more details on continuum mechanics. Note that because only essential boundary conditions are of interest, no tractions are prescribed to the model. The governing equation and its solution are typically discretized with the Finite Element Method (FEM), which can be substituted to Eqs.~\eqref{lambda}--\eqref{cal_R}. To minimize the objective in Eq.~\eqref{cal_R}, the standard Gauss–Newton algorithm is used, requiring the sensitivity fields $\partial {\bm u}({\bm X},\widehat{\bm \lambda}) / \partial \widehat{\lambda}_i$, often obtained numerically through finite differentiation \citep[cf., e.g.,][for more details]{Neggers:2016}.

\subsection{Boundary-Enriched Integrated Digital Image Correlation}
BE-IDIC is an IDIC methodology that considers material parameters $\widehat{\bm{\lambda}}_{\mathrm{mat}}$ and the boundary displacements $\widehat{\bm{\lambda}}_{\mathrm{kin}}$ as unknowns~\citep{rokos}, i.e., 
\begin{equation}
\widehat{\bm{\lambda}}= [\widehat{\bm{\lambda}}_{\mathrm{mat}}^{\sf T},\widehat{\bm{\lambda}}_{\mathrm{kin}}^{\sf T}]^{\sf T}, 
\label{kin-mat}
\end{equation}
where
\begin{align}
    \widehat{\bm \lambda}_{\mathrm{mat}} &= [G_1, K_1, \dots, G_{N_{\mathrm{mat}}}, K_{N_{\mathrm{mat}}}]^{\sf T}\in \mathbb{R}^{2N_{\mathrm{mat}}},\\ 
    \widehat{\bm \lambda}_{\mathrm{kin}} &= \widehat{\sf \bm u}_{\partial \Omega_{\mathrm{mve}}} \in \mathbb{R}^{N_{\mathrm{kin}}},
\end{align}
and where $G_i$, $K_i$ are micromechanical constitutive parameters of the $i$-th phase of a herein assumed hyperelastic material model (cf. Section~\ref{sec:constitutive_model} below), and $\widehat{\sf \bm u}_{\partial \Omega_{\mathrm{mve}}}$ is an array storing displacements of boundary nodes corresponding to chosen discretization of the boundary. The number of kinematic parameters can be very large, depending on the coarseness of this discretization. The cost functional ${\cal R}_{\text{dic}}(\widehat{\bm{\lambda}})$ of Eq.~\eqref{cal_R} is then minimized following the standard procedure, i.e., using the Gauss-Newton method. Note that the number of required sensitivity fields has to be expanded according to the number of kinematic parameters, so the computational cost increases significantly. On the other hand, this method allows for a higher precision, especially when the boundary data is noisy \citep{rokos,rokos:2023}.

\section{Stochastic Approach to Parameter Identification}\label{section:stochastic}

This section explains how a stochastic algorithm can be used to infer material and kinematic parameters of a micromechanical model. To this end, probability densities are first derived in Section~\ref{sec:bayes} using Bayesian theory, which are converted into probability integrals for obtaining relevant statistical quantities in Section~\ref{sec:mha}, computed subsequently through a suitable sampling method (MHA).

The deterministic IDIC approach is a relatively computationally inexpensive method that provides a single value for identified parameters, i.e., it is a deterministic method. On the other hand, obtaining precise enough MVE boundary conditions is a challenge on its own, and the method using GDIC to first identify boundary conditions with subsequent IDIC step \citep[discussed, e.g., in][]{shakoor} hence may be overly sensitive to the accuracy of the boundary conditions. This is mainly because the errors that occur during GDIC phase become fixed and cannot be corrected for in the subsequent IDIC step. 

To account for epistemic uncertainties in the model, such as image noise, Bayesian inference is typically used \citep{oberkampf}, which allows for updating probabilities as more data is gathered (i.e., Bayesian inference). The final answer to a parameter identification problem is then a posterior distribution, as opposed to a single value obtained from deterministic methods, and several uncertainty sources can be straightforwardly incorporated. This is why employing a stochastic method to the problem of parameter identification is potentially beneficial and has been already used in parameter identification of mechanical models, e.g., \citep{ROSIC2013179,BBDP2017, rappel,kucerova,Yue2022}. Although entire probability distribution functions are obtained from stochastic approaches containing much more information, typically certain descriptors are used for direct comparison with deterministic methods, such as mean values or modes \citep[cf.][]{rappel}, which is hereafter adopted as well in the results sections. Additionally, working with prior parameter distributions can regularize ill-posed problems otherwise not solvable by traditional deterministic methods.

\subsection{Parameter Estimation in Bayesian Statistics}\label{sec:bayes}

In the Bayesian statistics, parameter estimation is done by testing numerous hypotheses. 
The data set will then consist of individual tests that result in images $f({\bm X})$ and $g({\bm X})$ $\in [0, 255]$ observed in (virtual) experiments. 
Assuming that the brightness is conserved between two images and omitting the interpolation error, the following relation called brightness conservation, or optical flow equation, e.g., \citep{fernandez2012advanced}, holds
\begin{equation*}
f(\bm{X}) \cong g(\bm{X}+\bm{u}(\bm{X},\widehat{\bm \lambda})),
\end{equation*}
recall also Eq.~(\ref{cal_R}) for the deterministic approach. In practice, we work with the discretized version of $f$ and $g$, 
which we denote as real arrays ${\sf \bm f}$ and ${\sf \bm g}$ of $N_p$ elements, representing individual pixels. Accounting for the measurement error (i.e., image noise) we have
\begin{equation}
{\sf \bm f} + \eta \cong {\sf \bm g}(\widehat{\bm \lambda}) + \zeta,
\label{f-eta}
\end{equation}
where $\eta$ and $\zeta$ are random vectors representing measurement errors, normally distributed with zero mean and variance $\sigma_\eta^2$, i.e., $\eta,\,\zeta \in {\cal N}(0,\sigma_\eta^2)$. 
Eq.~\eqref{f-eta} can be rewritten as
\begin{equation}
\label{noise}
{\sf \bm f} \cong {\sf \bm g}(\widehat{\bm \lambda}) + \xi, \;\;\xi \in {\cal N}(0, 2\sigma_{\eta}^2).
\end{equation}
Because the measurement errors mostly stem from the image noise, we set $\sigma_\eta$ to be 1\% of the dynamic range of ${\sf \bm f}$, ${\sf \bm g}$ $\in [0, 255]$ \citep[see, e.g.,][]{Frank1975}. 

According to the Bayes' theorem, e.g., \citep{sep-bayes-theorem}, the posterior distribution $\pi(\widehat{\bm{\lambda}})$ is proportional to 
\begin{equation}
\label{pi_posterior}
\pi(\widehat{\bm{\lambda}}|{\sf \bm f}, {\sf \bm g}) \propto \pi(\widehat{\bm{\lambda}})\pi({\sf \bm f}, {\sf \bm g}|\widehat{\bm{\lambda}}), 
\end{equation}
where $\pi(\widehat{\bm{\lambda}})$ is the prior distribution of the identified parameters and $\pi({\sf \bm f}, {\sf \bm g}|\widehat{\bm{\lambda}})$ is the likelihood function.
We choose the prior distribution for material parameters as a normal distribution ${\cal N}(\tilde{{\bm \lambda}}_{\mathrm{mat}}^1, \sigma^2_{\mathrm{prior}})$, i.e., 
\begin{equation}
\label{pi_mat}
\pi(\widehat{\bm{\lambda}}_{\mathrm{mat}}) = \frac{1}{\sigma_\mathrm{prior}^{2N_\mathrm{mat}}\sqrt{(2\pi)^{2N_{\text{mat}}}}}\exp\left(-\frac{1}{2}\cdot \frac{||\widehat{\bm \lambda}_{\mathrm{mat}} - \tilde{{\bm \lambda}}_{\mathrm{mat}}^1||_2^2}{\sigma^2_{\mathrm{prior}}}\right),
\end{equation}
where $\tilde{\bm \lambda}_{\mathrm{mat}}^1$ is the initial guess and $\sigma^2_{\mathrm{prior}}$ represents our confidence in it. Note that this distribution can be adopted since the probability $p({\widehat{\bm \lambda}_\mathrm{mat}} < \boldsymbol{0})$ is sufficiently small to not cause any numerical issues. For the kinematic parameters, we choose uniform prior distributions with the mean of $\tilde{{\bm \lambda}}_{\mathrm{kin}}^1$, $ \tilde{\lambda}_{\mathrm{kin},n}^1 \in \mathcal{U}(\tilde{\lambda}^1_{\mathrm{kin},n}-e_\mathrm{bc}, \tilde{\lambda}^1_{\mathrm{kin},n}+e_\mathrm{bc})$, i.e.,
\begin{equation}
\label{pi_kin}
\pi(\widehat{\bm \lambda}_{\mathrm{kin}})=\prod_{n=1}^{N_{\mathrm{kin}}}\pi(\widehat{\lambda}_{\mathrm{kin},n}),
\end{equation}
\begin{equation}
\pi(\widehat{\lambda}_{\mathrm{kin},n})=\begin{cases}
                                \frac{1}{2e_{\mathrm{bc}}} & \text{if} \; |\widehat{ \lambda}_{\mathrm{kin},n}-\tilde{{\lambda}}_{\mathrm{kin},n}^1|\leq e_{\mathrm{bc}},\\
                                0 & \text{otherwise},\\
                                \end{cases}
\end{equation}
where $e_{\mathrm{bc}}$ is the assumed amplitude of the noise in the boundary conditions.
The likelihood  $\pi({\sf \bm f}, {\sf \bm g}|\widehat{\bm{\lambda}})$ is then computed using the probability density function of the noise $\xi$ from Eq.~\eqref{noise} as
\begin{equation}
\pi({\sf \bm f}, {\sf \bm g}|\widehat{\bm{\lambda}}) 
= \frac{1}{(2\sigma_{\eta})^{N_p}\sqrt{\pi^{N_p}} } \exp \left(-\frac{1}{2} \cdot \frac{||{\sf \bm f} - {\sf \bm g}(\widehat{\bm \lambda})||_2^2}{2\sigma_{\eta}^2}\right),
\label{pi_f_g_lambda}
\end{equation}
where $N_p$ is the total number of pixels within ROI. 
Assuming the prior probabilities $\pi(\widehat{\bm \lambda}_{\mathrm{mat}})$ and $\pi(\widehat{\bm \lambda}_{\mathrm{kin}})$ as independent, we can substitute into Eq.~\eqref{pi_posterior} from Eqs.~\eqref{pi_mat}--\eqref{pi_f_g_lambda}, to obtain the approximate posterior probability
\begin{equation}
\pi(\widehat{\bm{\lambda}}|{\sf \bm f}, {\sf \bm g}) = C\, \exp \left[ -\frac{1}{2} \cdot \left( \frac{||\widehat{\bm \lambda}_{\mathrm{mat}} - \tilde{{\bm \lambda}}_{\mathrm{mat}}^1||_2^2}{\sigma^2_{\mathrm{prior}}} + \frac{||{\sf \bm f} - {\sf \bm g}(\widehat{\bm \lambda})||_2^2}{2\sigma_{\eta}^2} \right) \right],
\end{equation}
where $C$ is a suitable normalization constant.

\subsection{Metropolis--Hastings Algorithm}
\label{sec:mha}
The problem \eqref{lambda} has two distinctive sets of parameters: the material parameters $\bm{\lambda}_{\mathrm{mat}}$ and boundary conditions $\bm{\lambda}_{\mathrm{kin}}$, recall Eq.~\eqref{kin-mat}. We are primarily interested in the estimation of $\bm{\lambda}_{\mathrm{mat}}$, although $\bm{\lambda}_{\mathrm{kin}}$ is important for the accuracy. To derive the marginal posterior distribution $\pi(\bm{\lambda}_{\mathrm{mat}}|{\sf \bm f}, {\sf \bm g})$ of the main quantity of interest~$\boldsymbol{\lambda}_\mathrm{mat}$, and not the entire joint distribution $\pi(\bm{\lambda}_{\mathrm{mat}},\bm{\lambda}_{\mathrm{kin}}|{\sf \bm f}, {\sf \bm g})$, the Markov Chain Monte Carlo (MCMC) sampling algorithm is employed \citep{brooks2011handbook}. The marginalization is performed as
\begin{equation}
    \pi(\bm{\lambda}_{\mathrm{mat}}|{\sf \bm f}, {\sf \bm g})=\int_{\mathcal{D}_\mathrm{kin}} \pi(\bm{\lambda}_{\mathrm{mat}},\bm{\lambda}_{\mathrm{kin}}|{\sf \bm f}, {\sf \bm g})\,\mathrm{d}\bm{\lambda}_{\mathrm{kin}},
\label{p}
\end{equation}
where $\mathcal{D}_\mathrm{kin}$ is the domain of integration for $\bm{\lambda}_{\mathrm{kin}}$. To evaluate this integral, the Metropolis--Hastings Algorithm~(MHA) is used.

In MHA, see \citep{Bayes_stats} and Alg.~\ref{alg:mha}, we start with an initial sample $\tilde{\bm{\lambda}}^{i}$ and set it as the current state $\widehat{\bm{\lambda}}^{i}=\tilde{\bm \lambda}^{i}$. Then each new proposal $\tilde{\bm \lambda}^{i+1}$ is generated based on the proposal distribution $q$, which must be symmetric  in the sense that $q(\widehat{\bm \lambda}^{i} \vert \tilde{\bm \lambda}^{i+1}) = q(\tilde{\bm \lambda}^{i+1}\vert \widehat{\bm \lambda}^{i})$ for all $\widehat{\bm \lambda}^{i}, \tilde{\bm \lambda}^{i+1} \in \mathbb{R}^{M=2N_{\mathrm{mat}}+N_{\mathrm{kin}}}$. Typically, $q$ is chosen as Gaussian, i.e., $q \in \mathcal{N}(\widehat{\bm \lambda}^{i}, \sigma_q^2)^{2N_{\mathrm{mat}}+N_{\mathrm{kin}}}$. Its variance $\sigma_q^2$, often referred to as the step size of the random walk, is chosen empirically such that the acceptance rate is around 30\%. The newly proposed state $\tilde{\bm\lambda}^{i+1}$  is accepted with the probability of $\min(1, p)$, where 
\begin{equation}
p = \frac{\pi(\tilde{\bm\lambda}^{i+1}|{\sf \bm f}, {\sf \bm g})}{\pi(\widehat{\bm{\lambda}}^{i}|{\sf \bm f}, {\sf \bm g})}.
\end{equation}
This procedure 
is usually implemented by generating a uniformly distributed random variable $\kappa\in\mathcal{U}(0,1)$. If $\kappa<p$, then $\tilde{\bm \lambda}^{i+1}$ is accepted, i.e., $\widehat{\bm\lambda}^{i+1}=\tilde{\bm\lambda}^{i+1}$, otherwise $\tilde{\bm\lambda}^{i+1}$ is rejected, i.e., $\widehat{\bm\lambda}^{i+1}=\widehat{\bm\lambda}^i$.
Continuing this way, one obtains a sequence $\widehat{\bm{\lambda}}^j$, $j=1,2,\dots,N$. After discarding the first $N_0$ elements (the so-called burn-in), we finally obtain approximation of the marginal (using Eq.~\eqref{p}) or joint posterior distributions of identified parameters $\bm \lambda$ as distributions of $\widehat{\lambda}^j_m$, $j=N_0+1,\dots,N$, for $m=1,\dots,M = 2N_{\mathrm{mat}}+N_{\mathrm{kin}}$.

\begin{algorithm}[H]
\caption{Metropolis--Hastings Algorithm~(MHA) to sample posterior distributions of Eq.~\eqref{p}.}
\label{alg:mha}
\centering
\vspace{-\topsep}
\begin{enumerate}
\item
Draw the initial state $\tilde{\bm\lambda}^{1}\in \mathbb{R}^{M=2N_{\mathrm{mat}}+N_{\mathrm{kin}}}$.\\
Set $\widehat{\bm{\lambda}}^{1}=\tilde{\bm\lambda}^{1}$.
\item
For $i=2,3,\dots,N$ do:
\begin{enumerate}
\item
draw the proposal $\tilde{\bm\lambda}^{i+1}\sim {\mathcal N}(\widehat{\bm{\lambda}}^{i},\sigma_q^2)$
\item
set $p = {\pi(\tilde{\bm\lambda}^{i+1}|{\sf \bm f}, {\sf \bm g})}/{\pi(\widehat{\bm{\lambda}}^{i}|{\sf \bm f}, {\sf \bm g})}$
\item
draw $\kappa \sim {\mathcal U}(0,1)$
\item 
if $\kappa <p$, set $\widehat{\bm{\lambda}}^{i+1}=\tilde{\bm\lambda}^{i+1}$ (i.e., accept $\tilde{\bm\lambda}^{i+1}$), else set $\widehat{\bm{\lambda}}^{i+1}=\widehat{\bm{\lambda}}^{i}$ (i.e., reject $\tilde{\bm\lambda}^{i+1}$)
\end{enumerate}
\item Discard $\widehat{\bm{\lambda}}^{1},\widehat{\bm{\lambda}}^{2},
\dots,\widehat{\bm{\lambda}}^{N_0}$ (burn-in).
\item Estimate statistical parameter characteristics
from $\widehat{\bm{\lambda}}^{N_0+1},
\dots,\widehat{\bm{\lambda}}^{N}$ (mean, standard deviation, etc.).
\end{enumerate}
\end{algorithm}

\section{Underlying Mechanical Model}\label{model}

This section elaborates in detail the mechanical model used in the adopted virtual experiments, describes constitutive model and corresponding micromechanical parameters to be identified, as well as generation of synthetic reference and deformed images for IDIC.

\subsection{Geometry}
To compare performance of the deterministic and stochastic identification procedures, a virtual experiment is performed on a specimen with prescribed material parameters $\bm{\lambda}_{\mathrm{mat,ref}}$. The specimen is assumed to occupy a full-scale domain $\Omega_{\mathrm{dns}}$, having the size of a $20\times20$ units of length, with a heterogeneous structure shown in Fig.~\ref{fig:sketch}. The microstructure consists of randomly distributed non-intersecting stiff circular inclusions with a diameter $d=1$ unit of length, and a surrounding compliant matrix. Although all geometric units and properties are dimensionless, they can be scaled to $\mu$m. 

\subsection{Constitutive Model and Governing Equations}
\label{sec:constitutive_model}
The material of the specimen is assumed to be nonlinear elastic. In particular, a compressible Neo--Hookean hyperelastic material is adopted, see, e.g., \citep{hackett2018hyperelasticity}, specified by the following elastic energy density
\begin{equation}
    \label{el_energy}
    W_{\alpha}(\bm{F}) = \frac{1}{2}G_{\alpha}(\overline{I}_1(\bm{F})-3)+\frac{1}{2}K_{\alpha}(\ln J(\bm{F}))^2,
\end{equation}
where $\bm{F} = [\bm{I}+{\bm \nabla}_0 \bm{u}(\bm{X})]^\mathsf{T}$ is the deformation gradient tensor, ${\bm \nabla}_0 \bm{u}(\bm{X})$ is the gradient of the displacement field, $J(\bm{F}) = \det \bm{F}$, and $\overline{I}_1(\bm{F}) = J^{-\frac{2}{3}}\tr(\bm{C})$ is the first modified invariant of the right Cauchy–Green deformation tensor $\bm{C} = \bm{F}^{\sf{T}}\bm{F}$. Individual materials are distinguished by the subscript $\alpha$, where $\alpha=1$ corresponds to the matrix and $\alpha=2$ to the inclusions. The prescribed values of the material parameters ${\bm \lambda}_{\mathrm{mat, ref}} = [G_1, K_1, G_2, K_2]$ are presented in Tab.~\ref{materials}. Because Dirichlet boundary conditions are applied on the entire boundary of the MVE, $\partial \Omega_{\mathrm{mve}}$, the problem is ill-posed, and only material parameter ratios can be obtained. In the IDIC procedure, one of the material parameters therefore needs to be fixed to an arbitrary value (exact, in our case of virtual experiments) for normalization purposes. The remaining parameters can be identified relative to that reference value. The fixed material parameter can be estimated by other means, i.e., using a force-based mechanical test or reliable experimental sources. Note that such a normalization can be performed in multiple ways, influencing the resulting accuracy of the identification. In contrast to IDIC, MHA does not require such normalization, as will be discussed more extensively below in Section~\ref{section:ill-posed}.

\begin{table}[!htbp]
\centering
\caption{Reference material parameters $\bm{\lambda}_{\mathrm{ref}}$ used in the virtual experiment.}
\begin{tabular}{c | c c} 
 
 \multirow{2}{*}{Physical parameters} & Matrix              & Inclusions \\
									  & $(\alpha = 1)$      & $(\alpha = 2)$ \\[0.5ex] 
 \hline
 Shear modulus, $G_{\alpha}$ & $1$ & $4$ \\ 
 Bulk modulus, $K_{\alpha}$ & $3$ & $12$ \\
 Poisson’s ratio, $\nu_{\alpha}=\frac{3K_{\alpha}-2G_{\alpha}}{2(3K_{\alpha}+G_{\alpha})}$ & $0.35$ & $0.35$ \\
 
\end{tabular}
\label{materials}
\end{table}
The first Piola--Kirchhoff stress tensor ${\bm P}$, introduced in the considered governing Eq.~\eqref{divP} of solid mechanical systems, now attains the form \begin{equation}
    {\bm P}({\bm u}(\bm{X}))=(1-\chi(\bm{X})) \frac{\partial W_1 ({\bm F})}{\partial {\bm F}} [\bm{I}+{\bm \nabla}_0 \bm{u}(\bm{X})]^\mathsf{T}+\chi(\bm{X}) \frac{\partial W_2 ({\bm F})}{\partial {\bm F}} [\bm{I}+{\bm \nabla}_0 \bm{u}(\bm{X})]^\mathsf{T},
\label{st_energy}
\end{equation}
where $\chi(\bm{X})$ is an indicator function of the inclusions, i.e., $\chi=1$ inside all inclusions and $\chi=0$ inside the matrix. Note that the explicit dependence on the identified parameters $\bm \lambda$ has been dropped for brevity.

\subsection{Applied Boundary Conditions}
Let us denote each side of the boundary of $\Omega_{\mathrm{dns}}$, $\partial \Omega_{\mathrm{dns}} =\Gamma$, as $\Gamma_i$, $i=1,\dots,4$, see Fig.~\ref{fig:sketch}. Two virtual mechanical tests are considered, one to introduce tension and another to introduce shear. Both are referred to as Direct Numerical Simulations (DNS) and provide the reference for the mechanical behavior of the system. The displacements prescribed at the specimen’s boundary are 
\begin{equation}
\label{BCs}
\begin{aligned}
\bm{u}(\bm{X})  &= (\overline{\bm{F}}-\bm{I})\bm{X}, && \bm{X} \in  \Gamma_2 \cup \Gamma_4,\\ 
\overline{\bm F} &= \bm{I} + 0.1 \bm{e}_1 \otimes \bm{e}_1, && \text{for tension}, \\
\overline{\bm F} &= \bm{I} + 0.1 \bm{e}_2 \otimes \bm{e}_1, && \text{for shear},
\end{aligned}
\end{equation}
where $\bm{e}_1=[1,0]^{\sf T}$ and $\bm{e}_2=[0,1]^{\sf T}$, while $\Gamma_1$, $\Gamma_3$ are left as free edges. 

\subsection{Numerical Solution and Discretization}
The solution of the mechanical problem in Eq.~\eqref{divP} is determined using the standard non-linear FEM, adopting the Total Lagrangian formulation \citep{borst}. The evolution of the system is solved incrementally, using the standard Newton--Raphson algorithm. 

The displacement field $\widehat{\bm{u}}(\bm{X},\widehat{\bm \lambda})$ is hence discretized as
\begin{equation}
\widehat{\bm u}(\bm{X},\widehat{\bm \lambda}) \approx \sum_{i=1}^{n_{\sf \bm u}/2}N_i(\bm  X)\widehat{\sf \bm u}_i(\widehat{\bm \lambda}),
\end{equation}
where $n_{\sf \bm u}$ denotes the total number of DOFs, $N_i(\bm{X})$ are the standard FE shape functions, and $\widehat{\sf \bm u}_i (\widehat{\bm \lambda}) = [\widehat{\sf u}_1^i (\widehat{\bm \lambda}),\widehat{\sf u}_2^i (\widehat{\bm \lambda})]^{\sf T} \in \mathbb{R}^2$ stores the horizontal and vertical displacements of an $i$-th node associated with a FE mesh. For later reference, these are stored in the array $\widehat{\sf \bm u} = [\widehat{\sf \bm u}^{\sf T}_1,\dots,\widehat{\sf \bm u}^{\sf T}_{n_{\sf \bm u}/2}]^{\sf T} \in \mathbb{R}^{n_{\sf \bm u}}$. 

\begin{figure}[!htbp]
	\centering
	\subfloat[A sketch of the $\Omega_{\mathrm{dns}}$ domain.]{\includegraphics[trim=6.9cm 2.1cm 6.3cm 1.5cm, clip=true, width=0.35\textwidth]{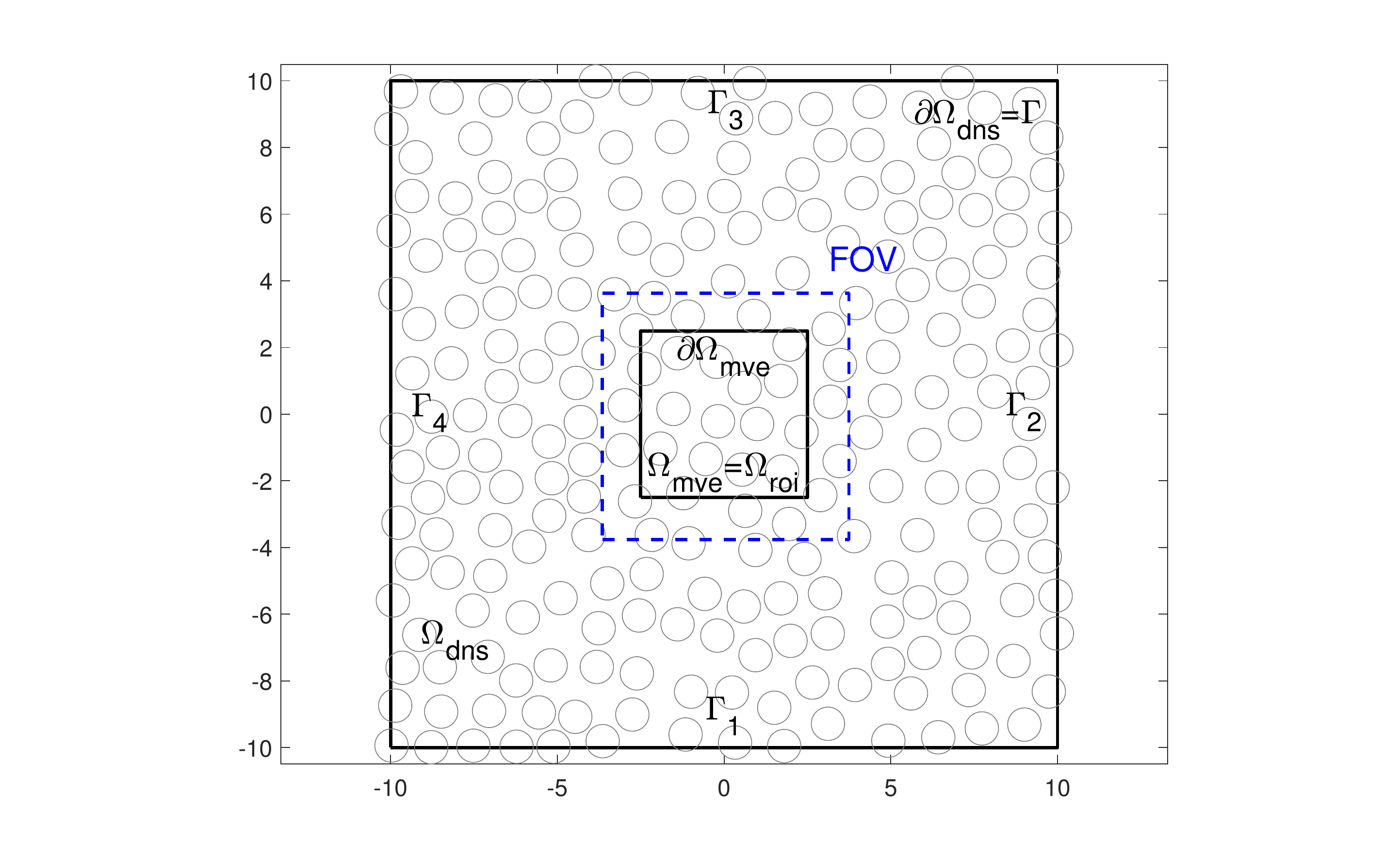}\label{fig:sketch}}
	\hspace{2em}
	\subfloat[Adopted speckle pattern.]{\includegraphics[trim=0cm 0cm -0.6cm 0cm, clip=true, width=0.32\textwidth]{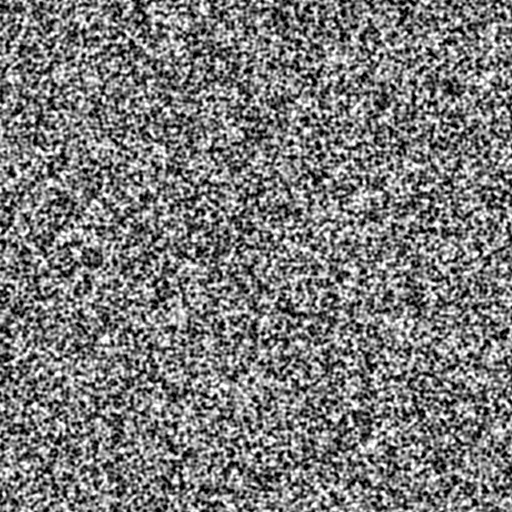}\label{fig:speckle}}\\
	\subfloat[$\Omega_{\mathrm{dns}}$ mesh close-up inside $\Omega_{\mathrm{mve}}$]{\includegraphics[trim=3.5cm 0.5cm 4cm 0.4cm, clip=true, width=0.35\textwidth]{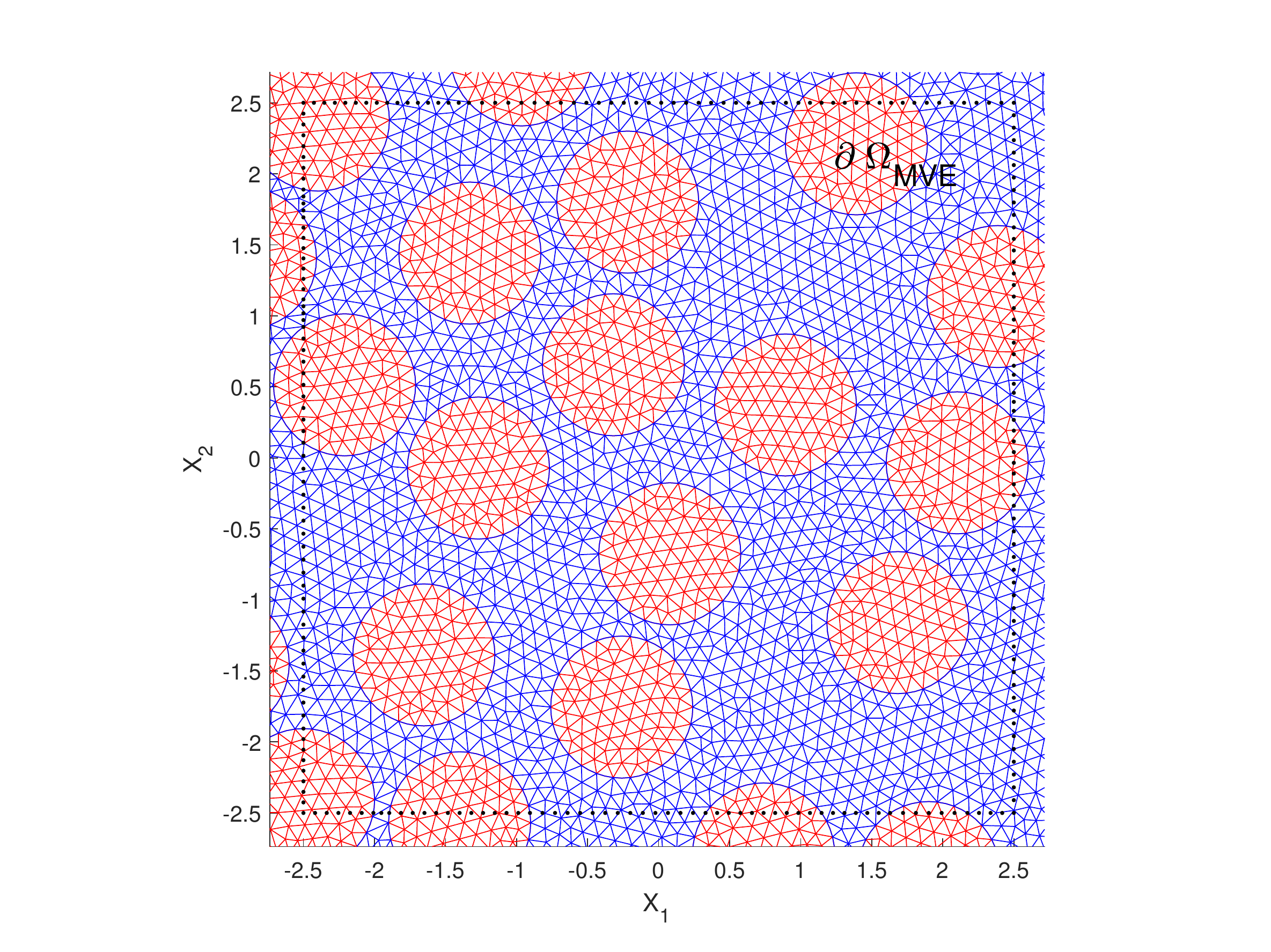}
    \label{fig:dns_mesh_zoom}}
	\hspace{1em}
	\subfloat[$\Omega_{\mathrm{mve}}$ mesh.]{\includegraphics[trim=6cm 0.6cm 6cm 1cm, clip=true, width=0.35\textwidth]{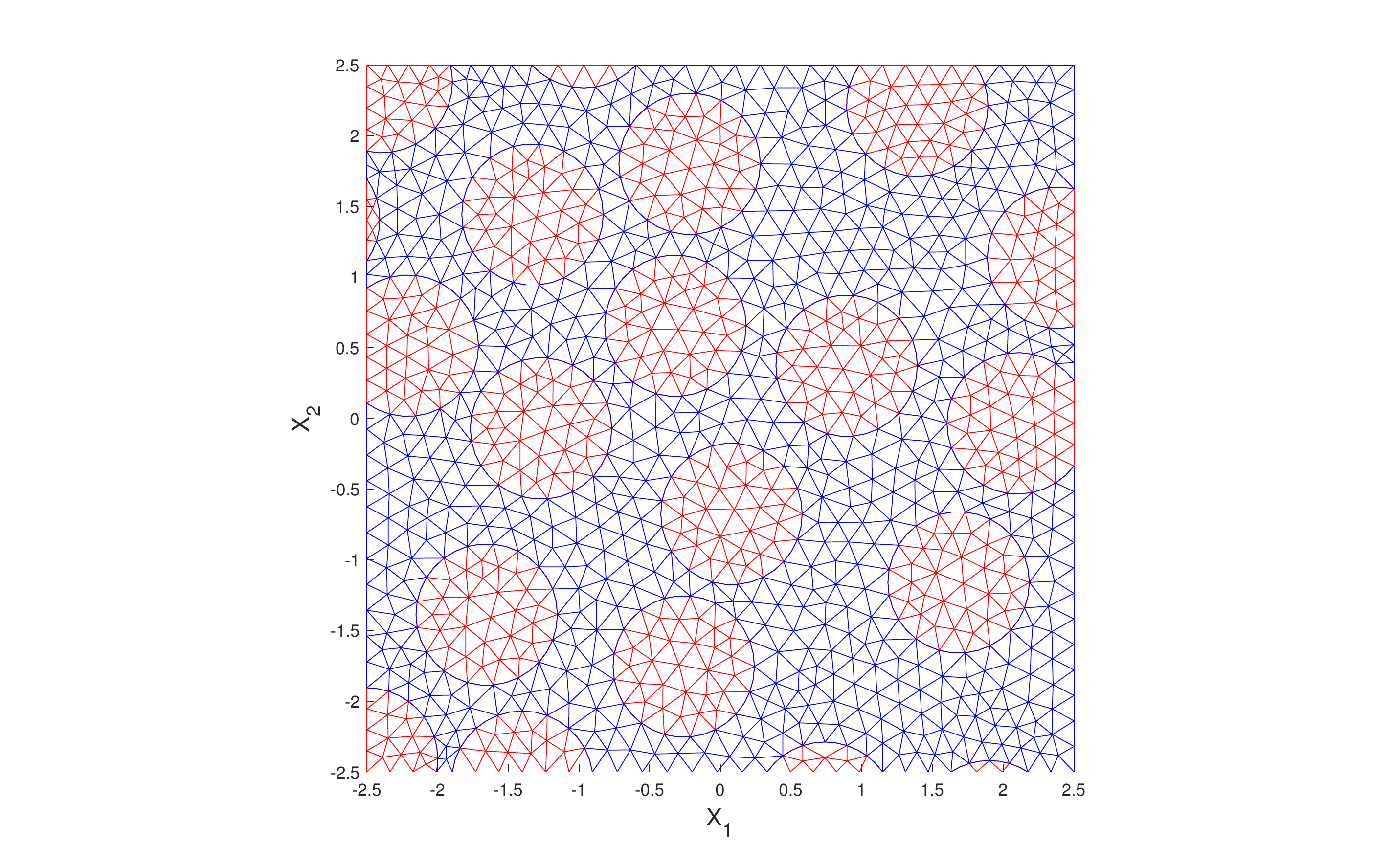}\label{fig:mve_mesh}}
\caption{The macroscopic domain $\Omega_\mathrm{dns}$ consists of stiff circular inclusions of diameter $d=1$ unit of length, embedded within a compliant matrix. (a) A sketch of the specimen's full square domain $\Omega_{\mathrm{dns}}$, microstructural volume element $\Omega_{\mathrm{mve}}$, and the FOV; (b) speckle pattern applied to $\Omega_{\mathrm{mve}}$, and (c) close-up on the MVE domain FEM mesh corresponding to the finely discretized full DNS system; (d) coarse discretization of the MVE model $\Omega_\mathrm{mve}$.}
\label{fig:dns_spec}
\end{figure}

For the FE solution, both $\Omega_{\mathrm{dns}}$ and $\Omega_{\mathrm{mve}}$ domains are discretized with the Gmsh mesh generator \citep{gmsh}, using quadratic isoparametric triangular elements with the three-point Gauss quadrature rule to approximate the integrals appearing in the weak form. All calculations were programmed and performed in MATLAB \citep{matlab}, using an in-house FEM library for hyperelastic materials with computationally heavy parts implemented in external C/C++ mex files for efficiency reasons. For the DNS, the fine mesh shown in Fig.~\ref{fig:dns_mesh_zoom} is used, whereas a coarser MVE triangulation is shown in Fig.~\ref{fig:mve_mesh}. Because the reference Poisson’s ratios for both materials are significantly smaller than 0.5 (recall Tab.~\ref{materials}), and because deformations in the simulations are moderate, incompressibility issues do not occur. 

To track the deformation in a real life experiment, a speckle pattern needs to be applied on the specimen \citep{DIC_guide}. The reference image $\sf f$, representing the applied speckle pattern, has been adopted from \citep[][``medium pattern size"]{bornert}, shown in Fig.~\ref{fig:speckle}. Its resolution is $512\times512$ pixels inside FOV, which corresponds approximately to $340\times340$ pixels inside ROI, when ROI is set equal to MVE. The DNS displacements obtained from all mechanical tests are interpolated from the FE mesh to the regular image mesh by inversion of the elements' isoparametric maps, and the resulting displacement fields are used to map the deformed image $\sf g$ back into the reference configuration. The deformed image is then interpolated at the pixel positions using the bi-cubic polynomial interpolation \citep{interp}. 

\section{Robustness of MHA with Respect to Fixed Errors in Applied Boundary Conditions}\label{section:experiments}
This section presents several experiments that quantify the robustness of the proposed MHA in comparison with the IDIC method \citep{leclerc, ruybalid,  buljac, shakoor} with respect to errors in the applied boundary conditions. The kinematic degrees of freedom $\bm{\lambda}_{\mathrm{kin}}$ are fixed with an error, so the effect of the error on the material parameter identification can be directly observed. The systematic errors, like the smoothing of kinematic fields by the GDIC, are studied first, followed by the effect of the uncorrelated random noise, which is typically observed in the local DIC or the global DIC with a very fine discretization \citep[see][for a similar study]{rokos}.

\subsection{Sensitivity with Respect to Systematic Errors}\label{section:smooth}
\begin{figure}[tbp]
	\centering
\subfloat[Smoothed boundary conditions.]{
    \includegraphics[trim=4.5cm 9cm 5cm 9cm, clip=true, width=0.45\textwidth]{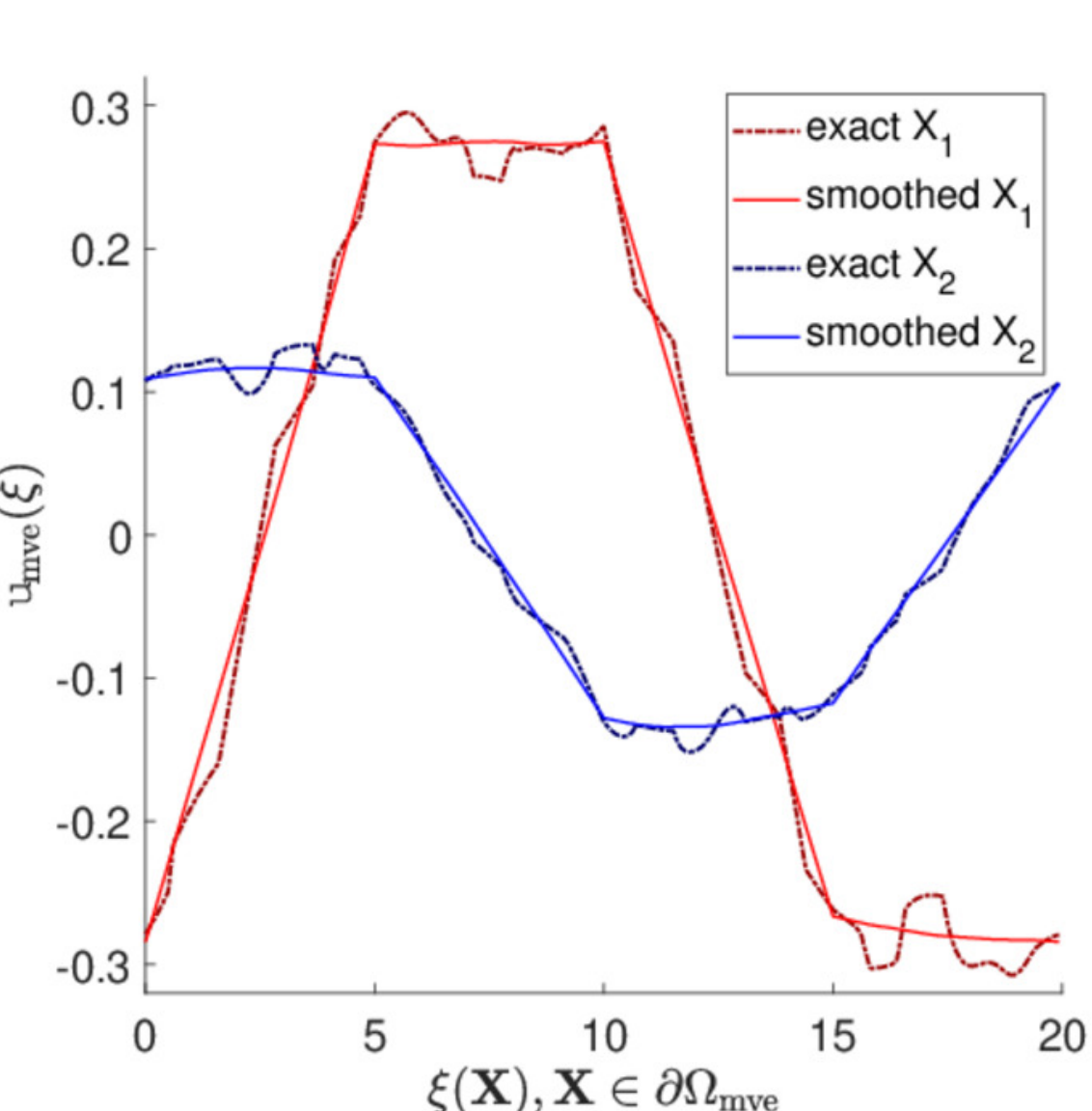}
	\label{fig:smooth_BCs}}
\subfloat[Noisy boundary conditions.]{
    \includegraphics[trim=4.5cm 9cm 5cm 9cm, clip=true, width=0.45\textwidth]{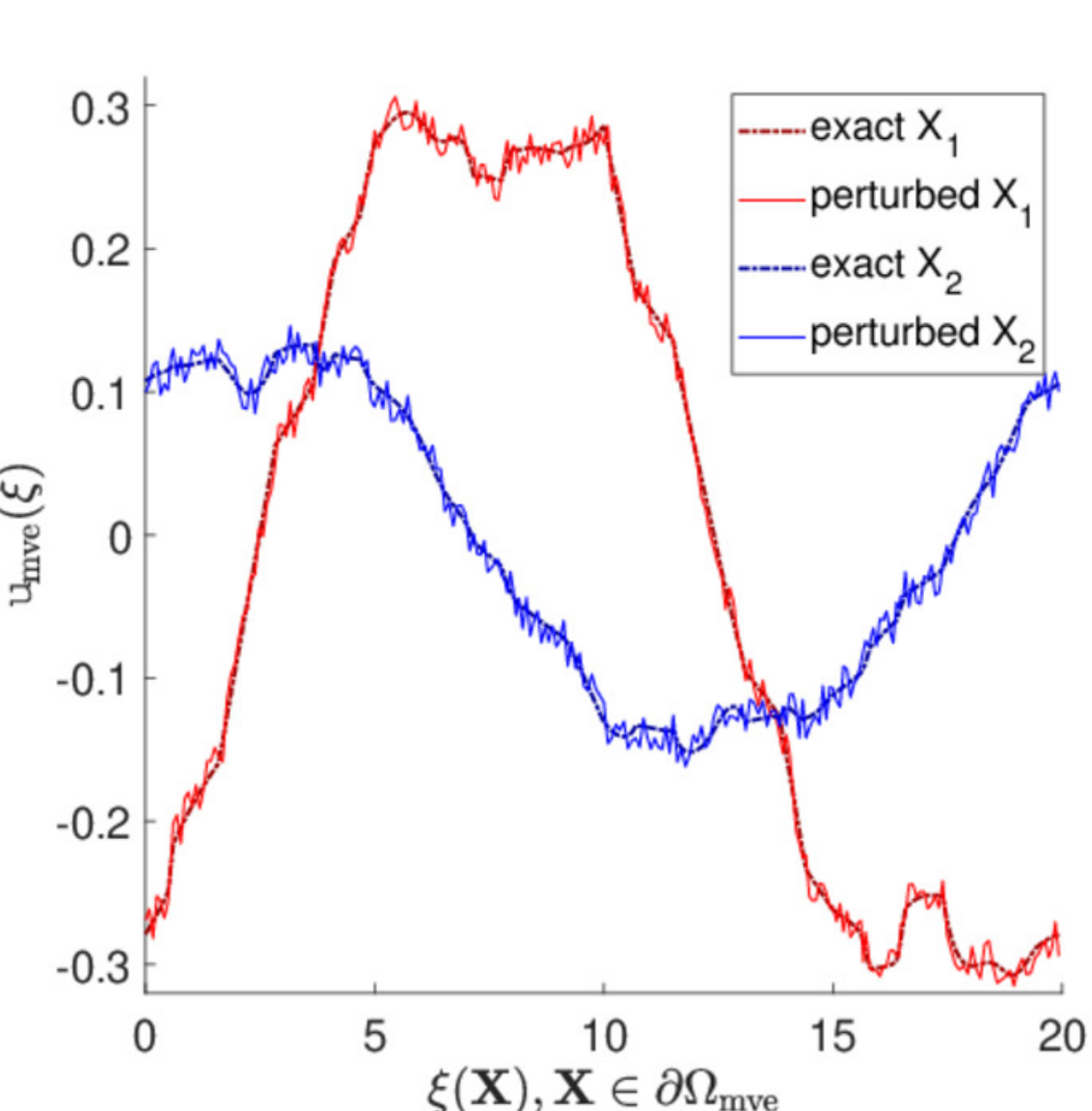}
	\label{fig:BCs_noise}}
\caption{Examples of applied boundary conditions. (a) Smoothed boundary conditions using the pillbox-shaped kernel (red line), $\varepsilon~=~5$, cf. Eq.~\eqref{smooth}, and (b) with superimposed noise according to Eq.~\eqref{noise_imp} with the amplitude $\sigma_{\mathrm{bc}}=0.1$. The scalar coordinate $\xi$ specifies the position at the boundary according to Fig.~\ref{fig:sketch}.}
\end{figure}

The GDIC with a coarse interpolation can have a smoothing effect on the boundary conditions \citep{leclerc2012}. To quantify the effect of smoothed boundary conditions on the identification of the material parameters, the Dirichlet boundary conditions used in the identification are obtained by interpolating the DNS displacements at the nodal positions of the MVE boundary $\partial \Omega_{\mathrm{mve}}$, and the resulting exact displacement field is smoothed with a pillbox-shaped kernel $h_{\varepsilon}$ as, e.g.,~\citep{smith2013digital}
\begin{equation}
 \tilde{\bm u}_{\mathrm{dns}}(\bm{X}) = \int_{\Omega_{\mathrm{dns}}}{\bm u}_{\mathrm{dns}}(\bm{Y})h_{\varepsilon}(\bm{Y}-\bm{X})\,\mathrm{d}\bm{Y},
\label{smooth}
\end{equation}
where $\varepsilon \geq 0$ is a dimension-less diameter  (normalized by the inclusion’s diameter $d=1$). The smoothing effect for the extreme kernel $\varepsilon = 5$ can be observed in Fig.~\ref{fig:smooth_BCs}. The smoothed data are then prescribed as the boundary conditions to the FEM MVE model, i.e.,
\begin{equation}
     {\sf \bm u}_{\mathrm{mve}}(\bm{X}) = \tilde{\sf \bm u}_{\mathrm{dns}}(\bm{X}), \quad \bm{X} \in \partial \Omega_{\mathrm{mve}}.
	 \label{u_eq}	
\end{equation}
\begin{figure}[!h]
	\centering
    \subfloat[Tensile test.]{
        \includegraphics[trim=6.0cm 0.25cm 7.5cm 1.1cm, clip=true, width=0.45\textwidth]{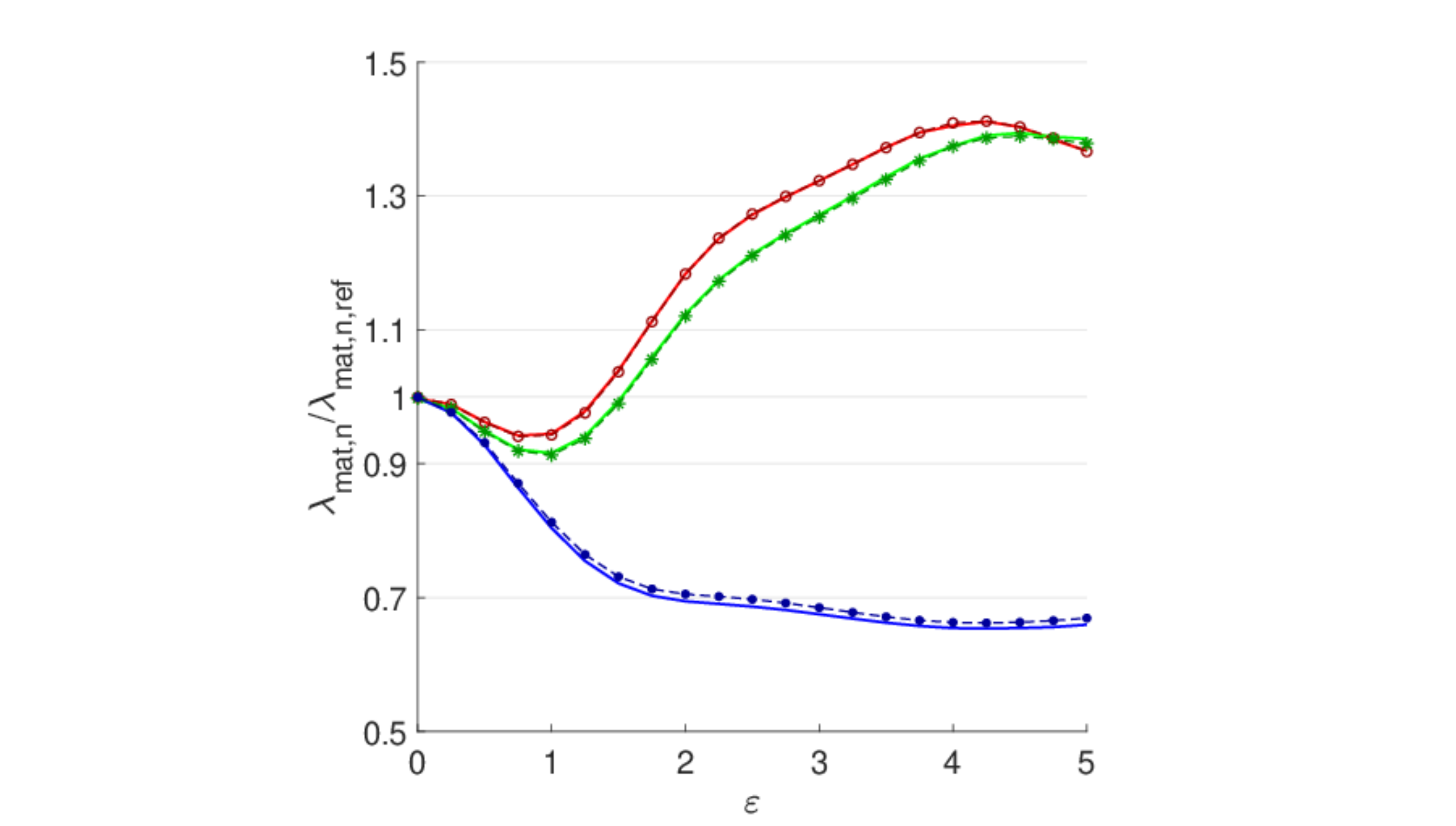}}
    \subfloat[Shear test.]{
        \includegraphics[trim=6.0cm 0.25cm 7.5cm 1.1cm, clip=true, width=0.45\textwidth]{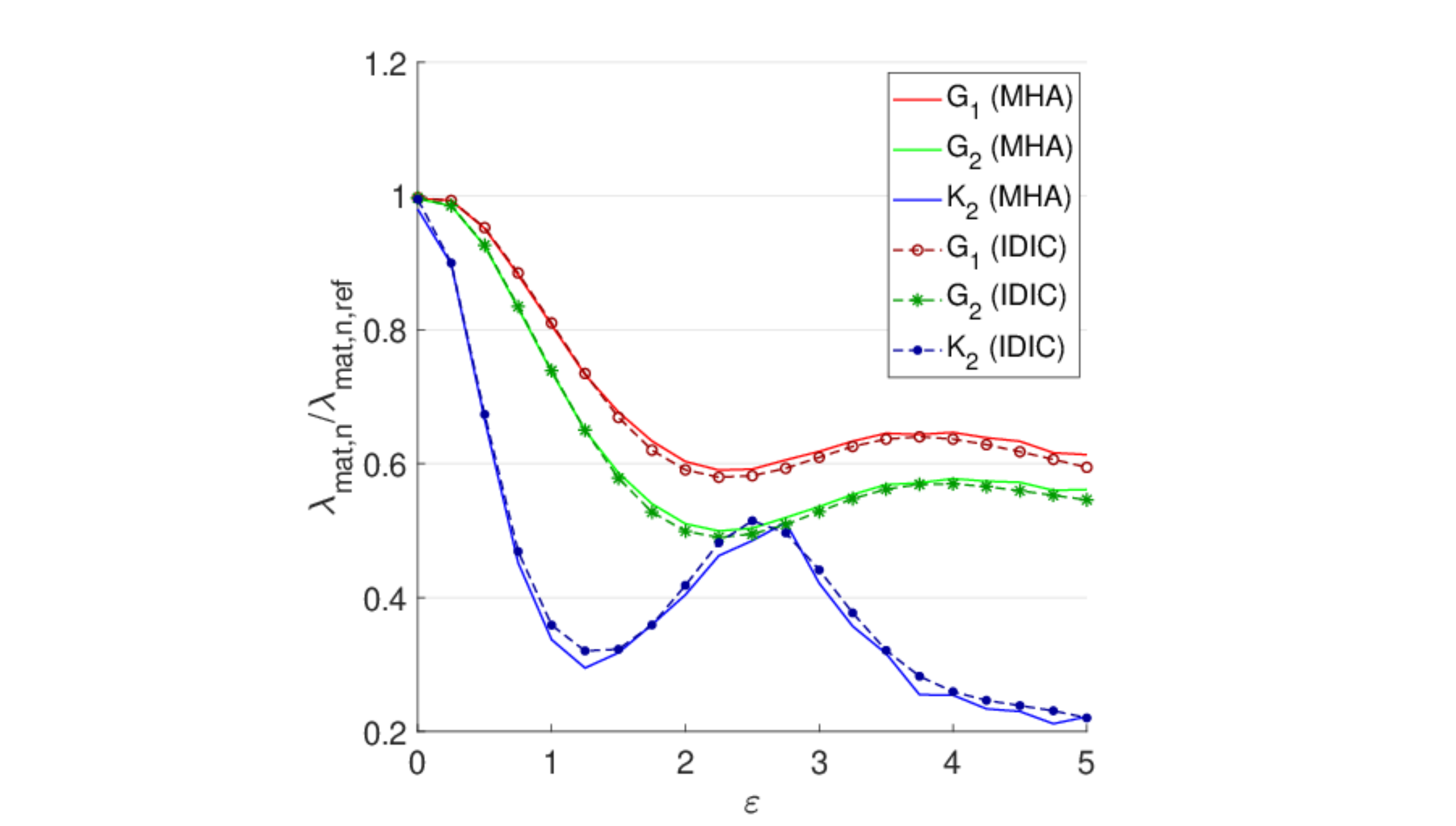}}
\caption{Averaged modes of normalized material parameter distributions $\lambda_{\mathrm{mat,n}}/\lambda_{\mathrm{mat,n,ref}}$ identified by the stochastic MHA method compared against the deterministic IDIC method. Boundary displacements are smoothed using the pillbox-shaped kernel with the parameter $\varepsilon$ and fixed according to Eq.~\eqref{smooth} for (a) tensile and (b) shear tests. Parameters identified by MHA correspond to $N=8\,000$ steps and burn-in $N_0=6\,000$ steps. Parameter $K_1$ is fixed to the reference value to ensure uniqueness.}
\label{fig:smooth}
\end{figure}
In Eq.~\eqref{u_eq}, $\tilde{\sf \bm u}_{\mathrm{dns}}(\bm{X}),\bm{X} \in \partial \Omega_{\mathrm{mve}}$, is a column of displacements of $\tilde{\bm u}_{\mathrm{dns}}$ evaluated at the MVE boundary nodes. For easier implementation, the integral in Eq.~\eqref{smooth} is calculated at discrete pixel positions numerically, while the corresponding displacements at the MVE boundary are linearly interpolated. The bulk modulus of the matrix $K_1$ is chosen as the normalization parameter in both IDIC as well as MHA, so it is fixed at its reference value. Therefore, the dimension of the MHA sampling is reduced only to the remaining material parameters $\bm{\lambda}_{\mathrm{mat}}$, so the sampling is performed in a three-dimensional space. The maximum number of steps for the MHA was set to $N = 8\,000$ with 75\% burn-in. The prior distributions for the material parameters were chosen as normal with mean values set to the initial guess $0.9\lambda_{\mathrm{ref},i}$ and variances $\sigma_{\text{mat,prior}}^2 = 1$. The step size $\sigma_q$ was set as 1\% of $\frac{\lambda_{\mathrm{ref},i}}{\sum_{j=1}^N\lambda_{\mathrm{ref},j}} \sigma_{\text{mat,prior}}$.

Fig.~\ref{fig:smooth} shows obtained modes of the resulting posterior distributions for the tension and shear tests from 50 Monte Carlo~(MC) runs, and compares them to the IDIC method. Here we first observe that the MHA with fixed boundary conditions performs almost identically to the IDIC. A slight difference between the MHA and IDIC can be observed in the bulk modulus of the inclusions $K_2$  in the shear test. The resulting MHA intervals within one standard deviation from the mean are invisible in the scale of the figure, their widths ranging from 5 \textperthousand~to 1\% of the identified value. Next, in accordance with previous observations, see \citep{rokos}, smoothing of the boundary conditions has a considerable negative impact on the identification accuracy. We can therefore conclude that both approaches provide practically identical results for this particular test and show comparable accuracy and robustness.

\subsection{Random Errors in Applied Boundary Conditions}\label{section:noise}
The GDIC can also introduce random errors to the boundary conditions, for too fine meshes. To quantify the effect of random noise in boundary conditions, the following test is performed. Uncorrelated random noise is superimposed on the exact boundary displacement, i.e., 
\begin{equation}
{\sf \bm u}_{\mathrm{mve}}(\bm{X}) = {\sf \bm u}_{\mathrm{dns}}(\bm X) + \sigma_{\mathrm{bc}} \, \underset{\bm{Y} \in \partial\Omega_\mathrm{dns}}{\max} (||{\sf \bm u}_{\mathrm{dns}}(\bm{Y}) ||_2)\,{\sf \bm u}_\mathrm{rnd}, \quad {\bm X} \in \partial \Omega_{\mathrm{mve}},   
\label{noise_imp}
\end{equation}
where ${\sf \bm u}_{\mathrm{mve}}(\bm{X})$ is a vector that stores nodal displacements of the MVE boundary nodes, ${\sf \bm u}_{\mathrm{dns}}(\bm{Y})$ does the same for the DNS boundary nodes, ${\sf \bm u}_\mathrm{rnd}$ is the corresponding column of independent and identically distributed random variables with a uniform distribution over $[-0.5,0.5]$, i.e., from $\mathcal{U}(-0.5,0.5)$, and $\sigma_{\mathrm{bc}} \in [0, 0.1]$ is the standard deviation of the prescribed random noise.

All three  material parameters are identified again for both the tensile and shear tests, using 50 MC realizations of the noise for each value of $\sigma_{\mathrm{bc}}$. An example of boundary data with $\sigma_{\mathrm{bc}}=0.1$ is shown in Fig.~\ref{fig:BCs_noise}. The number of steps for the MHA is set to $N = 8\,000$ with 75\% burn-in, with the rest of hyper-parameters and prior parameter distributions set the same as in the previous example. The boundary conditions are kept fixed with $K_1$ used as the normalization parameter again.

\begin{figure}[!htbp]
\centering
	\subfloat[Tensile test.]{
	\includegraphics[trim=6.5cm 0cm 7.5cm 0.5cm, clip=true, width=0.45\textwidth]{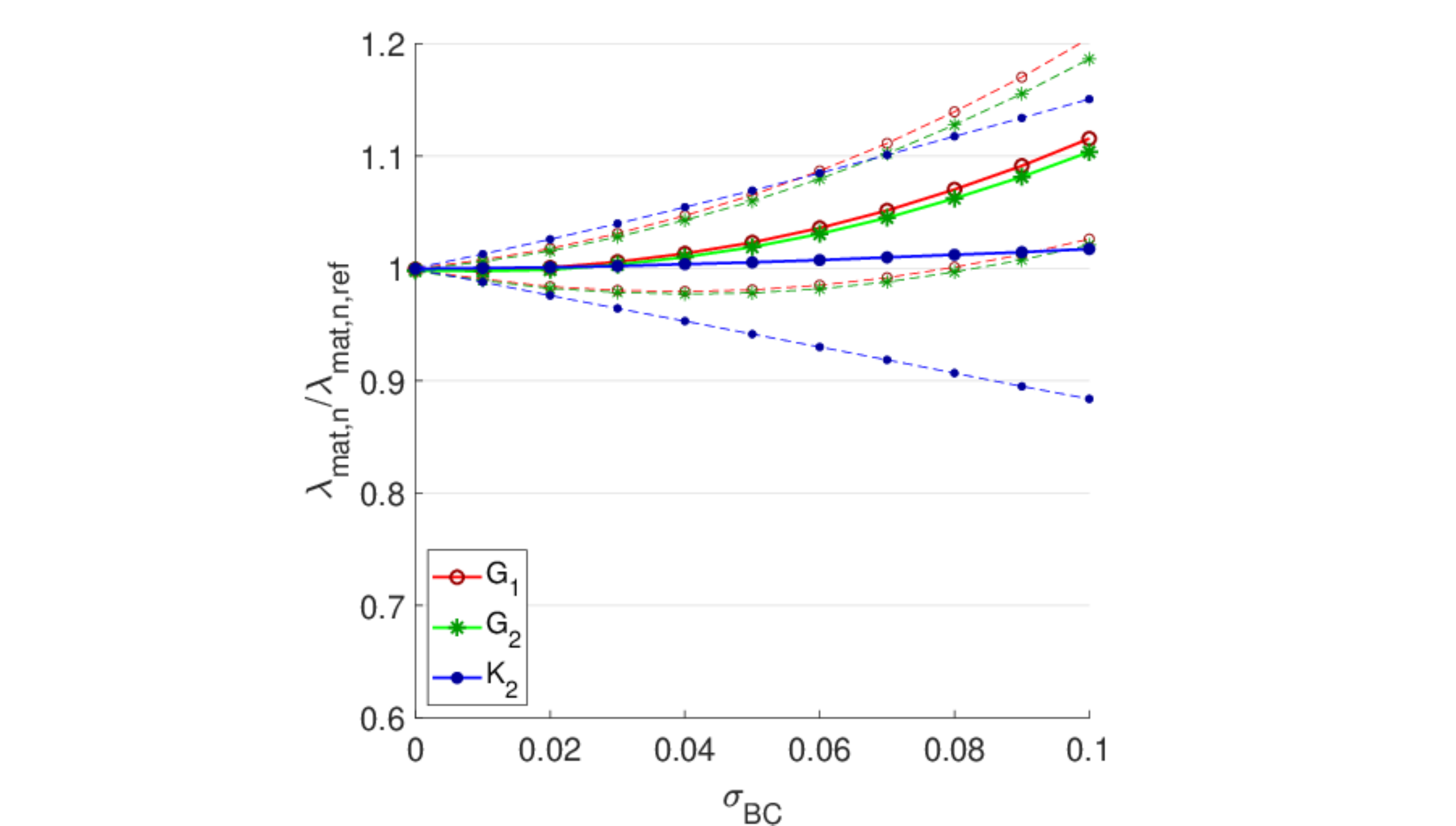}}
	\subfloat[Shear test.]{
	\includegraphics[trim=6.5cm 0cm 7.5cm 0.5cm, clip=true, width=0.45\textwidth]{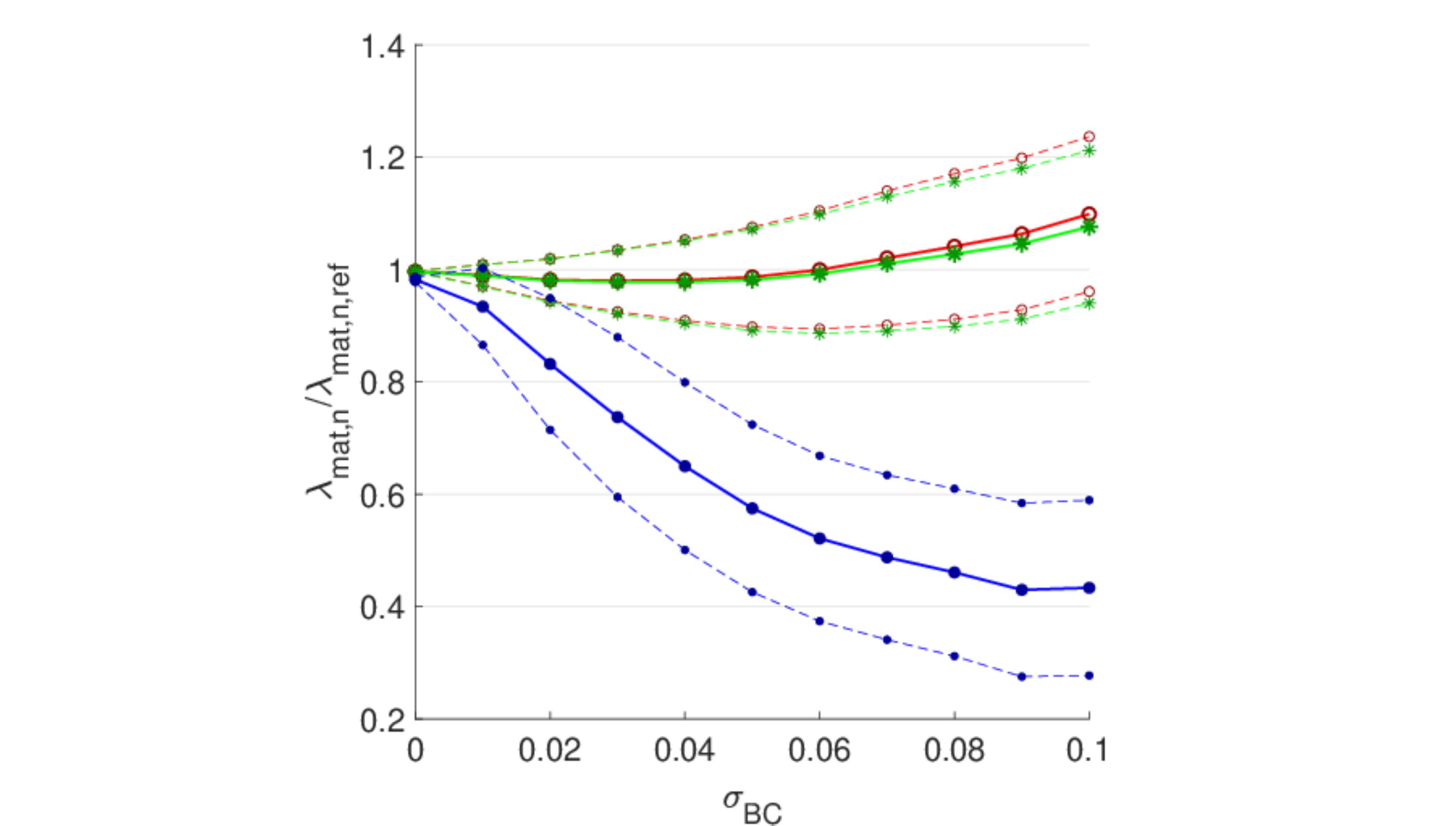}}
\caption{Averaged modes of normalized material parameter distributions $\lambda_{\mathrm{mat,n}}/\lambda_{\mathrm{mat,n,ref}}$ identified by MHA for the fixed random noise in boundary conditions, cf Eq.~\eqref{noise_imp}, with a noise of standard deviation $\sigma_{\mathrm{bc}}$ for (a) tensile and (b) shear tests. The means across all sampled values of all iterations (after burn-in) are plotted with solid lines and are complemented with $\pm$ standard deviations (dashed lines) with $N=8\,000$ and burn-in $N_0=6\,000$ steps. Parameter $K_1$ is fixed to the reference value to ensure uniqueness.}
\label{fig:ns_all}
\end{figure}

\begin{figure}[!htbp]
\centering
\subfloat[Tension, $G_1$.]{
		\includegraphics[trim=6cm 0cm 8cm 0cm, clip=true, width=0.32\textwidth]{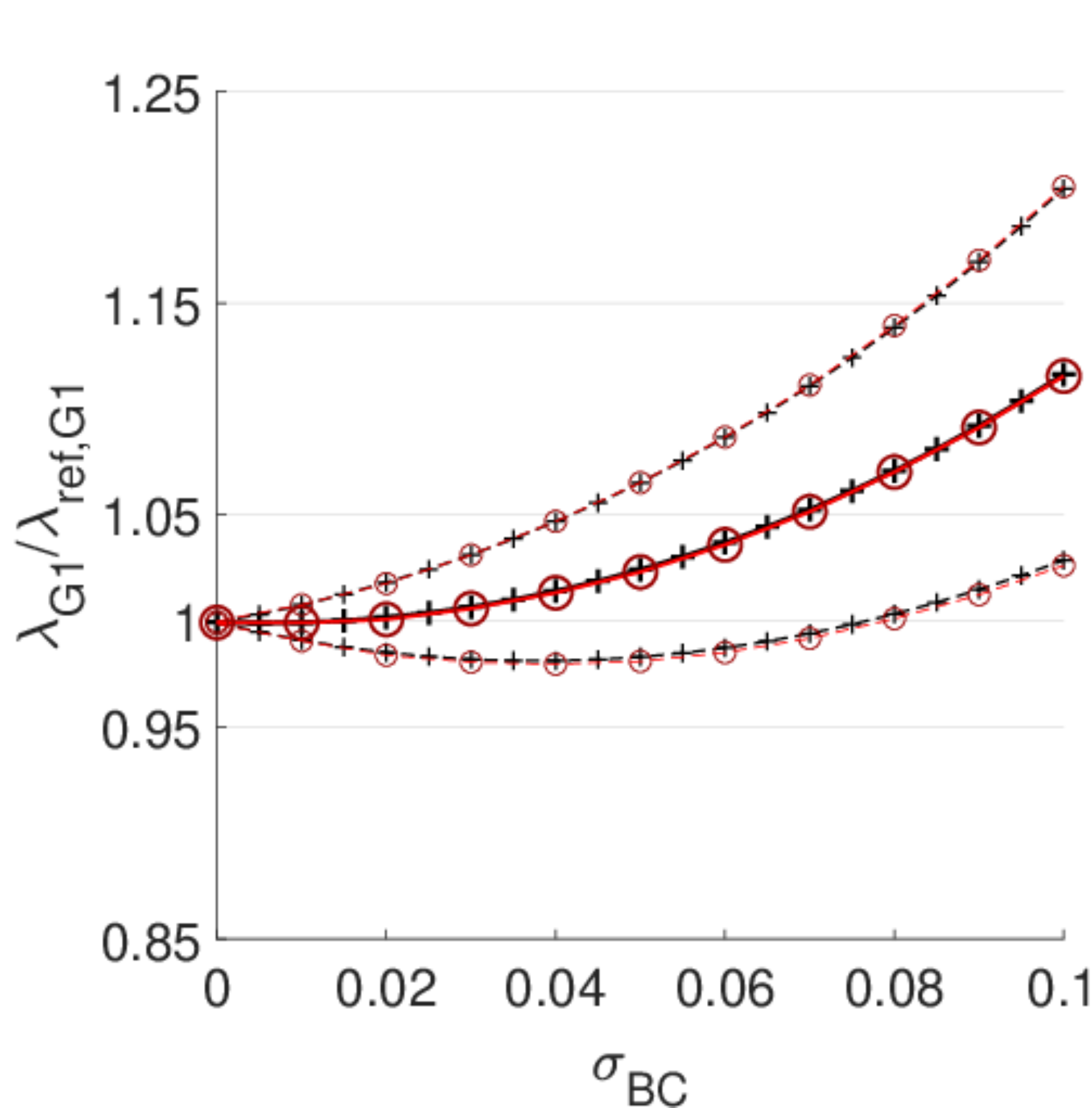}
		\label{fig:G1}}
\subfloat[Tension, $G_2$.]{
		\includegraphics[trim=6cm 0cm 8cm 0cm, clip=true, width=0.32\textwidth]{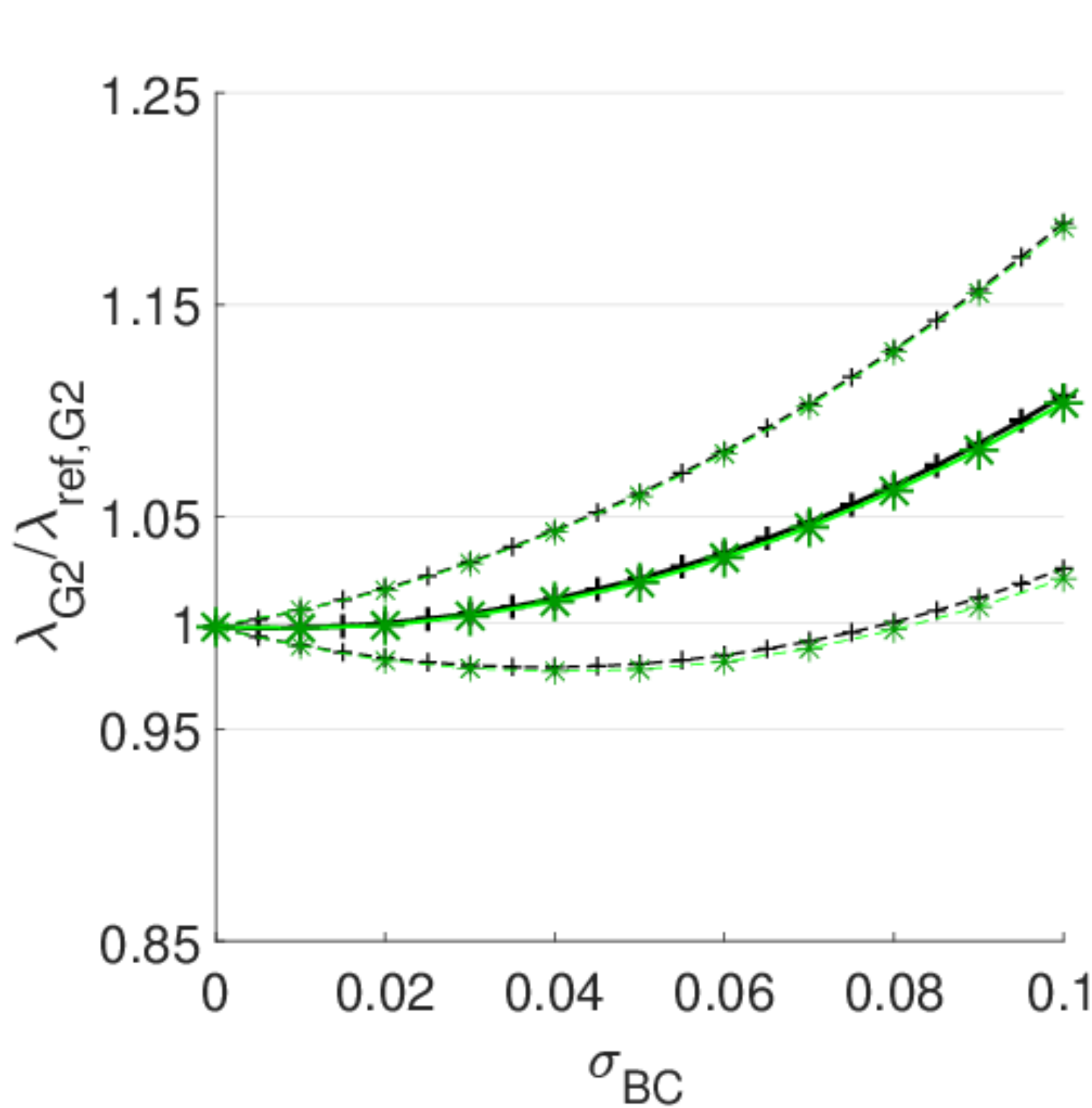}
		\label{fig:G2}}
\subfloat[Tension, $K_2$.]{
		\includegraphics[trim=6cm 0cm 8cm 0cm, clip=true, width=0.32\textwidth]{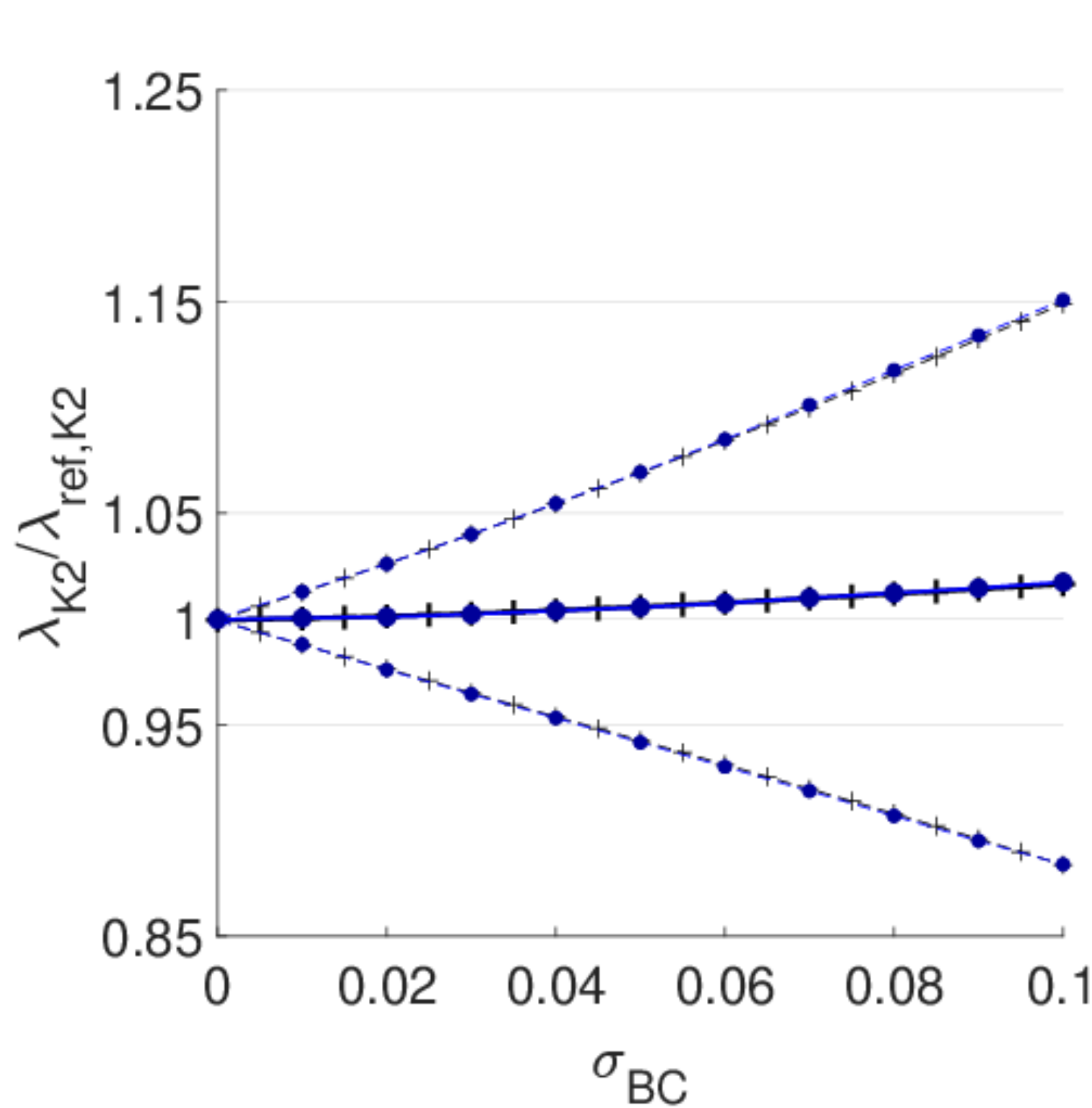}
		\label{fig:K2}}\\
\subfloat[Shear, $G_1$.]{
		\includegraphics[trim=6cm 0cm 8cm 0cm, clip=true, width=0.32\textwidth]{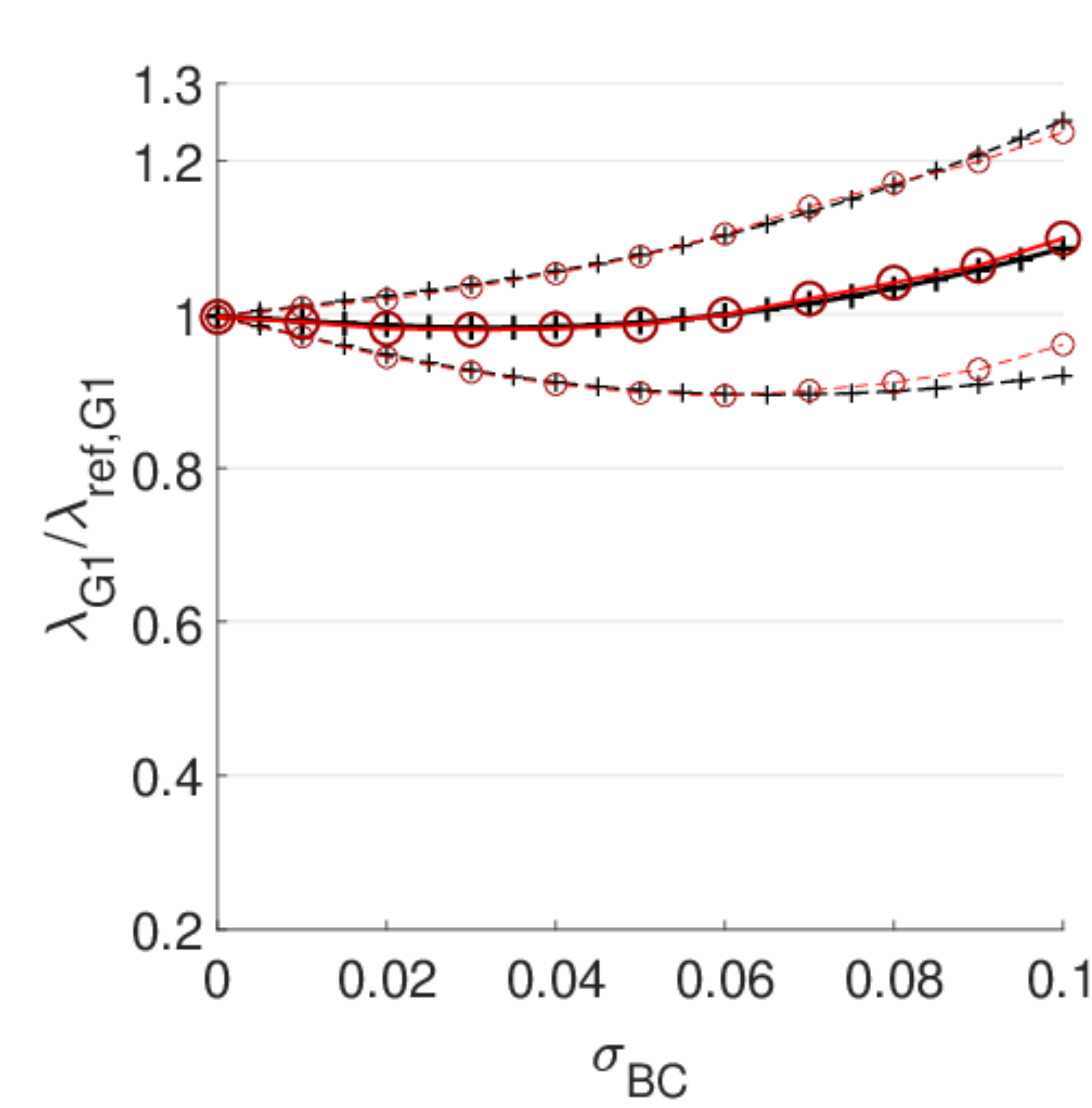}
		\label{fig:G1_shear}}
\subfloat[Shear, $G_2$.]{
		\includegraphics[trim=6cm 0cm 8cm 0cm, clip=true, width=0.32\textwidth]{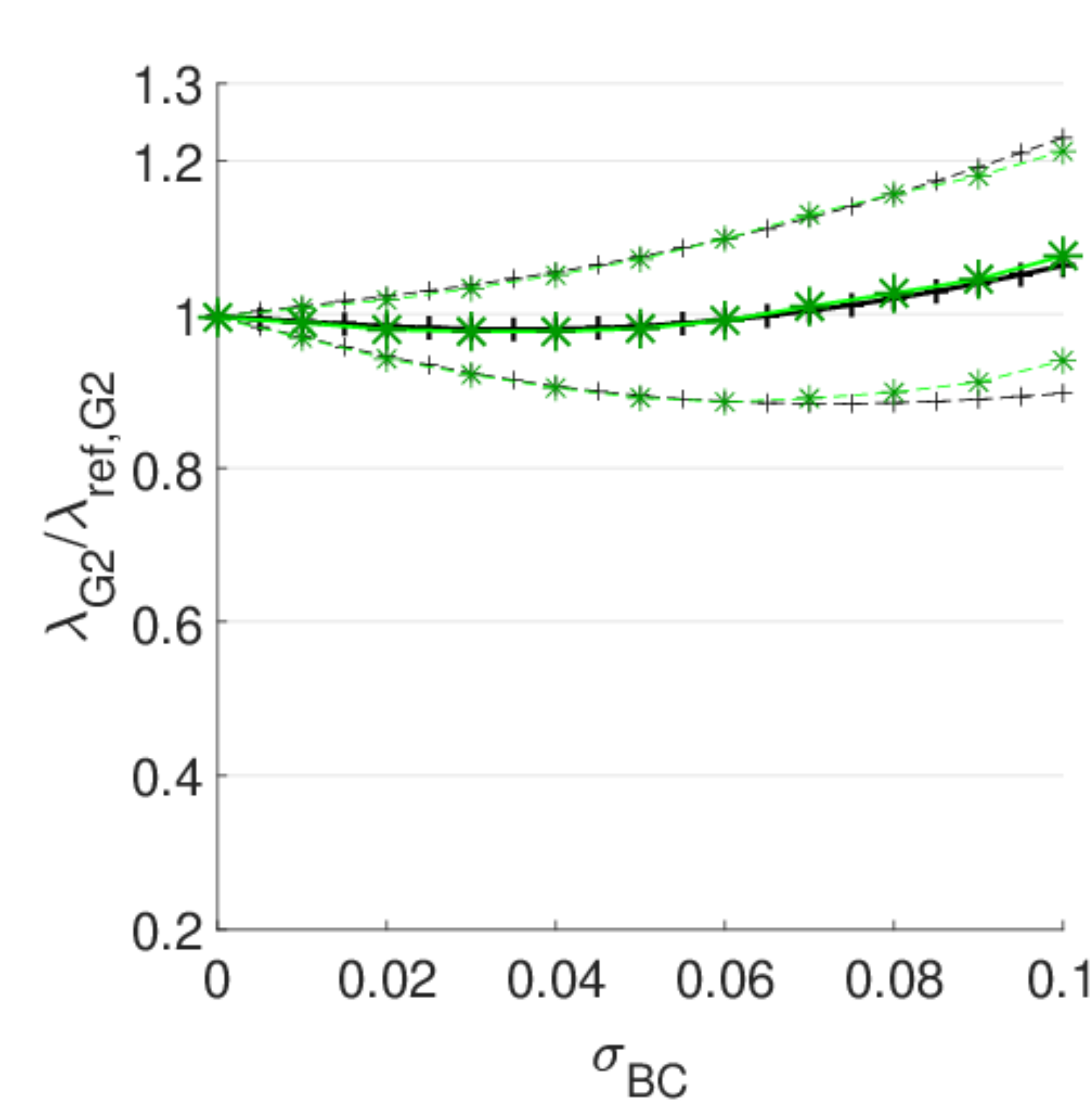}
		\label{fig:G2_shear}}
\subfloat[Shear, $K_2$.]{
		\includegraphics[trim=6cm 0cm 8cm 0cm, clip=true, width=0.32\textwidth]{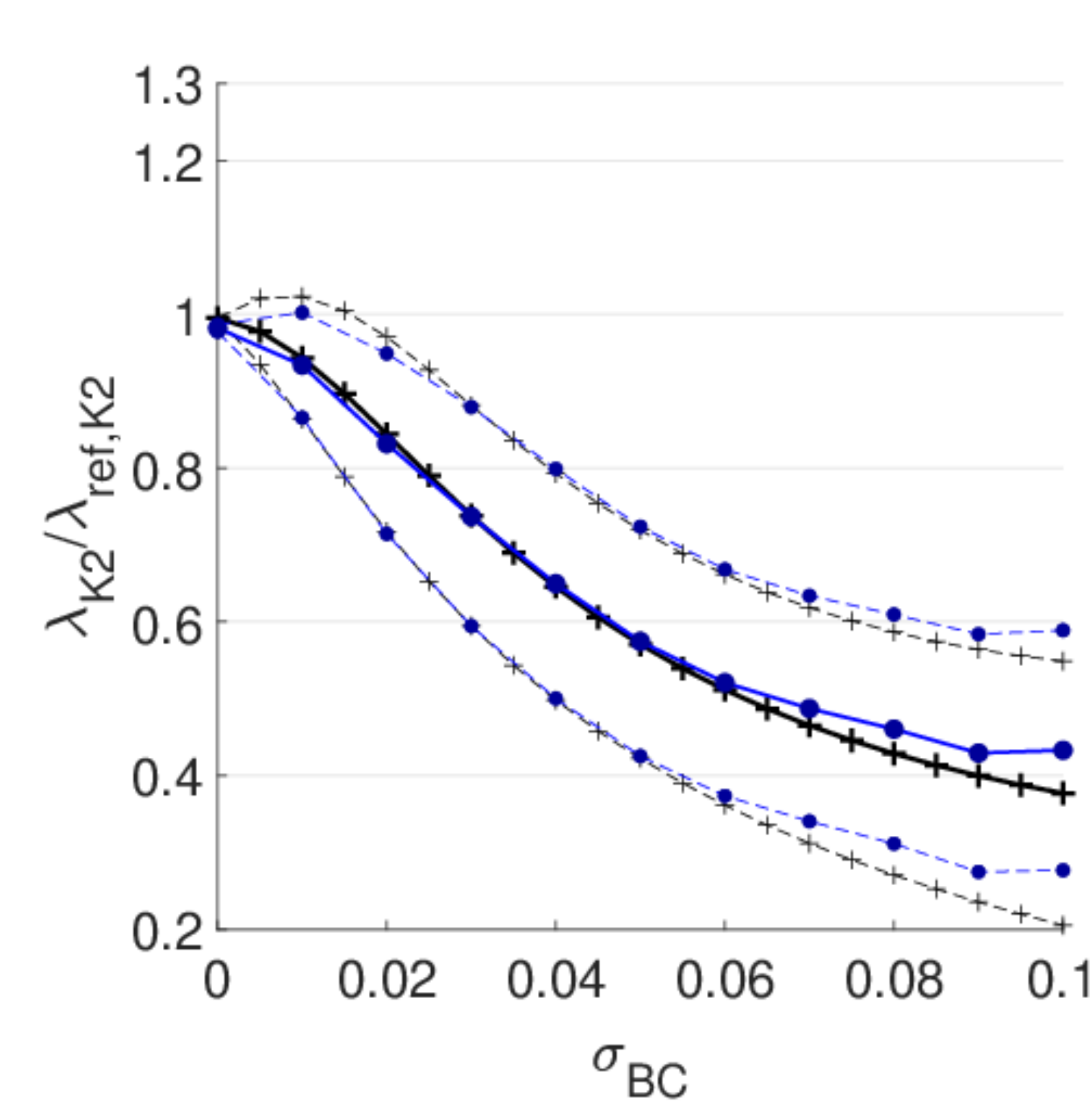}
		\label{fig:K2_shear}}
\caption{Comparison of the averaged modes of normalized distributions of identified (a, d) matrix shear modulus $G_1$, (b, e) fibre shear modulus $G_2$, and (c, f) fiber bulk modulus $K_2$ obtained from the tensile (top) and shear (bottom) tests under fixed boundary conditions with applied noise according to Eq.~\eqref{noise_imp} of standard deviation $\sigma_{\mathrm{bc}}$ for the IDIC (black) and the MHA (color) methods. The means across all sampled values of all iterations (after burn-in) are plotted with solid lines and are complemented with $\pm$ standard deviations (dashed lines). The results correspond to all 50 MC realizations with $N=8\,000$ and burn-in $N_0=6\,000$ steps. Parameter $K_1$ is fixed  to the reference value to ensure uniqueness.}
\label{fig:ns_ind}
\end{figure}

For each of the resulting posterior distributions, the mode was again taken as the measure of the identified parameters and then averaged across all 50 realizations with the same noise amplitude. The derived mean values for all the material parameters are shown in Fig.~\ref{fig:ns_all} with thick lines, while the standard deviations are superimposed over the mean values (dashed lines). Similarly to the previous example the confidence intervals for each iteration are negligible.

Fig.~\ref{fig:ns_ind} shows that the MHA (color lines) and the IDIC (black lines) again deliver almost identical results. The slight deviation for the values of $\sigma_{\mathrm{bc}} > 0.05$ in the shear test can be explained by the insufficient convergence of the MHA for the given number of steps.

We conclude that the MHA with fixed errors in boundary conditions has the same robustness as the IDIC method, and that the error in boundary conditions has a significant effect on the values of the identified material parameters, especially in the shear test.

\section{MHA with Relaxed Boundary Conditions}\label{section:relaxed_BC}
One of the benefits of the MHA is that the number of sampled parameters can be increased to include boundary conditions with no direct additional computational effort, as compared to the case with fixed boundary conditions. A straightforward way to implement new DOFs is to set the boundary condition parameters $\bm{\lambda}_{\mathrm{kin}}$ as displacements in the FE nodes along $\partial \Omega_{\mathrm{mve}}$. For the configuration described in Section \ref{model} (and used for all numerical experiments), this would mean adding extra $2 \times 244$ parameters, considering that 244 FE nodes (on a randomly generated domain) store displacements along the horizontal and vertical directions. We expect this approach to provide more accurate estimates for the material parameters, since it also minimizes the error in the boundary conditions. On the other hand, this approach is generally expected to introduce the so-called ``curse of dimensionality'' \citep{Au,zuev}. Even though calculating the solution of the mechanical system for each new sample requires the same computational effort for any number of the employed DOFs, the algorithm might require more steps to find the high probability region, as the number of sampled parameters grows. One way to address this problem is to reduce MHA sampling to manageable dimensions by approximating boundary conditions with some basis functions, e.g., using Fourier transform or Karhunen--Loève expansion, both widely used in signal enhancement \citep{hermus}. However, the MHA sampling in the resulting spaces is not uniform and thus requires substantial additional tuning. That is why, for the purpose of this paper, a simple reduction of nodes in the FE basis was ultimately chosen, where only some of the existing DOFs on the boundary $\partial\Omega_{\mathrm{mve}}$ were used to reduce the parameters' dimension. The rest of the parameters are dependent on those parameters (by interpolation). The initial MVE boundary mesh was constructed by assigning nodes at equal intervals with displacements interpolated from the DNS, and every $n$-th node is then used for the boundary approximation. The displacement error of the boundary condition is quantified as
\begin{equation}
  \epsilon^{\mathrm{bc}}_{\mathrm{rel}} = \frac{||{\sf \bm u}_{\mathrm{mve}}^t(\bm{X})-{\sf \bm u}_{\mathrm{dns}}(\bm{X})||_2}{||{\sf \bm u}_{\mathrm{dns}}(\bm{X})||_2}, \quad \bm{X} \in \partial \Omega_{\mathrm{mve}},
  \label{error-fe}
\end{equation}
where ${\sf \bm u}_{\mathrm{mve}}^t(\bm{X})$ is the boundary approximation with $2\times t$ DOFs, and ${\sf \bm u}_{\mathrm{dns}}(\bm{X})$ is the exact boundary displacement obtained by the DNS. 

The boundary condition parameter $\bm{\lambda}_{\mathrm{kin}}$ can be accordingly relaxed in the MHA sampling. The number of total employed kinematic DOFs was $2\times61$ (25\%), $2\times122$ (50\%), $2\times183$ (75\%), and $2\times244$ (100\%). To examine the MHA's robustness with respect to the noise in the initial applied boundary conditions, the initial guess of the kinematic parameter $\tilde{\bm{\lambda}}_{\mathrm{kin}}^1$ was assumed with an increasing noise amplitude $\tilde{\sigma}_{\mathrm{bc}}$. The starting point for the material parameter $\bm{\lambda}_{\mathrm{mat}}^1$ was assumed as $\lambda_{\mathrm{ref},i}$, $i=1,\dots,4$, so the influence of the relaxed boundary conditions can be assessed.  The bulk modulus of the matrix $K_1$ was used as the normalization parameter. The prior distributions for the material parameters were again chosen as normal, with the mean values set to their reference value as the initial guess $\lambda_{\mathrm{ref},i}$, with the variances $\sigma_{\text{mat,prior}}^2 = 1$. The material parameter step size $\sigma_{q, i}$ was set as 5\textperthousand~of $\frac{\lambda_{\mathrm{ref},i}}{\sum_{j=1}^N\lambda_{\mathrm{ref},j}} \sigma_{\text{mat,prior}}$. The prior distributions for the kinematic parameters were set as uniform according to Eq.~(\ref{pi_kin}), where $e_{\mathrm{bc}} = 0.1 \underset{\bm{Y} \in \partial\Omega_{\mathrm{mve}}}{\max} (||{\sf \bm u}_{\mathrm{dns}}(\bm{Y}) ||_2)$ is the maximum possible error of the initial guess for the boundary node coordinate. The kinematic parameter step size was then chosen relative to the material parameter step size as 0.4\% of $\sum_{i=1}^N\sigma_{q, i}$. The total number of steps was set as $N = 24\,000$ with 92\% burn-in (chosen during post-processing) in all experiments for practicality. 

\subsection{Tensile Test}
\begin{figure}[!htbp]
\centering
\subfloat[$t=2\times244$ DOFs (100\% of total).]{
	\includegraphics[trim=6.5cm 0cm 7.5cm 1cm, clip=true, width=0.4\textwidth]{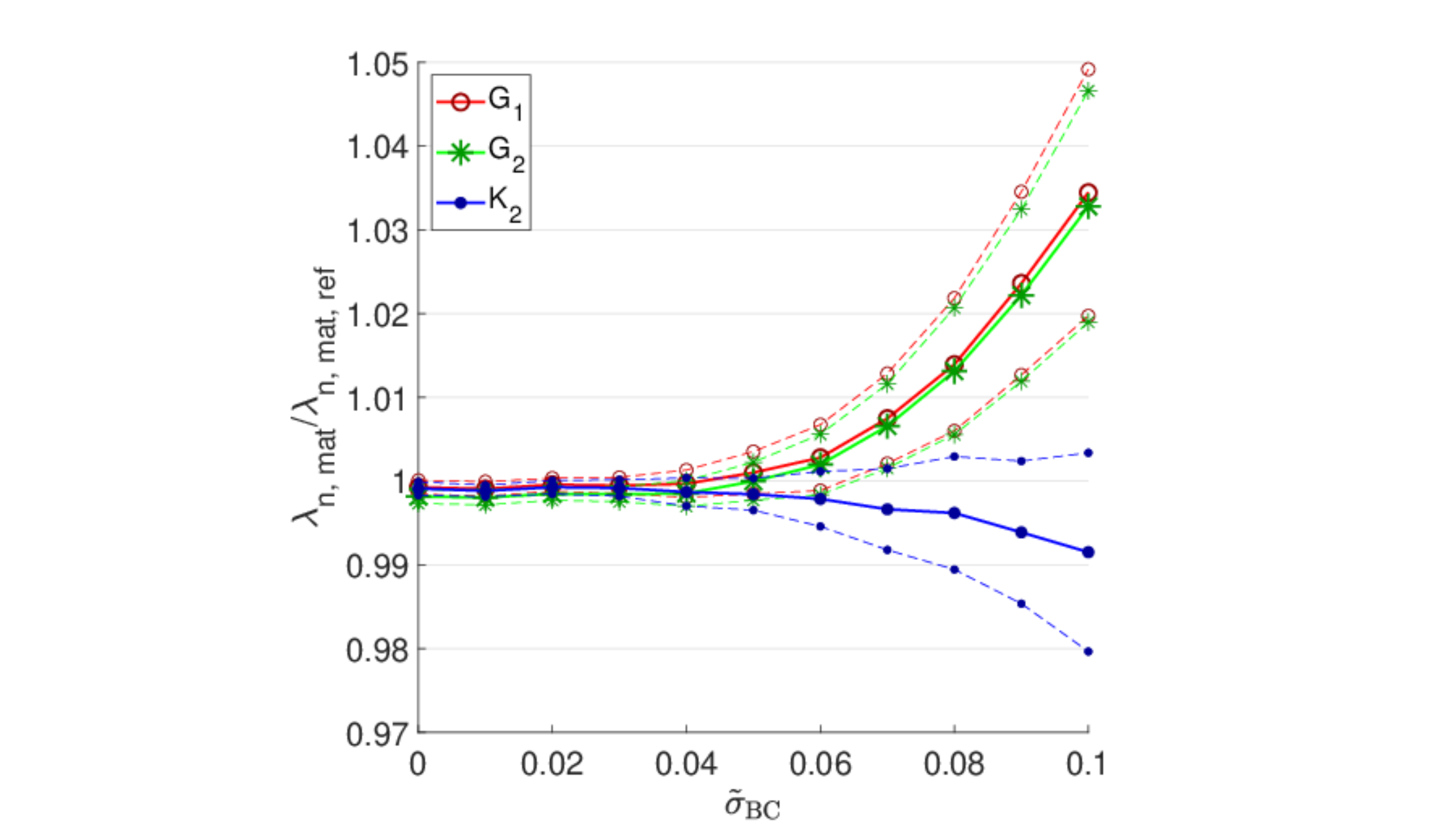}}
	\hspace{0.5em}
\subfloat[$t=2\times183$ DOFs (75\% of total).]{
	\includegraphics[trim=6.5cm 0cm 7.5cm 1cm, clip=true, width=0.4\textwidth]{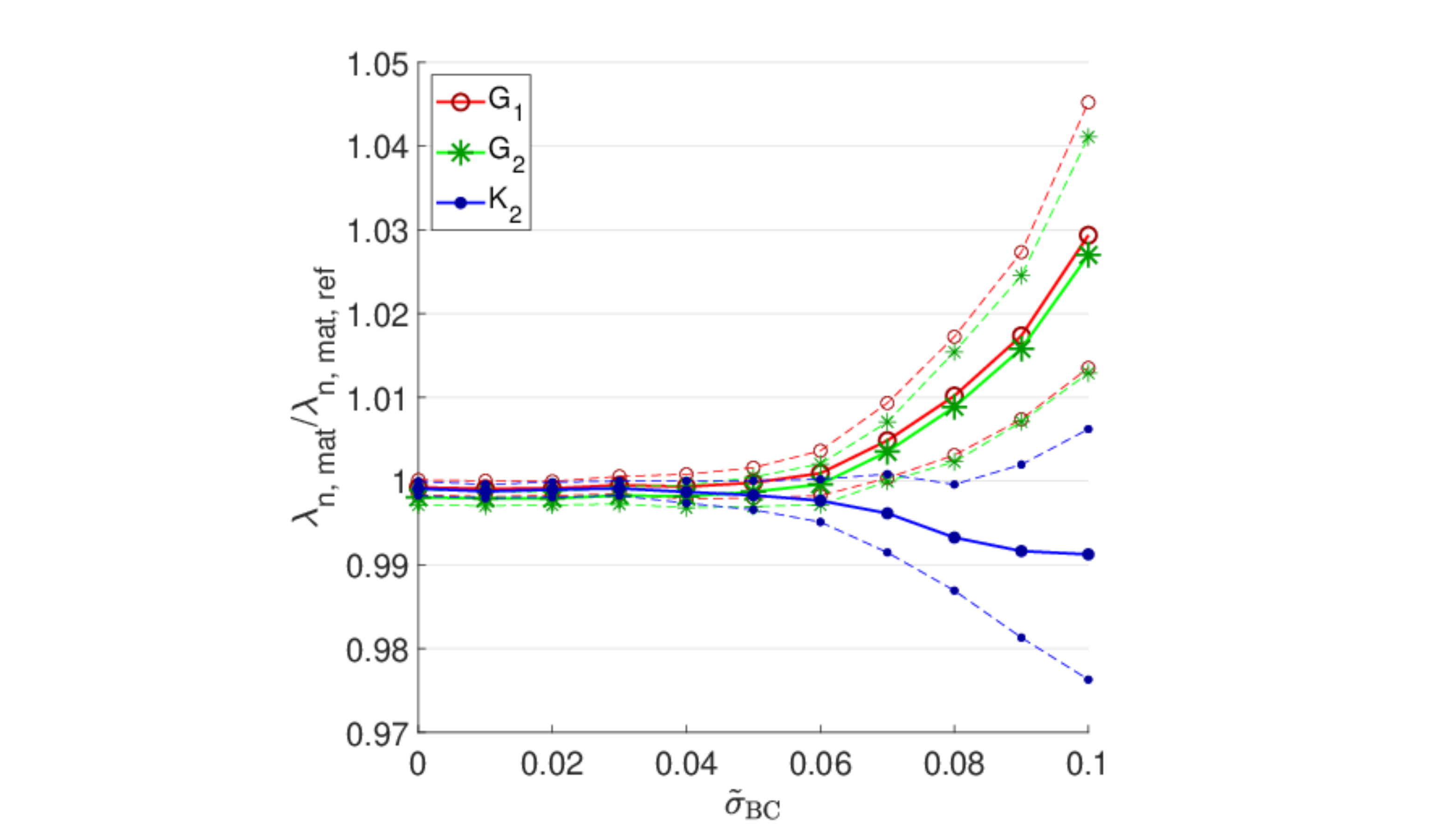}}\\
\subfloat[$t=2\times122$ DOFs (50\% of total).]{
	\includegraphics[trim=6.5cm 0cm 7.5cm 1cm, clip=true, width=0.4\textwidth]{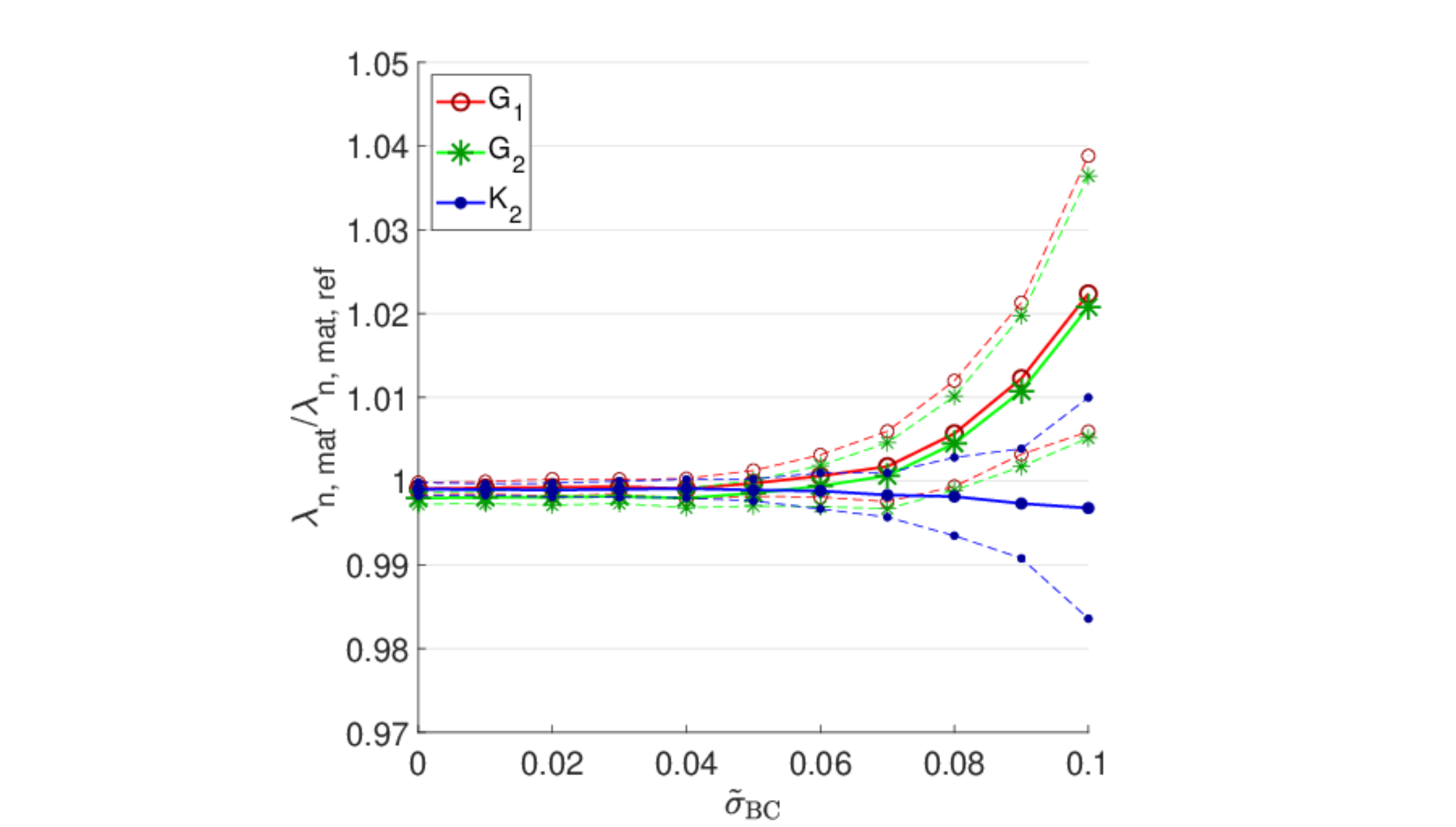}}
	\hspace{0.5em}
\subfloat[$t=2\times61$ DOFs (25\% of total).]{
	\includegraphics[trim=6.5cm 0cm 7.5cm 1cm, clip=true, width=0.4\textwidth]{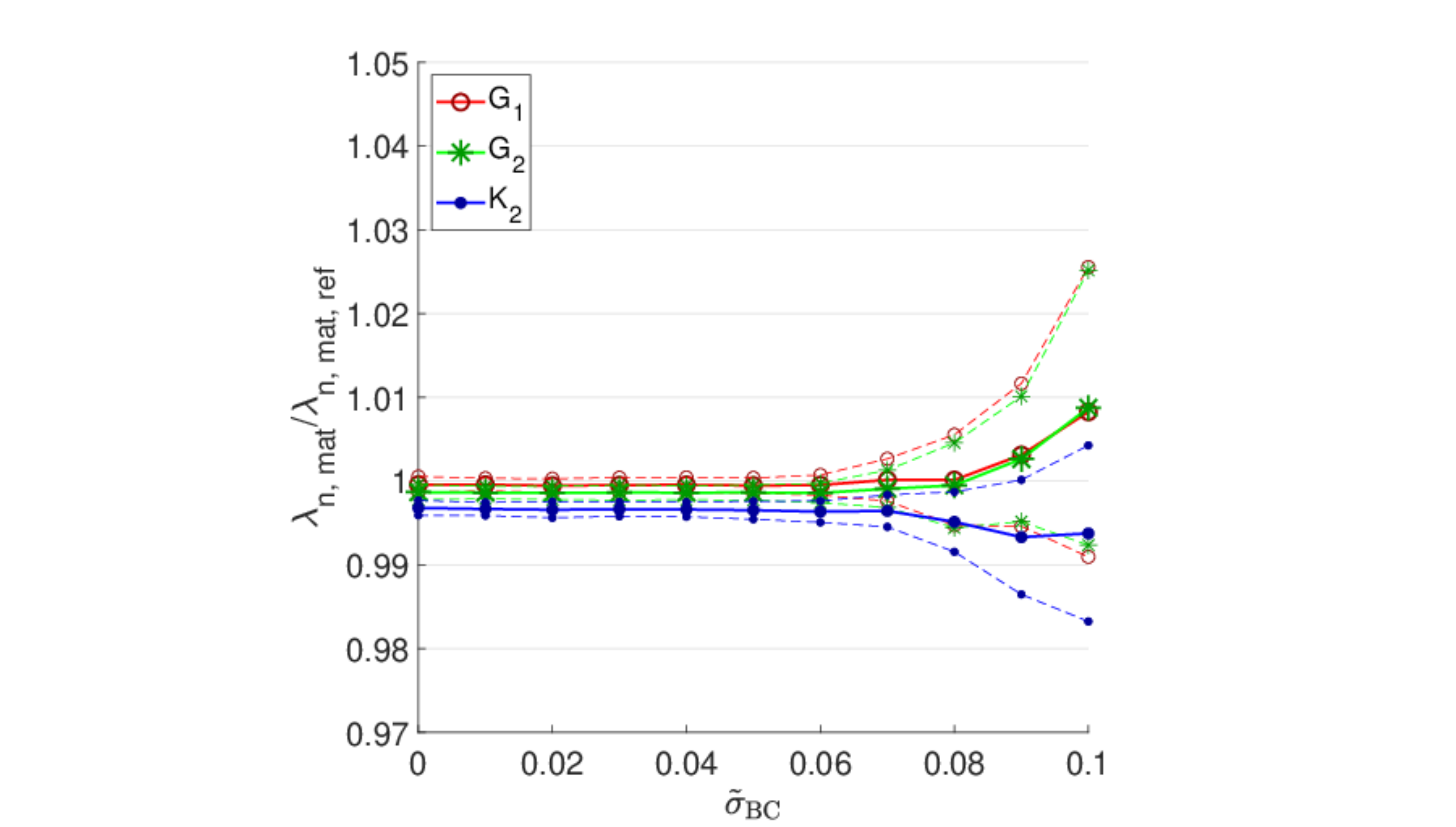}}
\caption{Averaged modes of normalized material parameter distributions $\lambda_{\mathrm{mat},n}/\lambda_{\mathrm {mat},n,\mathrm{ref}}$ for relaxed boundary conditions, cf. Eq.~\eqref{noise_imp}, with the standard deviation $\tilde{\sigma}_{\mathrm{bc}}$ of the starting noise for the tensile test. The means across all sampled values of all iterations (after burn-in) are plotted with solid lines and are complemented with $\pm$ standard deviations (dashed lines) with $N=24\,000$ and burn-in $N_0=22\,000$ steps. Parameter $K_1$ is fixed to the reference value to ensure uniqueness.}
\label{fig:mats_FE_tension}
\end{figure}

In the tensile test (Fig.~\ref{fig:mats_FE_tension}) the resulting error in the material parameters was decreased at least twofold for all the considered numbers of the kinematic DOFs, compared to the fixed boundary MHA (cf. Fig.~\ref{fig:ns_all}; note also the difference in scale between the two figures). The starting boundary error still has a considerable negative impact on the overall accuracy. It should be noted that the chosen fixed number of steps of MHA, although relatively high, was insufficient, as the convergence was not reached in most of the random chains for the starting noise amplitude $\tilde{\sigma}_{\mathrm{bc}} \geq 0.04$ in the experiments with the number of kinematic DOFs equal to 50, 75, and 100\% of the total, and $\tilde{\sigma}_{\mathrm{bc}} \geq 0.06$ for the 25\% of the total kinematic DOFs.

\begin{figure}[!htbp]
\centering
\subfloat[$t=2\times244$ (100\% of total), $\tilde{\sigma}_{\mathrm{bc}} = 0$.]{
	\includegraphics[trim=5.5cm 0cm 7.5cm 1cm, clip=true, width=0.41\textwidth]{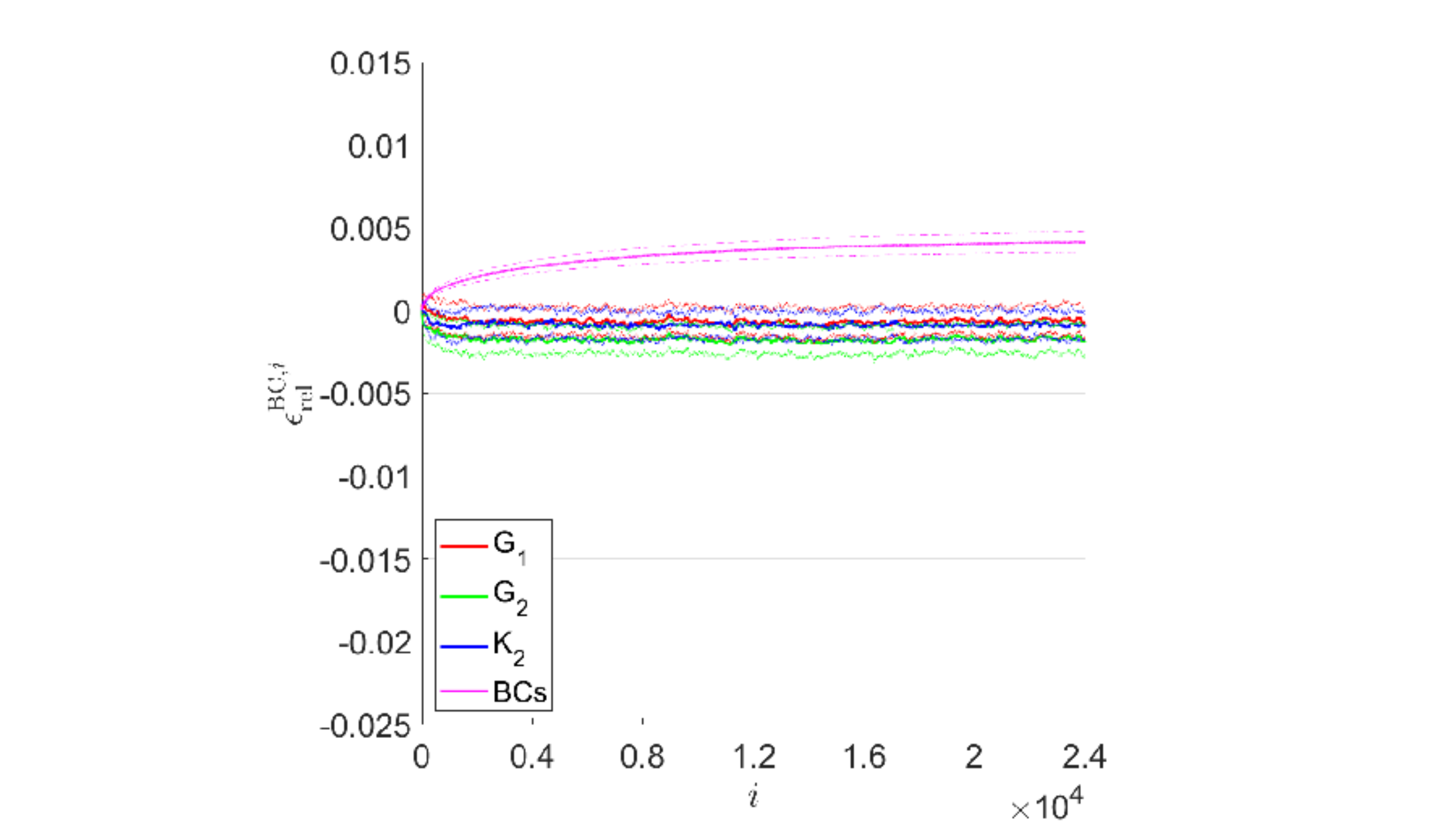}
	\label{fig:MH_steps_sigma_tension_1}}
\subfloat[$t=2\times61$ (25\% of total), $\tilde{\sigma}_{\mathrm{bc}} = 0$.]{
	\includegraphics[trim=5.5cm 0cm 7.5cm 1cm, clip=true, width=0.41\textwidth]{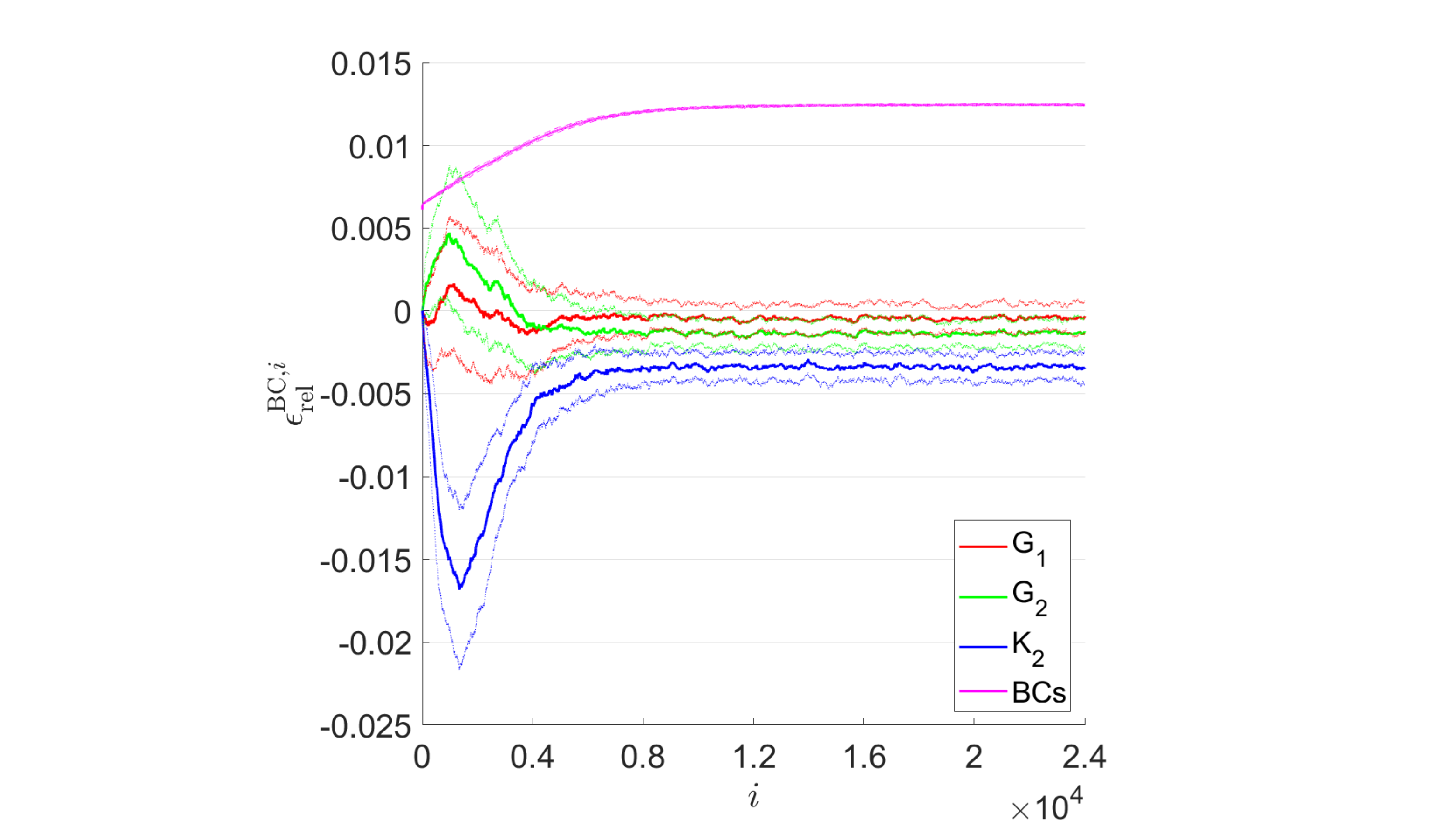}
	\label{fig:MH_steps_sigma_tension_2}}\\
\subfloat[$t=2\times244$ (100\% of total), $\tilde{\sigma}_{\mathrm{bc}} = 0.1$.]{
	\includegraphics[trim=5.5cm 0cm 7.5cm 1cm, clip=true, width=0.41\textwidth]{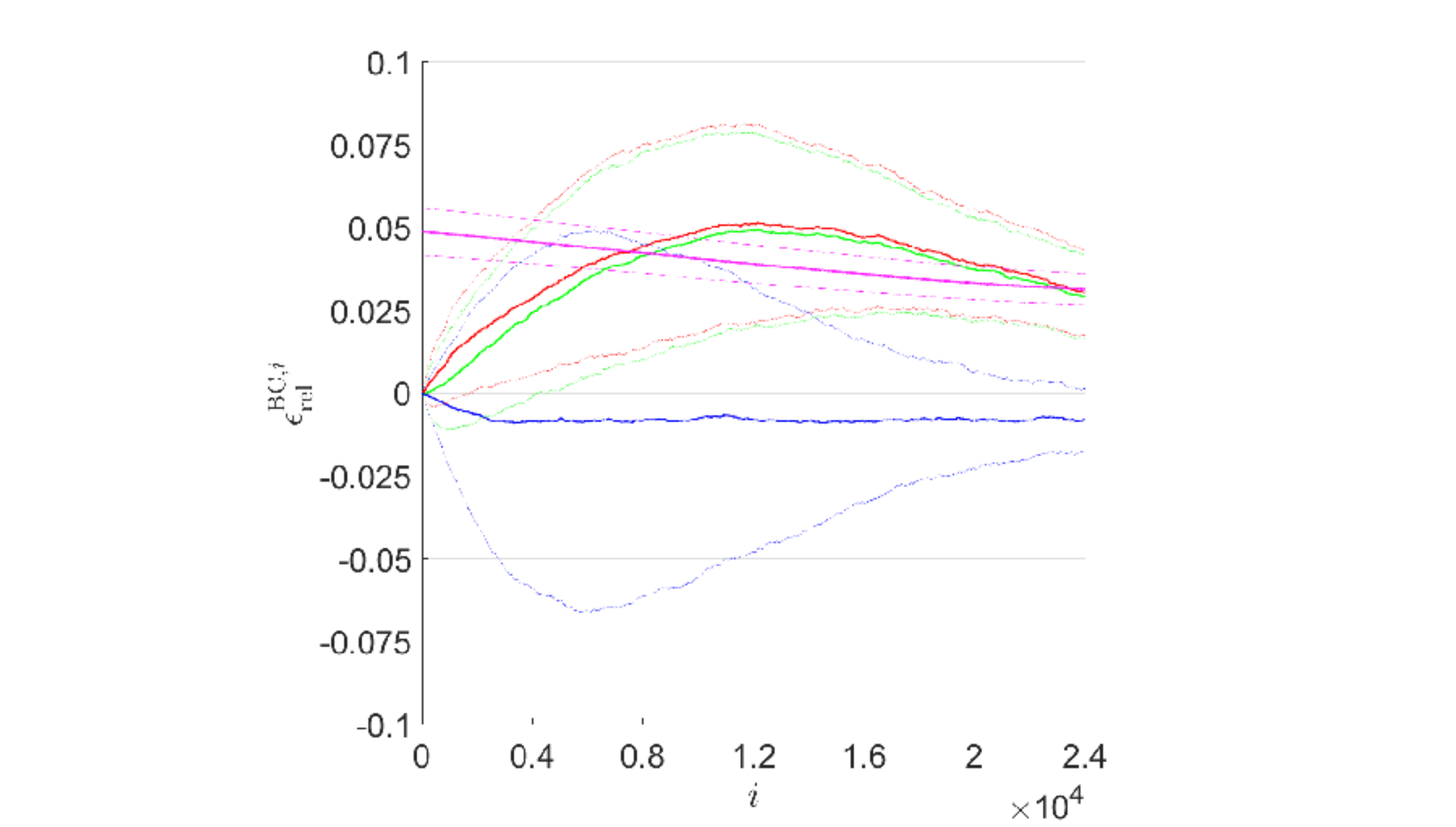}
	\label{fig:MH_steps_sigma_tension_3}}
\subfloat[$t=2\times61$ (25\% of total), $\tilde{\sigma}_{\mathrm{bc}} = 0.1$.]{
	\includegraphics[trim=5.5cm 0cm 7.5cm 1cm, clip=true, width=0.41\textwidth]{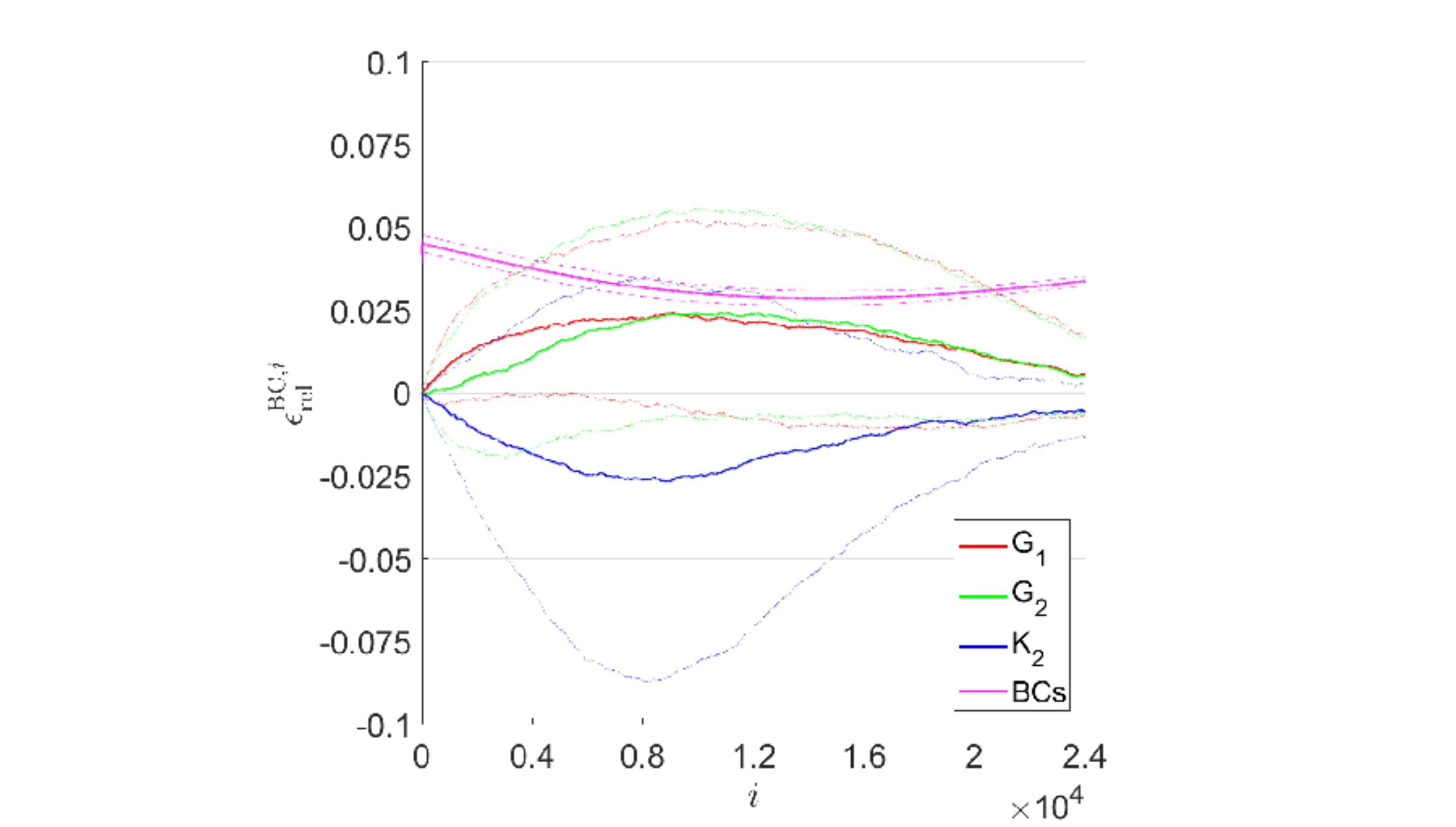}
	\label{fig:MH_steps_sigma_tension_4}}
\caption{The evolution of the relative error in the identified material parameters and boundary conditions, cf. Eqs.~\eqref{error-total-mat} and~\eqref{error-total-kin}, as a function of MHA step $i$ for the tensile test. The means across all sampled values of all iterations (after burn-in) are plotted with solid lines and are complemented with $\pm$ standard deviations (dashed lines). Parameter $K_1$ is fixed to the reference value to ensure uniqueness.}
\label{fig:MH_steps_sigma_tension}
\end{figure}

To better illustrate the mutual influence of the material and kinematic parameters, each parameter's average relative error for a given step across all iterations is shown in Fig.~\ref{fig:MH_steps_sigma_tension}. For brevity, we focus only on the cases with either minimum or maximum noise, with the relative number of kinematic nodes equal to 100\% and 25\%. For the material parameters the error is calculated as
\begin{equation}
  \epsilon^{n,i}_{\mathrm{rel}} = \frac{\widehat{\lambda}_{n}^i-\lambda_{n, \mathrm{ref}}}{\lambda_{n,\mathrm{ref}}}, \;\; n=1,\dots,n_{\mathrm{mat}},
  \label{error-total-mat}
\end{equation}
while for the kinematic DOFs as
\begin{equation}
  \epsilon^{\mathrm{bc},i}_{\mathrm{rel}} = \frac{||\widehat{\bm{\lambda}}_{\mathrm{kin}}^i-\bm{\lambda}_{\mathrm{kin,ref}}||_2}{||\bm{\lambda}_{\mathrm{kin,ref}}||_2},
  \label{error-total-kin}
\end{equation}
where the hatted variables with index $i$ denote $i$-th accepted sample. 

For the zero starting noise in the boundary conditions, the material parameters in the experiment with the full number of DOFs converge very fast to a high-likelihood region close to the initialization point (i.e., to the reference material values, see Figs.~\ref{fig:MH_steps_sigma_tension_1}, \ref{fig:MH_steps_sigma_tension_2}). The average error in the boundary conditions slightly grows throughout the experiment. This is caused by the numerical errors in the process, such as image interpolation.
The experiment with the reduced number of DOFs starts with a higher error in the boundary conditions due to the interpolation error, and after a period of linear growth stabilizes. The growth of the average error in the boundary conditions in this experiment coincides with the rapid spread in the material parameters that stabilize in local optima around the same time as the kinematic parameters. 

In the experiment with the highest starting noise, $\tilde{\sigma}_{\mathrm{bc}} = 0.1$, shown in Figs.~\ref{fig:MH_steps_sigma_tension_3} and~\ref{fig:MH_steps_sigma_tension_4}, the initial average boundary condition error is roughly the same in both cases with full and reduced number of kinematic DOFs. The higher regularity of the boundary, that results from the interpolation, has a different effect on each material parameter: while the highly correlated matrix and inclusions shear moduli $G_1$ and $G_2$ reach a higher interim error in the experiment with the full number of kinematic DOFs, the parameter $K_2$ on average stabilizes faster than for the reduced number of kinematic DOFs, although with a high variance between the chains. At first, the boundary error decreases faster for the reduced number of DOFs, but it starts increasing around the halfway of the total running time. The boundary error in the experiment with full boundary DOFs decreases linearly. It is unclear if it would show a similar behavior as its counterpart, as the error at the final ($25\,000$th) step of the non-reduced problem was comparable  to the one observed  for the reduced system around the 12~500th step.

\begin{figure}[!htbp]
\centering
	\subfloat[$t=2\times244$ kinematic DOFs (100\% of total).]{
	\includegraphics[trim=0cm -0.1cm 0cm 0cm, clip=true, width=0.48\textwidth]{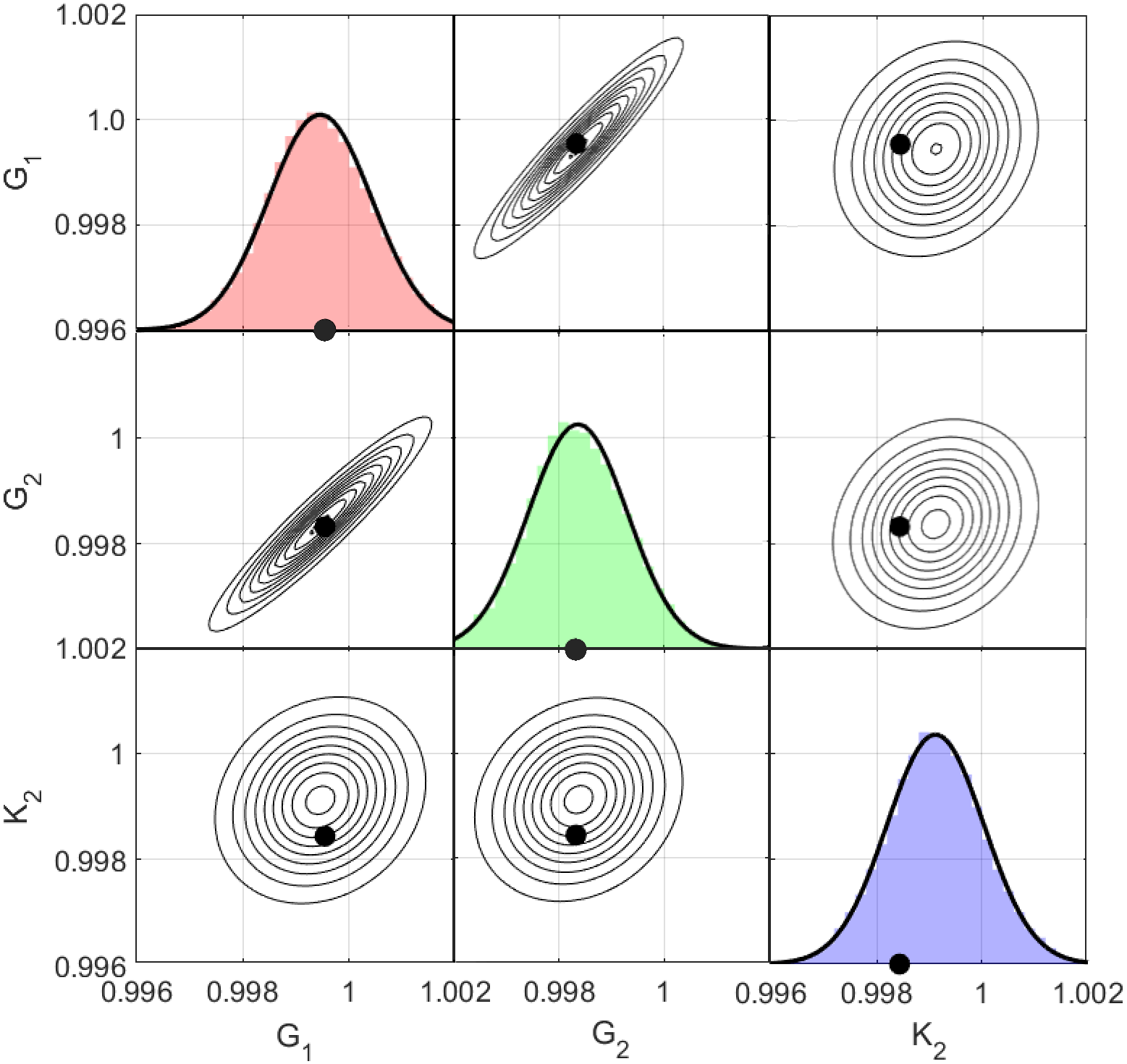}
    \label{fig:hist_q=1.0}}
	\subfloat[$t=2\times61$ kinematic DOFs (25\% of total).]{
    \includegraphics[trim=0cm 0cm 0cm 0cm, clip=true, width=0.48\textwidth]{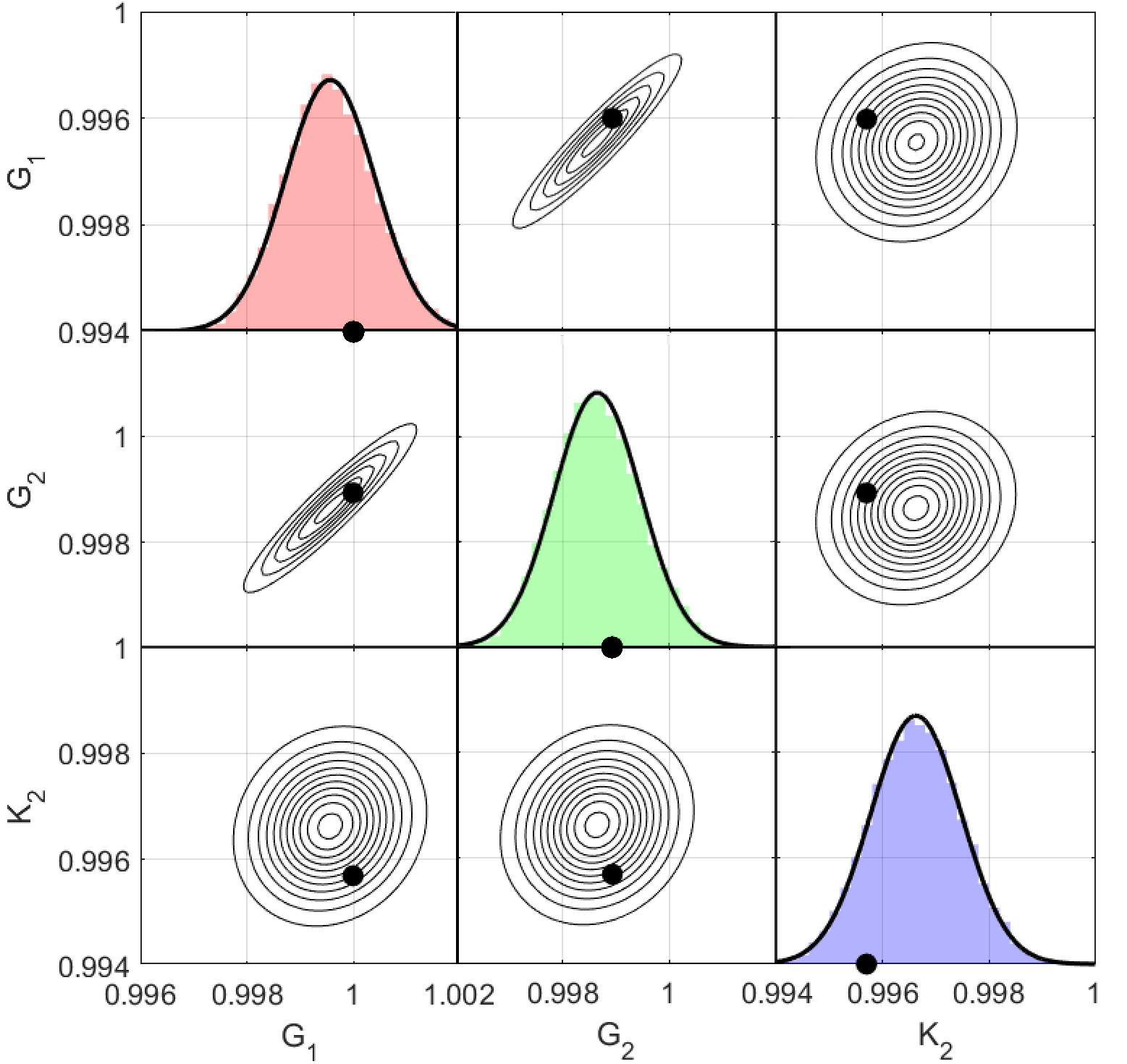}
	\label{fig:hist_q=0.25}}
    \caption{Posterior probability distributions for the material parameters identified by MHA from the tensile test,  and point estimations given by BE-IDIC (shown as black dots). Parameter $K_1$ is fixed to the reference value to ensure uniqueness.}
\end{figure}

Finally, all the converged chains from the above experiments (those with a lower amount of initial noise) are congregated to estimate the posterior probability density functions for the material parameters. Those, in turn, can be compared with the results obtained by the deterministic BE-IDIC method. In Figs.~\ref{fig:hist_q=1.0} and \ref{fig:hist_q=0.25} it is shown that all the resulting posterior distributions are approximately normal, with the experiments with reduced number of boundary nodes resulting in a slightly narrower confidence intervals. Solutions provided by the deterministic method (denoted by black dots) tend to coincide with the modes of the histograms (or means of the estimated PDFs), except for the parameter $K_2$, where the average of the stochastic method is slightly closer to the true parameter value in both cases. 

\subsection{Shear Test}
\begin{figure}[!htbp]
\centering
	\subfloat[$t=2\times244$ DOFs (100\% of total).]{
	\includegraphics[trim=6.5cm 0cm 7.5cm 1cm, clip=true, width=0.4\textwidth]{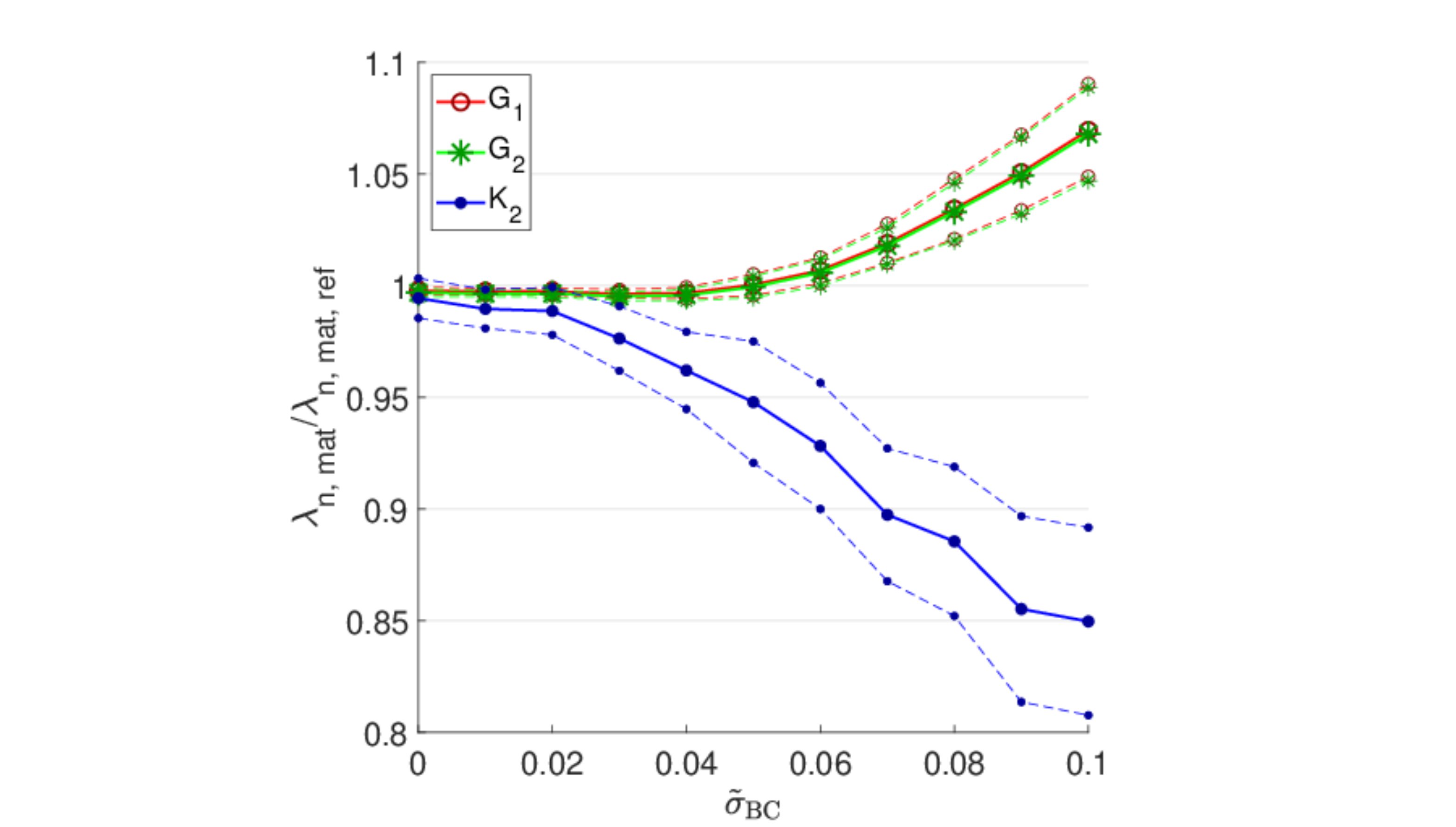}}
	\subfloat[$t=2\times183$ DOFs (75\% of total).]{\includegraphics[trim=6.5cm 0cm 7.5cm 1cm, clip=true, width=0.4\textwidth]{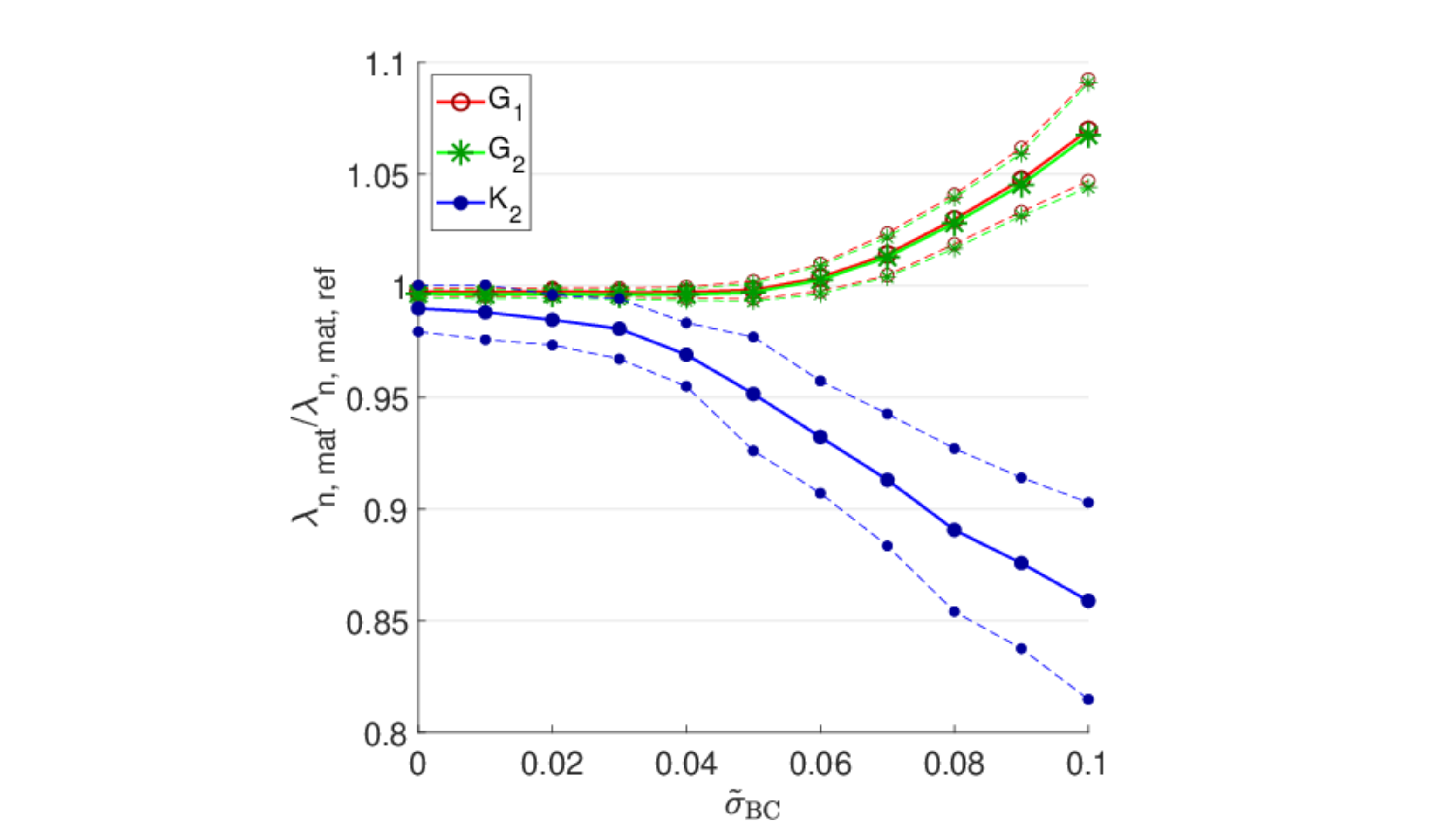}}\\
	\subfloat[$t=2\times122$ DOFs (50\% of total).]{\includegraphics[trim=6.5cm 0cm 7.5cm 1cm, clip=true, width=0.4\textwidth]{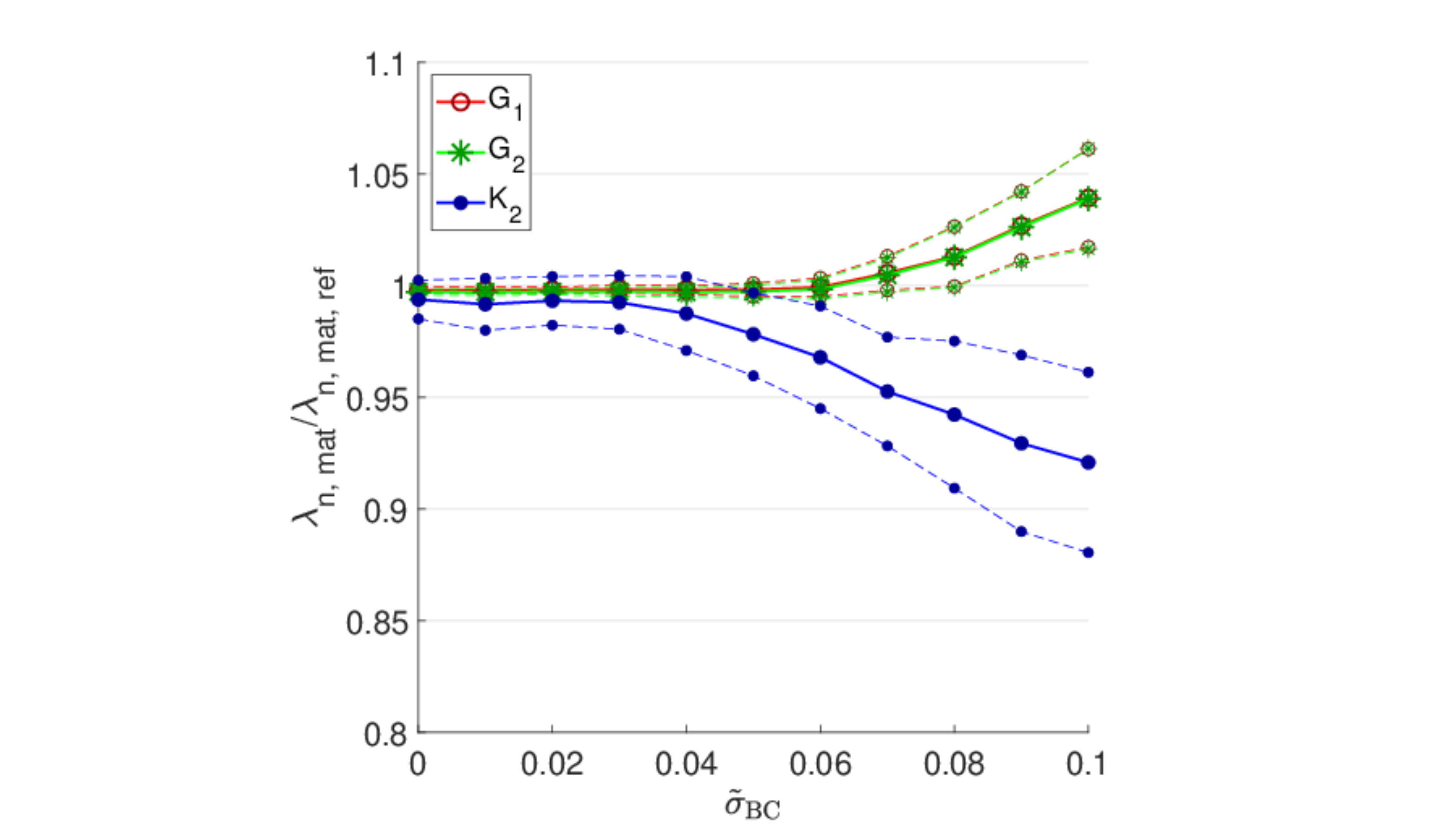}}
	\subfloat[$t=2\times61$ DOFs (25\% of total).]{
	\includegraphics[trim=6.5cm 0cm 7.5cm 1cm, clip=true, width=0.4\textwidth]{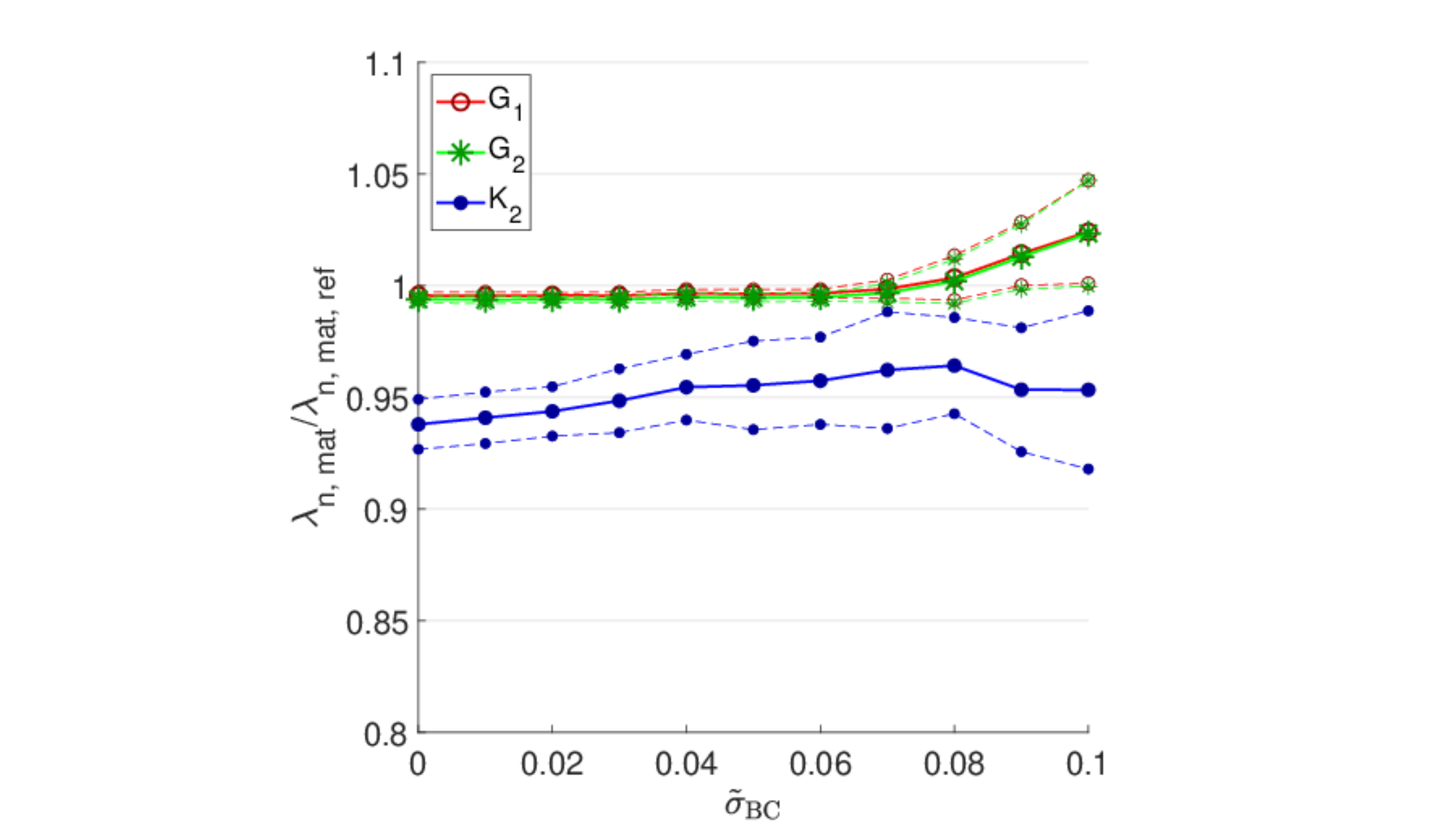}}
	\caption{Averaged modes of normalized material parameter distributions $\lambda_{\mathrm{mat},n}/\lambda_{\mathrm{mat},n,\mathrm{ref}}$ for relaxed boundary conditions, cf. Eq.~\eqref{noise_imp}, with the starting noise standard deviation $\tilde{\sigma}_{\mathrm{bc}}$ for the shear test. The means across all sampled values of all iterations (after burn-in) are plotted with solid lines and are complemented with $\pm$ standard deviations (dashed lines) with $N=24\,000$ and burn-in $N_0=22\,000$ steps. Parameter $K_1$ is fixed to the reference value to ensure uniqueness.}
	\label{fig:mats_FE_shear}
\end{figure}

\begin{figure}[!htbp]
\centering
    \subfloat[$t=2\times244$ (100\% of total), $\tilde{\sigma}_{\mathrm{bc}}=0$.]{
	\includegraphics[trim=6cm 0cm 8cm 1cm, clip=true, width=0.41\textwidth]{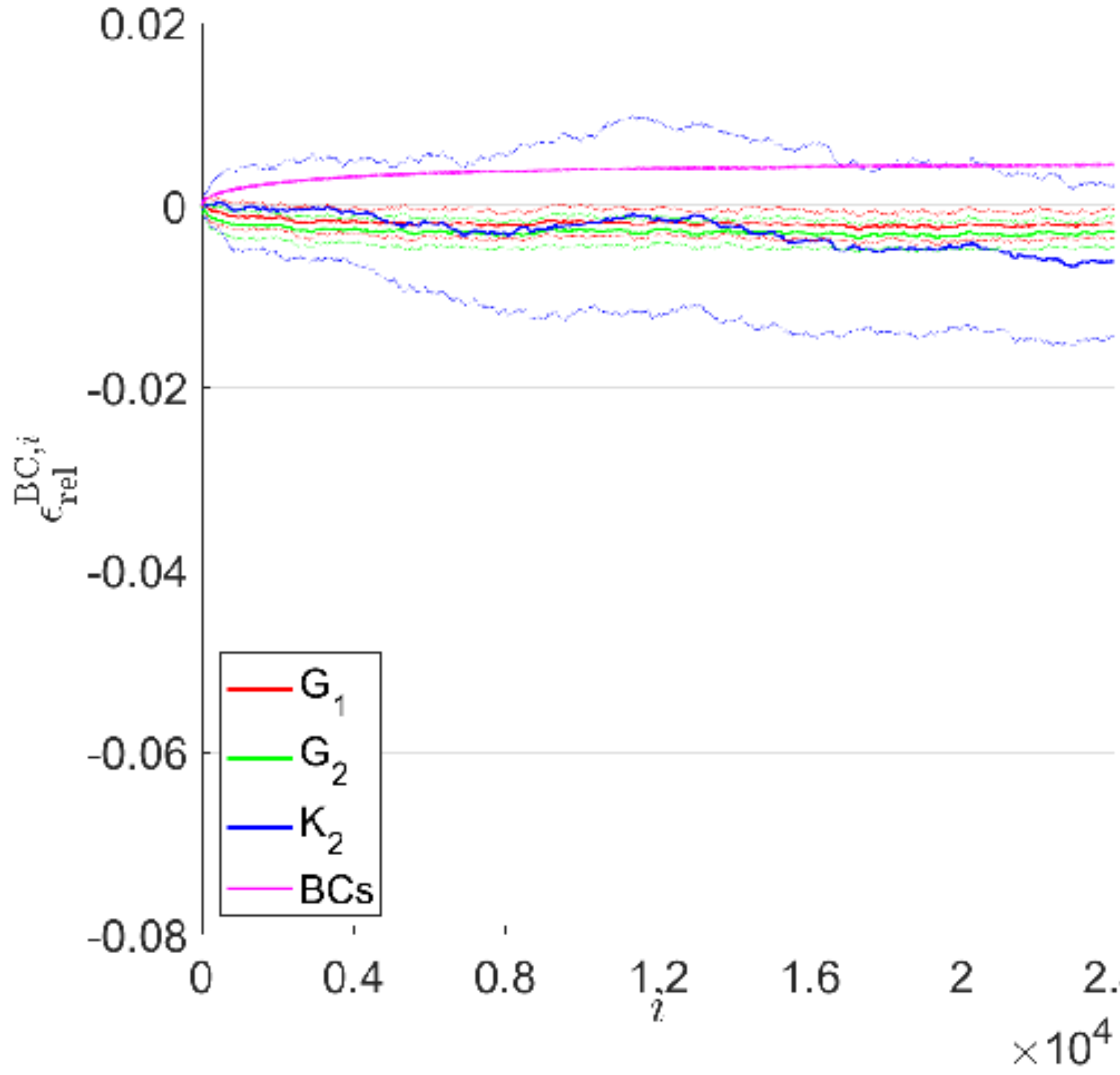}}
	\subfloat[$t=2\times61$ (25\% of total), $\tilde{\sigma}_{\mathrm{bc}}=0$.]{
	\includegraphics[trim=6cm 0cm 8cm 1cm, clip=true, width=0.41\textwidth]{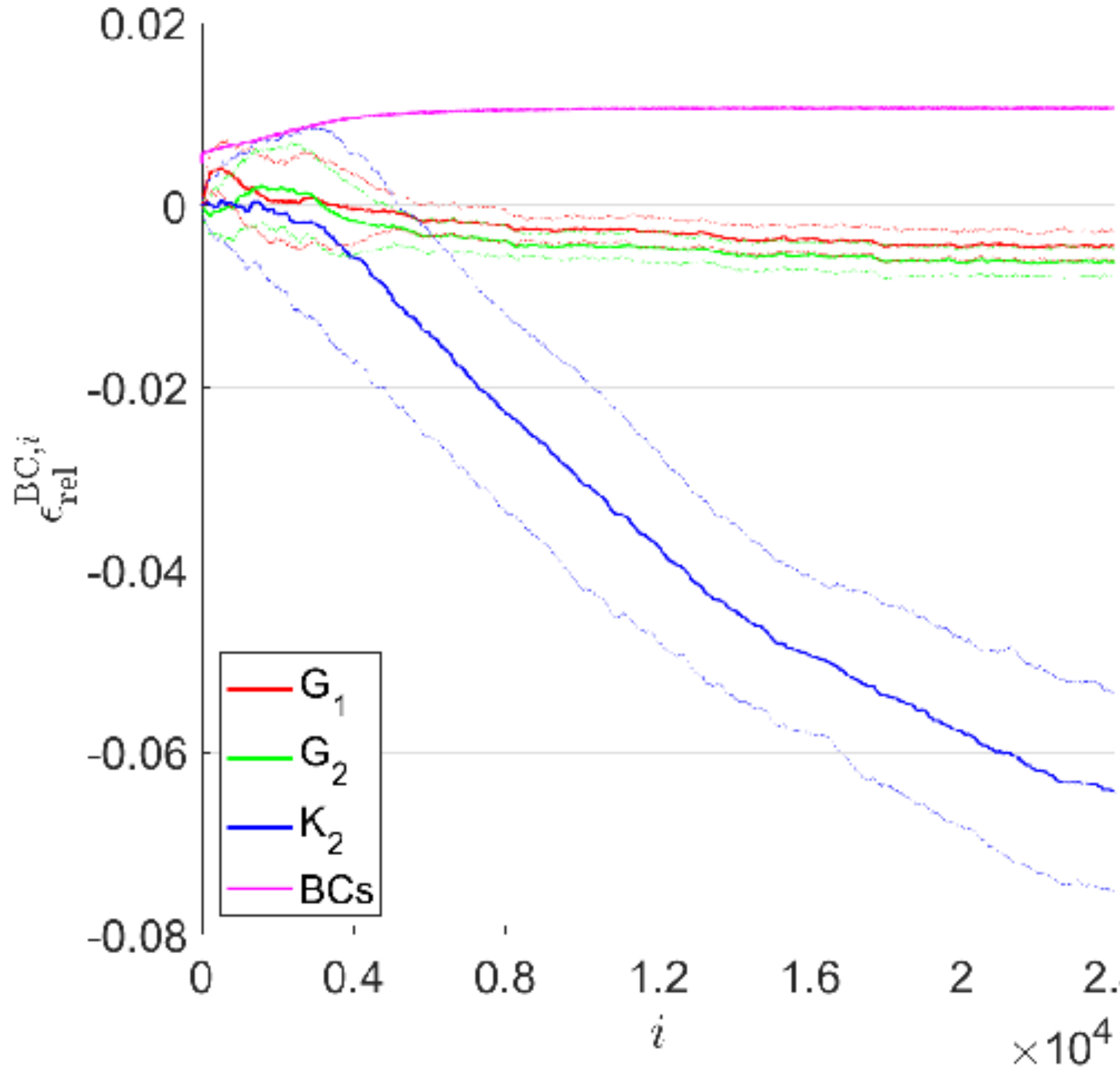}}\\
	\subfloat[$t=2\times244$ (100\% of total), $\tilde{\sigma}_{\mathrm{bc}}=0.1$.]{
	\includegraphics[trim=6cm 0cm 8cm 1cm, clip=true, width=0.41\textwidth]{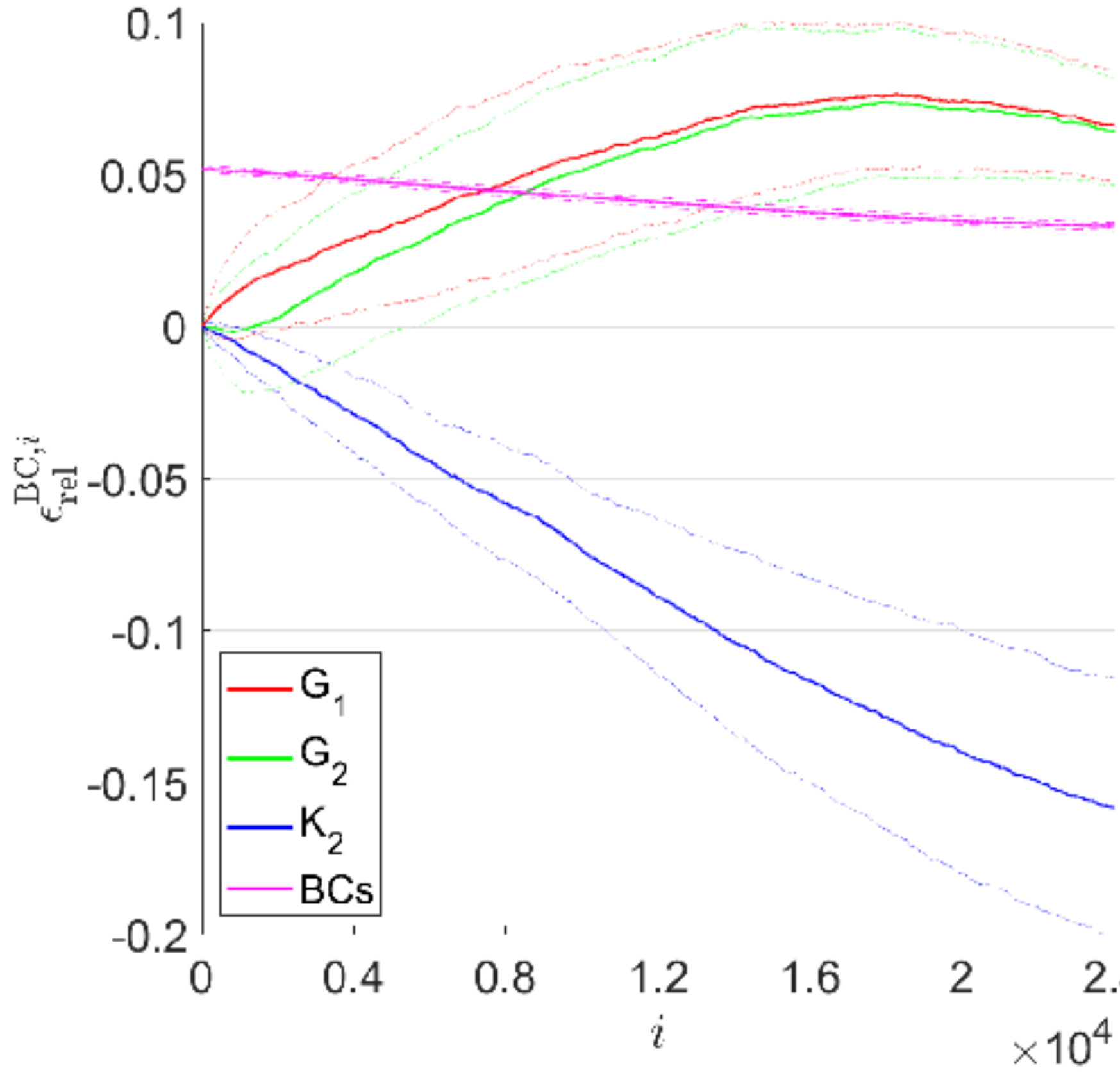}}
	\subfloat[$t=2\times61$ (25\% of total),$\tilde{\sigma}_{\mathrm{bc}}=0.1$.]{
	\includegraphics[trim=6cm 0cm 8cm 1cm, clip=true, width=0.41\textwidth]{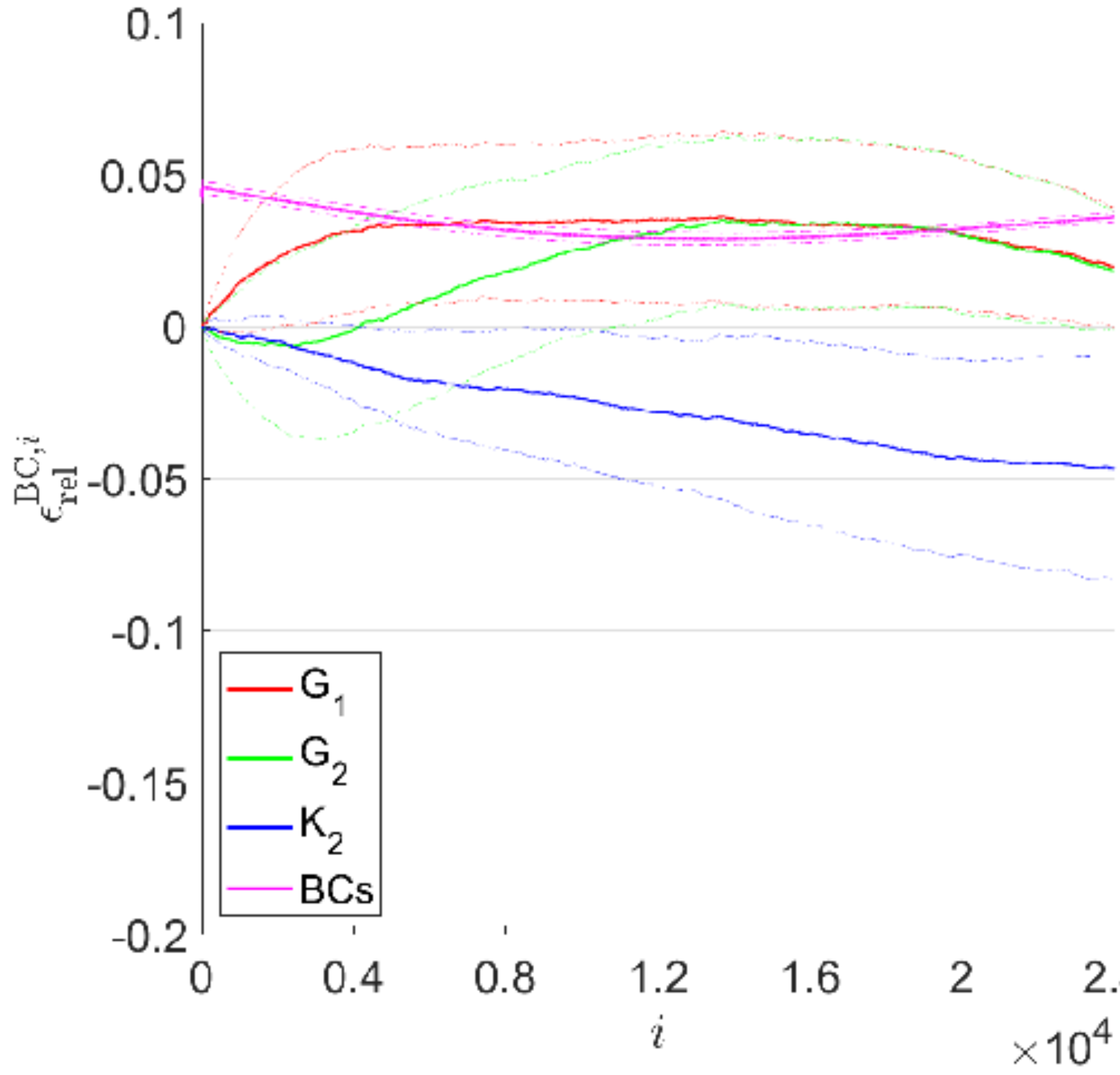}}
\caption{The evolution of the relative error in the identified material parameters and boundary conditions, cf. Eqs.~\eqref{error-total-mat} and~\eqref{error-total-kin}, as a function of MHA step $i$ for the shear test, the means across all sampled values of all iterations (after burn-in) are plotted with solid lines and are complemented with $\pm$ standard deviations (dashed lines). Parameter $K_1$ is fixed to the reference value to ensure uniqueness.}
\label{fig:MH_steps_sigma_shear}
\end{figure}

The shear test behaves quite differently from the tensile test for all numbers of the kinematic DOFs, see Fig.~\ref{fig:mats_FE_shear}. This test seems to have an overall low sensitivity to the bulk modulus of the inclusions $K_2$. The accuracy decreases for the rest of the material parameters as well. Nevertheless, Fig.~\ref{fig:MH_steps_sigma_shear} shows that the average relative error in boundary conditions decreases at the same rate as for the analogous tensile experiments. 

\begin{figure}[!htbp]
\centering
    \subfloat[$t=2\times244$ kinematic DOFs (25\% of total).]{
	\includegraphics[trim=0cm 0cm 0cm 0cm, clip=true, width=0.48\textwidth]{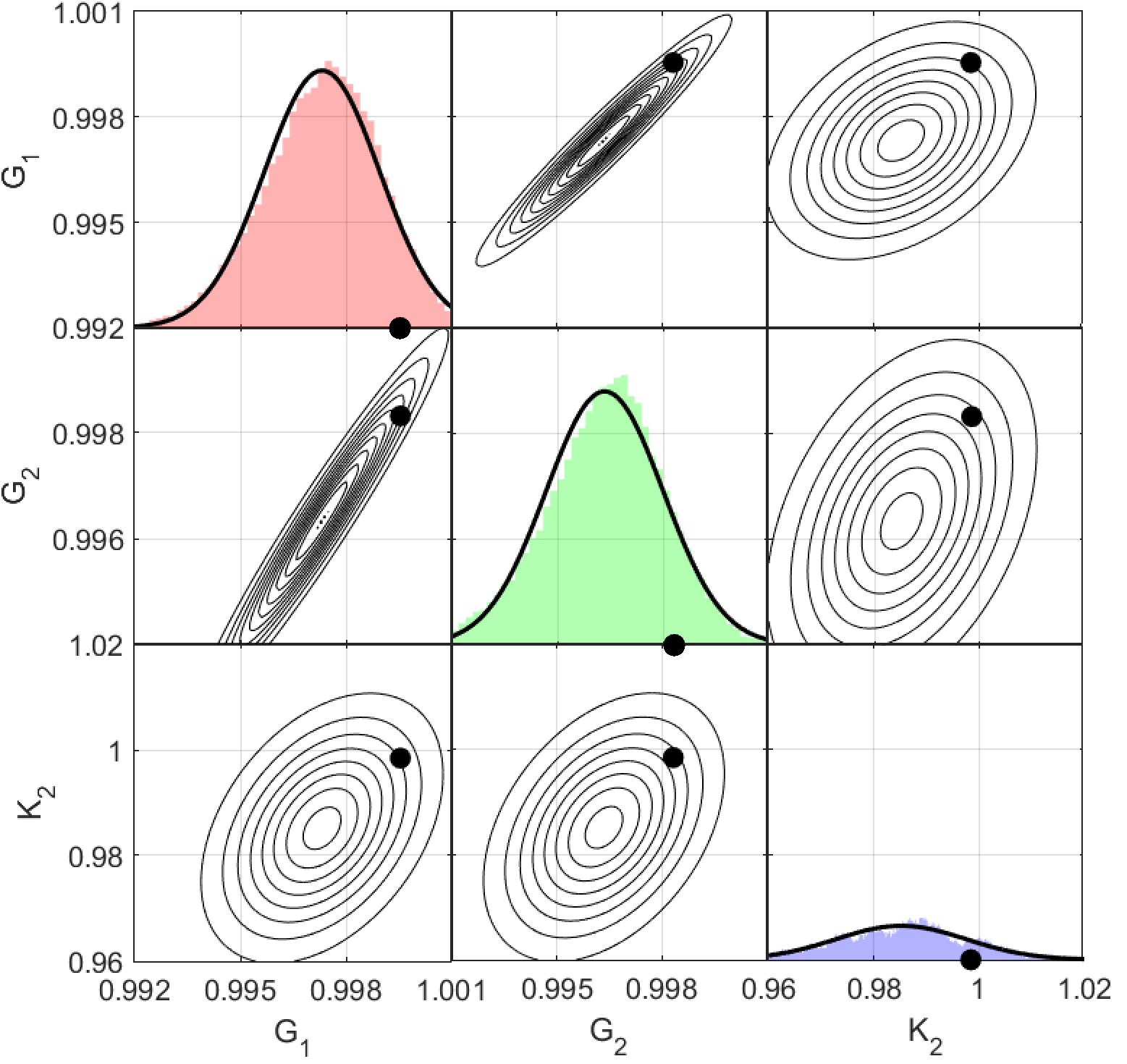}
    \label{fig:hist_q=1.0_shear}}
    \subfloat[$t=2\times61$ kinematic DOFs (25\% of total).]{
    \includegraphics[trim=0cm 0cm 0cm 0cm, clip=true, width=0.48\textwidth]{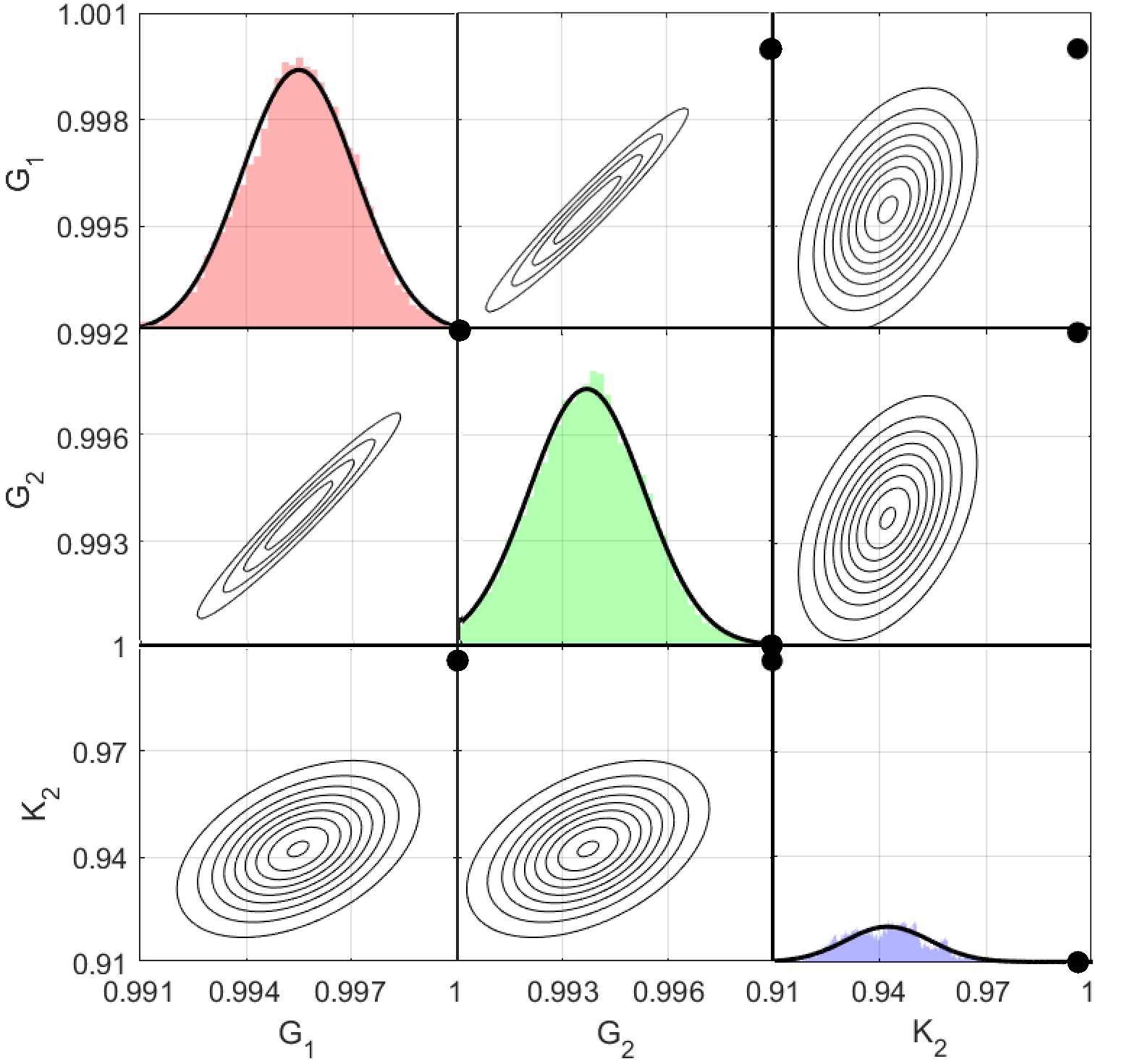}
    \label{fig:hist_q=0.25_shear}}
    \caption{Posterior probability distributions for the material parameters identified from the shear test, and point estimations given by BE-IDIC (shown as black dots). Parameter $K_1$ is fixed to the reference value to ensure uniqueness.}
\end{figure}

For the shear test, the histograms and estimated uni- and multi-variate posterior probability distributions obtained from the converged chains are shown in Figs.~\ref{fig:hist_q=1.0_shear} and \ref{fig:hist_q=0.25_shear}. The reduced number of kinematic DOFs results in a higher systematic error than for the tensile test, as well as wider confidence intervals. The means of the posterior distributions are less accurate as the estimators of the material parameters than the solution obtained by the deterministic BE-IDIC method.

\section{Application in Ill-Posed Problems}\label{section:ill-posed}
\begin{figure}[!htbp]
\centering
	\subfloat[$t=2\times244$ DOFs (100\% of total).]{\includegraphics[trim=6cm 0cm 7.5cm 1cm, clip=true, width=0.4\textwidth]{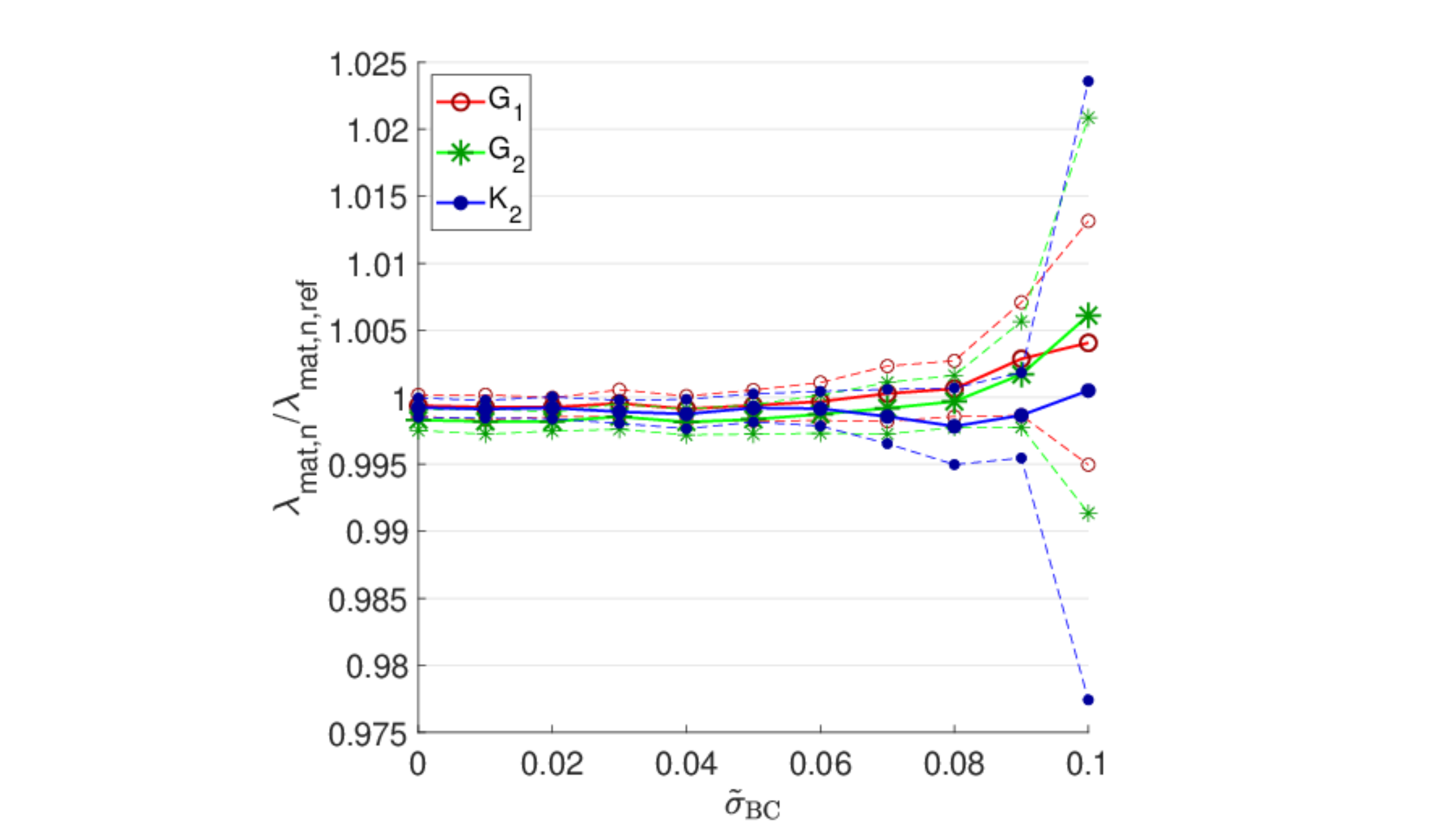}}
	\subfloat[$t=2\times61$ DOFs (25\% of total).]{\includegraphics[trim=6cm 0cm 7.5cm 1cm, clip=true, width=0.4\textwidth]{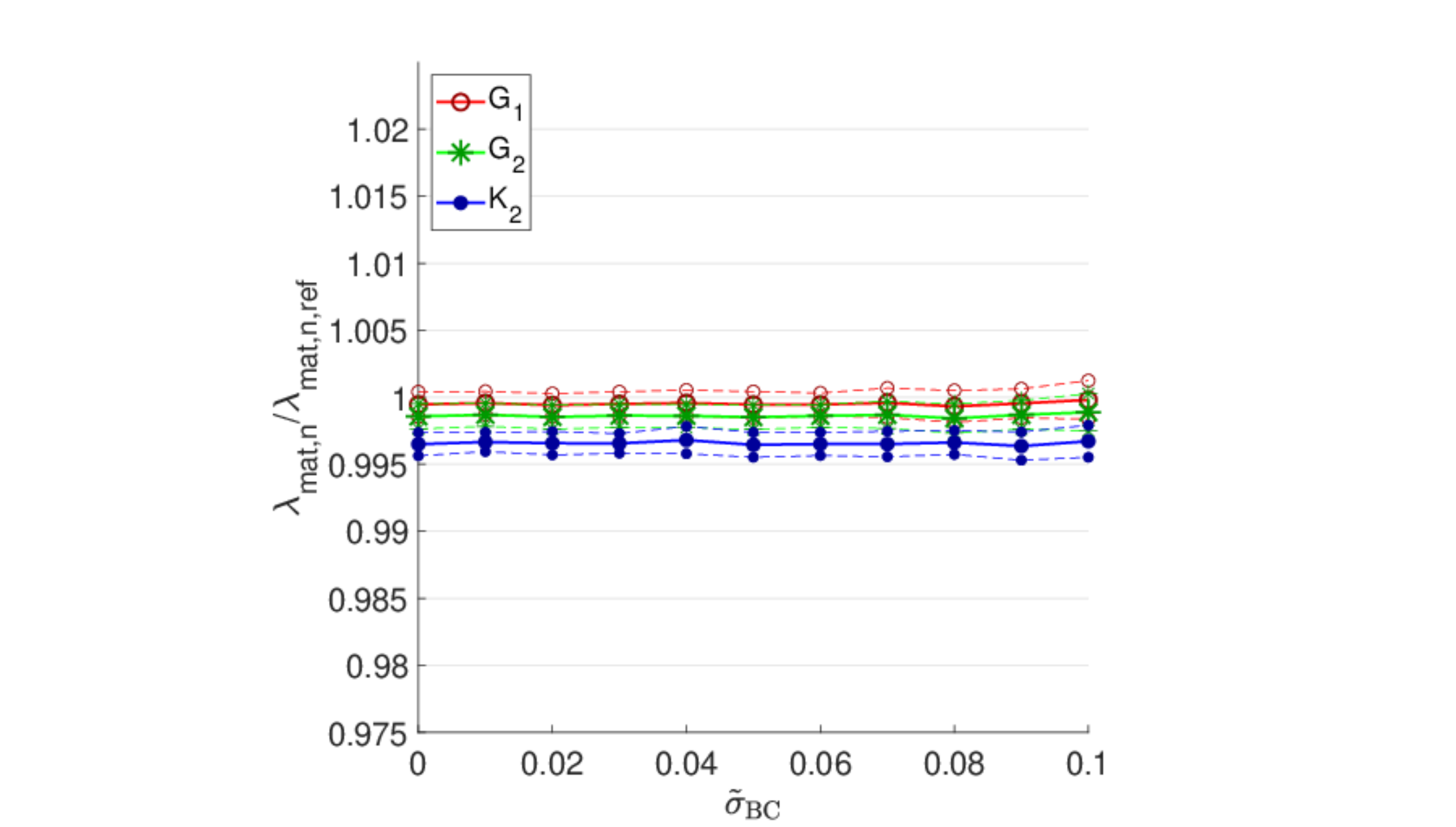}}
\caption{Averaged modes of normalized material parameter distributions $\lambda_{\mathrm{mat},n}/\lambda_{\mathrm{mat},n,\mathrm{ref}}$ for relaxed boundary conditions, cf. Eq.~\eqref{noise_imp}, with the starting noise standard deviation $\tilde{\sigma}_{\mathrm{bc}}$ for the tensile test. The means across all sampled values of all iterations (after burn-in) are plotted with solid lines and are complemented with $\pm$ standard deviations (dashed lines) with $N=24\,000$ and burn-in $N_0=22\,000$ steps. Identified parameter $K_1$ is used as a normalization factor, and is not fixed during identification.}
\label{fig:mats_FE_tension_ill}
\end{figure}

As was discussed earlier, the inverse problem with the Dirichlet boundary conditions requires one of the parameters to be fixed at its reference value, so the deterministic optimizers can only find the correct ratios of the material parameters. The choice of such normalization parameter influences the final accuracy of the identification \citep[cf.][]{rokos}. This, however, is not necessary for the MHA, and the normalization can be done in the post-processing step. 
The sampling can be performed in all four material parameter dimensions instead of three, with the same kinematic parameters if needed. To compare the accuracy of the identification between the ill-posed and the normalized problems, the results from the non-normalized MHA are also normalized in the post-processing step by the parameter $K_1$, which is not fixed at a given value this time, but is updated concurrently with the other parameters. All the other hyper-parameters remain the same, such as material and kinematic parameter step sizes and the prior distributions. After the same number of steps ($24\,000$) and burn-in ($92\%$), the identification accuracy of material parameters is generally much higher for the non-normalized MHA, see Figs.~\ref{fig:mats_FE_tension_ill} and \ref{fig:mats_FE_shear_ill}, where convergence of all material parameters is achieved for a much higher starting noise in the same number of steps.

\subsection{Tensile Test}

\begin{figure}[!htbp]
\centering
	\subfloat[Normalized MHA, $\tilde{\sigma}_{\mathrm{bc}}=0$.]{
	\includegraphics[trim=6cm 0cm 7.5cm 0.5cm, clip=true, width=0.4\textwidth]{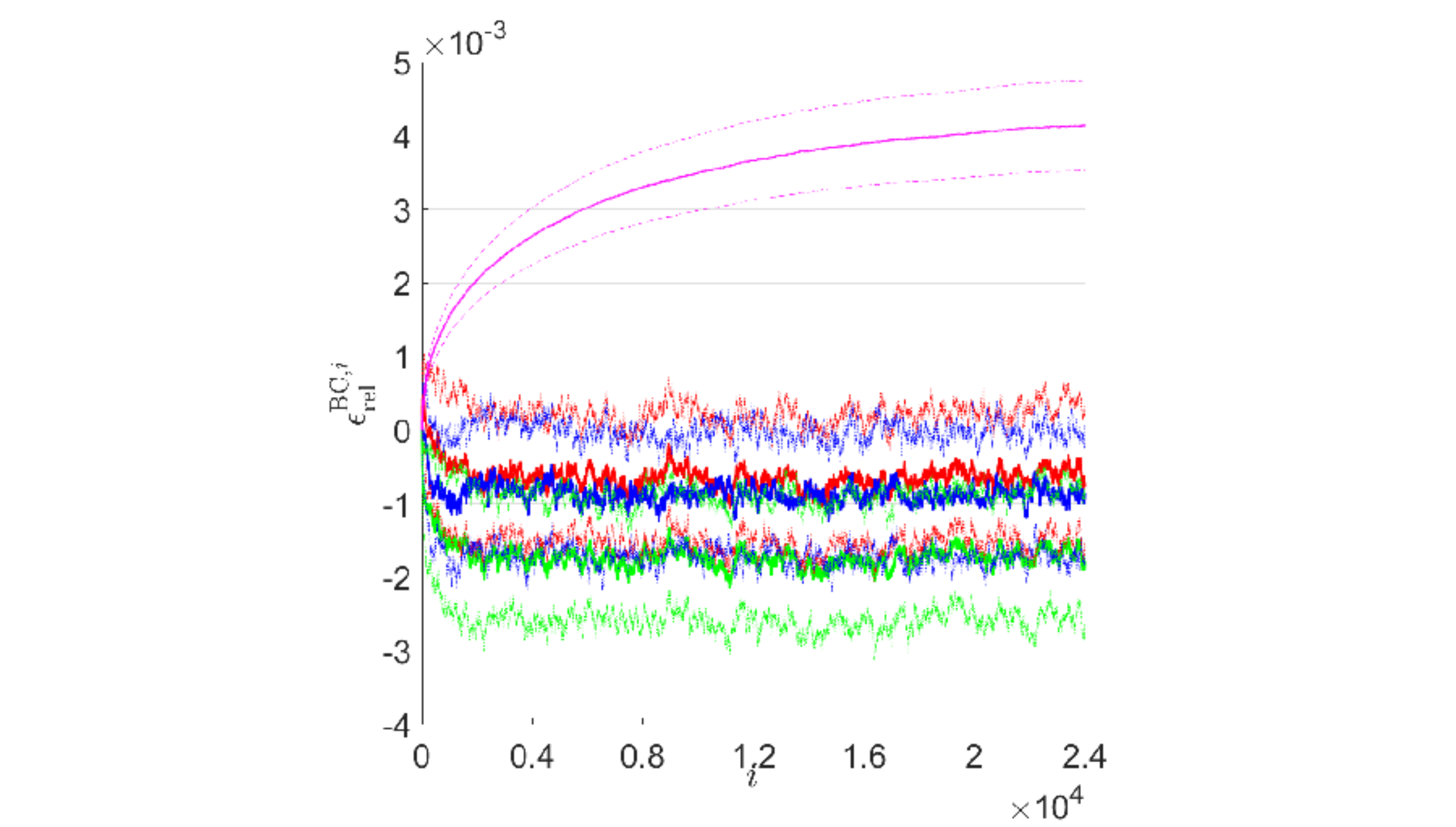}}
    \subfloat[Non-normalized MHA, $\tilde{\sigma}_{\mathrm{bc}}=0$.]{
	\includegraphics[trim=6cm 0cm 7.5cm 0.5cm, clip=true, width=0.4\textwidth]{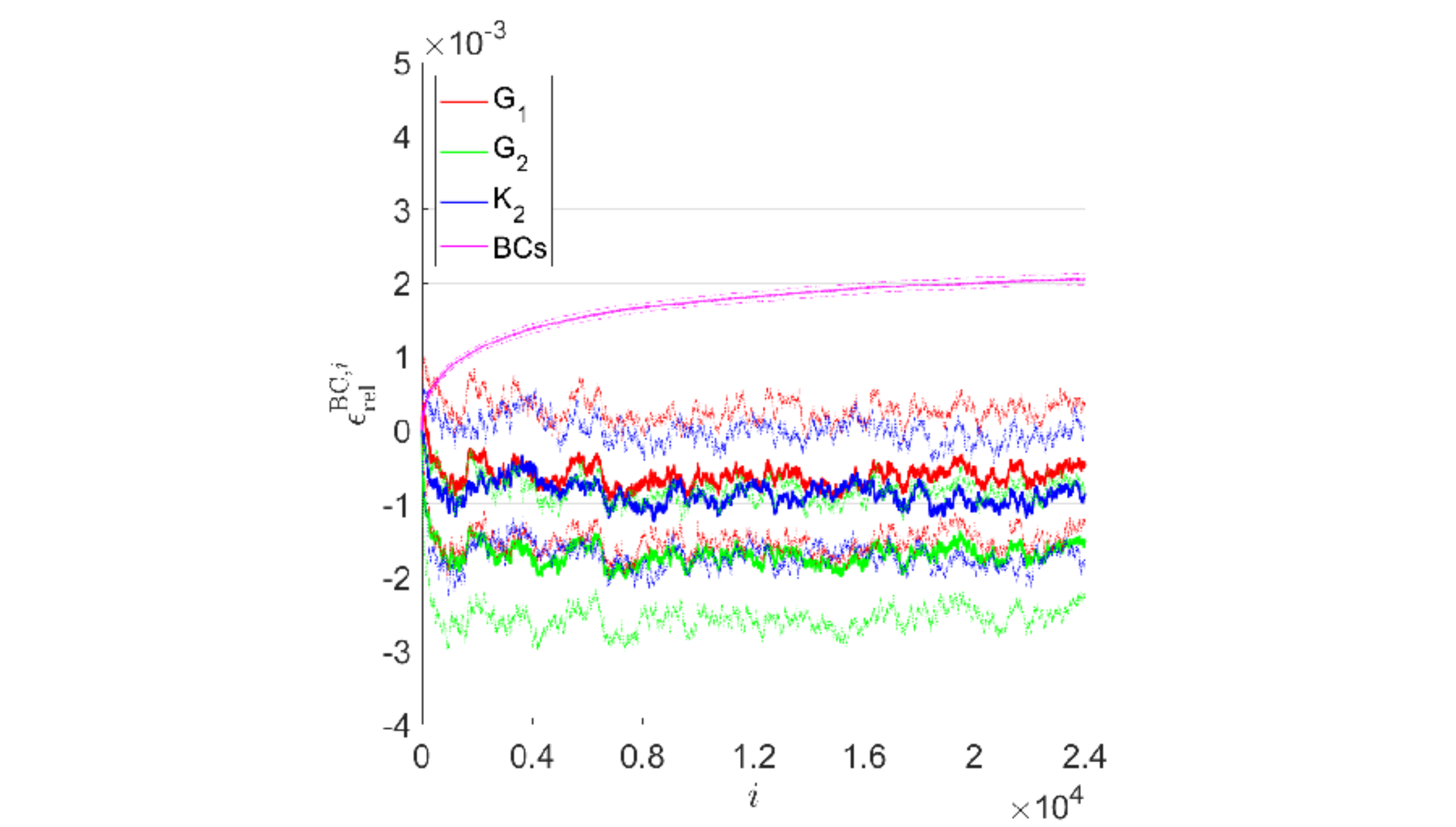}}\\
	\subfloat[Normalized MHA, $\tilde{\sigma}_{\mathrm{bc}}=0.1$.]{
	\includegraphics[trim=6cm 0cm 7.5cm 0.5cm, clip=true, width=0.4\textwidth]{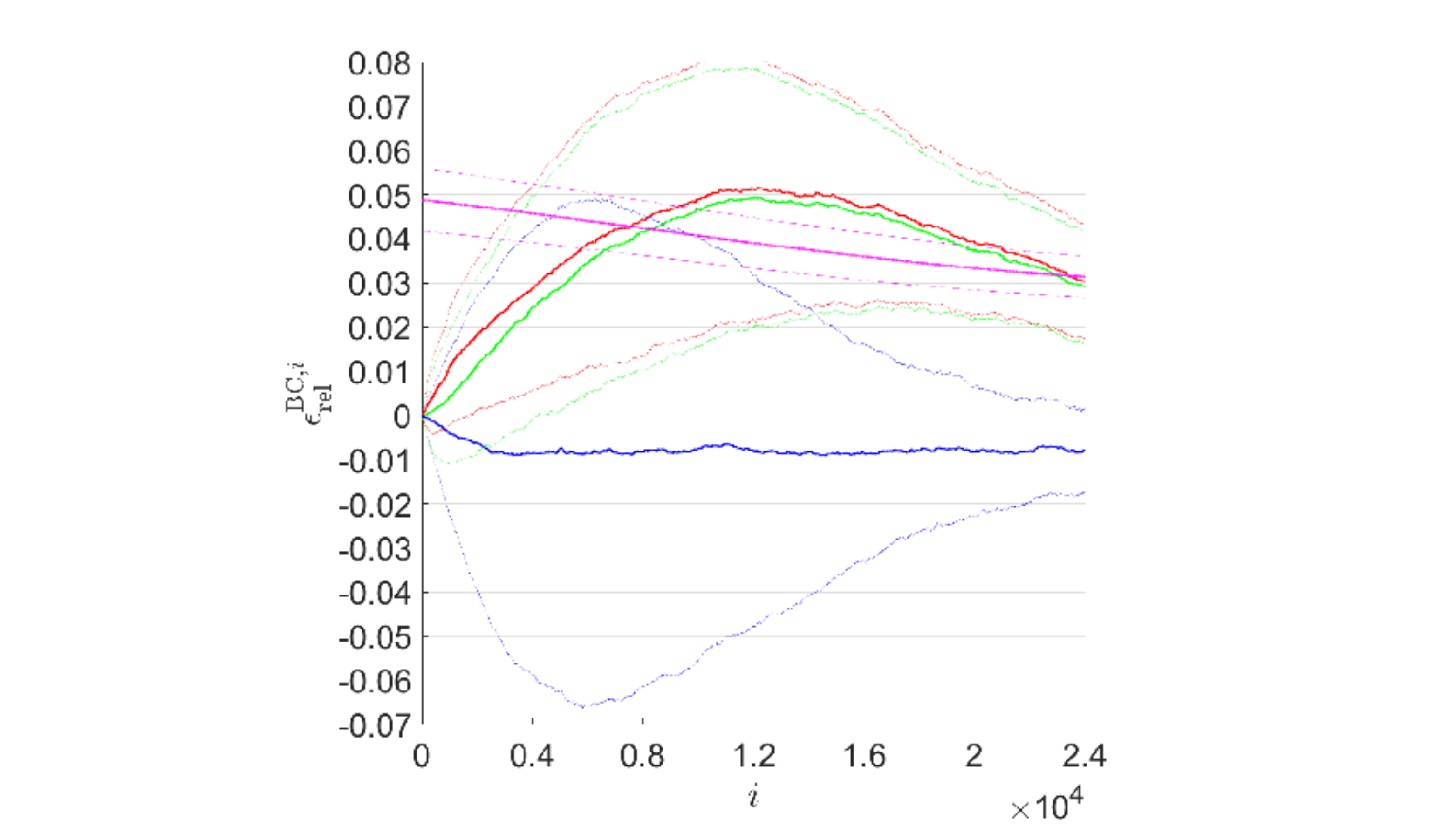}}
    \subfloat[Non-normalized MHA, $\tilde{\sigma}_{\mathrm{bc}}=0.1$.]{
	\includegraphics[trim=6cm 0cm 7.5cm 0.5cm, clip=true, width=0.4\textwidth]{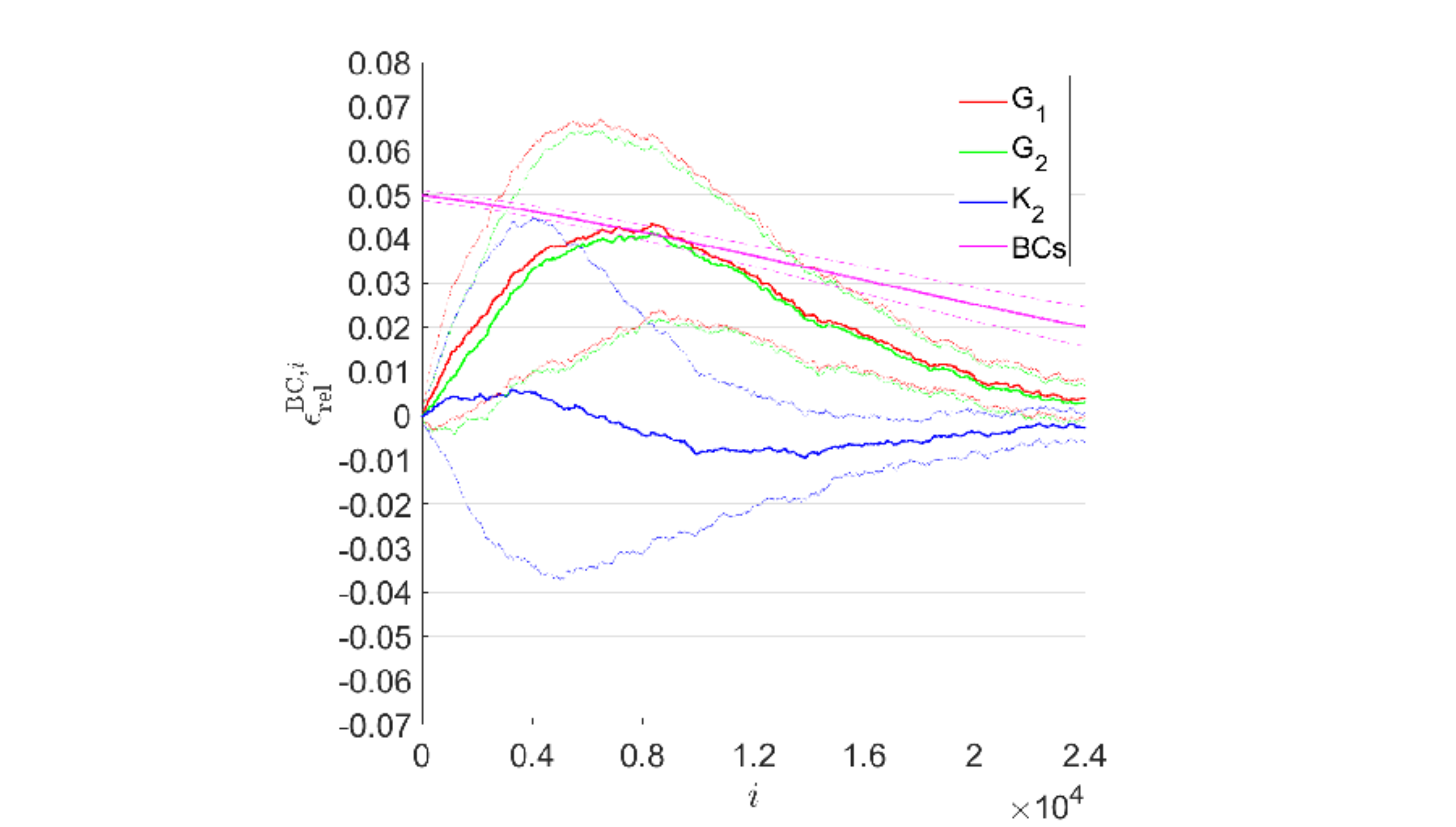}}
\caption{The evolution of the relative error in the identified material parameters and boundary conditions, cf. Eqs.~\eqref{error-total-mat} and~\eqref{error-total-kin}, as a function of MHA step $i$ for the tensile test, the means across all sampled values of all iterations (after burn-in) are plotted with solid lines and are complemented with $\pm$ standard deviations (dashed lines), $t=2\times244$ (100\% of total). Identified parameter $K_1$ is used as a normalization factor, and is not fixed during identification.}
\label{fig:MH_steps_sigma_tension_comparison}
\end{figure}

The convergence speed is higher for all the considered parameters, including the boundary conditions, where the minimum achieved average error is roughly two times smaller than in the normalized version of MHA, see Fig.~\ref{fig:MH_steps_sigma_tension_comparison}.

\begin{figure}[!htbp]
    \centering
	\includegraphics[width=0.4\textwidth]{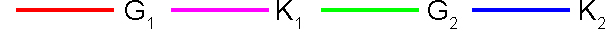}\\
	\subfloat[$t=2\times244$ (100\% of total).]{
	\includegraphics[trim=0cm 0cm 0cm 0cm, clip=true, width=0.49\textwidth]{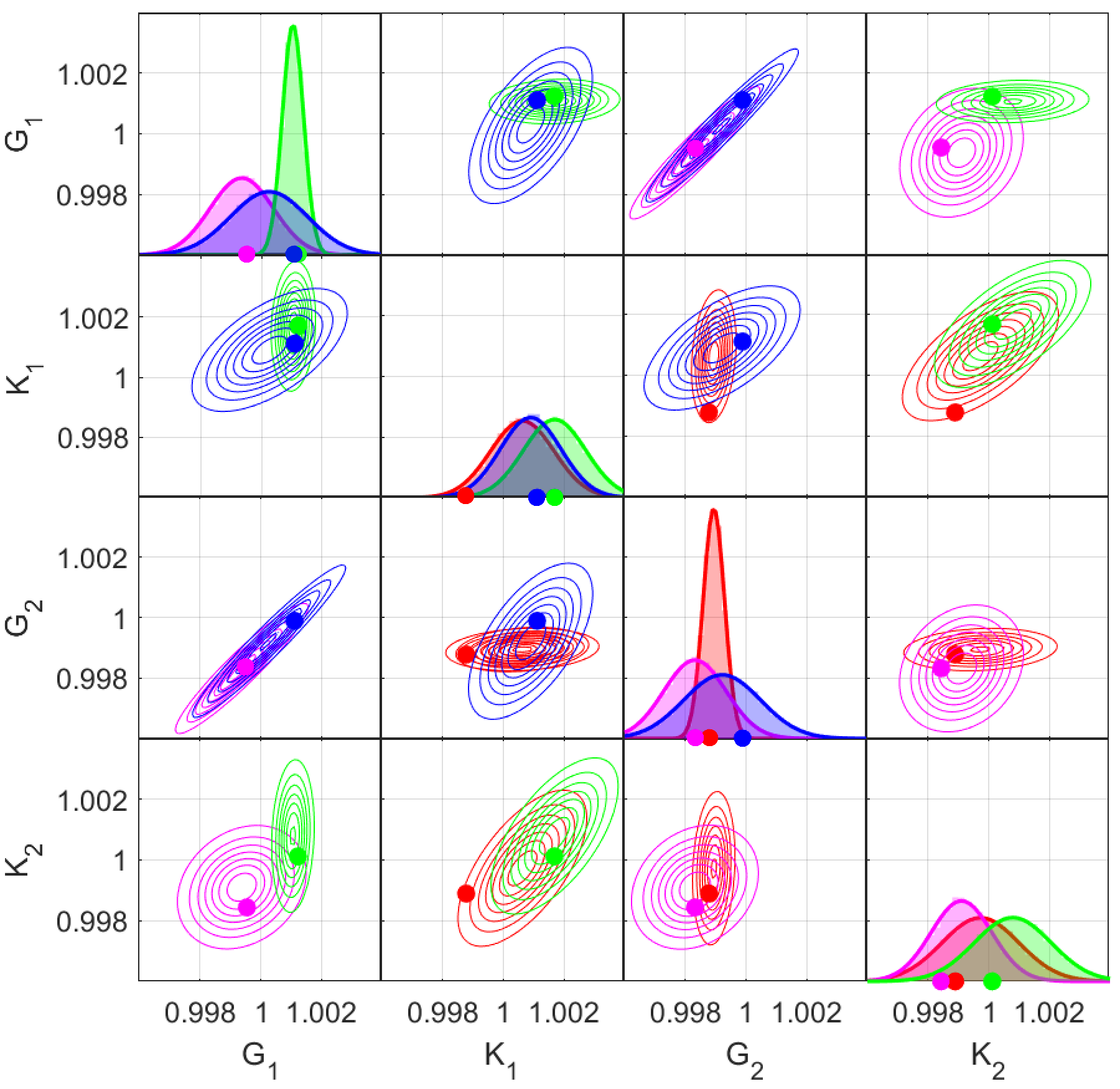}
    \label{fig:hist_ill_q=1.0_tension}}
    \subfloat[$t=2\times61$ (25\% of total).]{
	\includegraphics[trim=0cm 0cm 0cm 0cm, clip=true, width=0.49\textwidth]{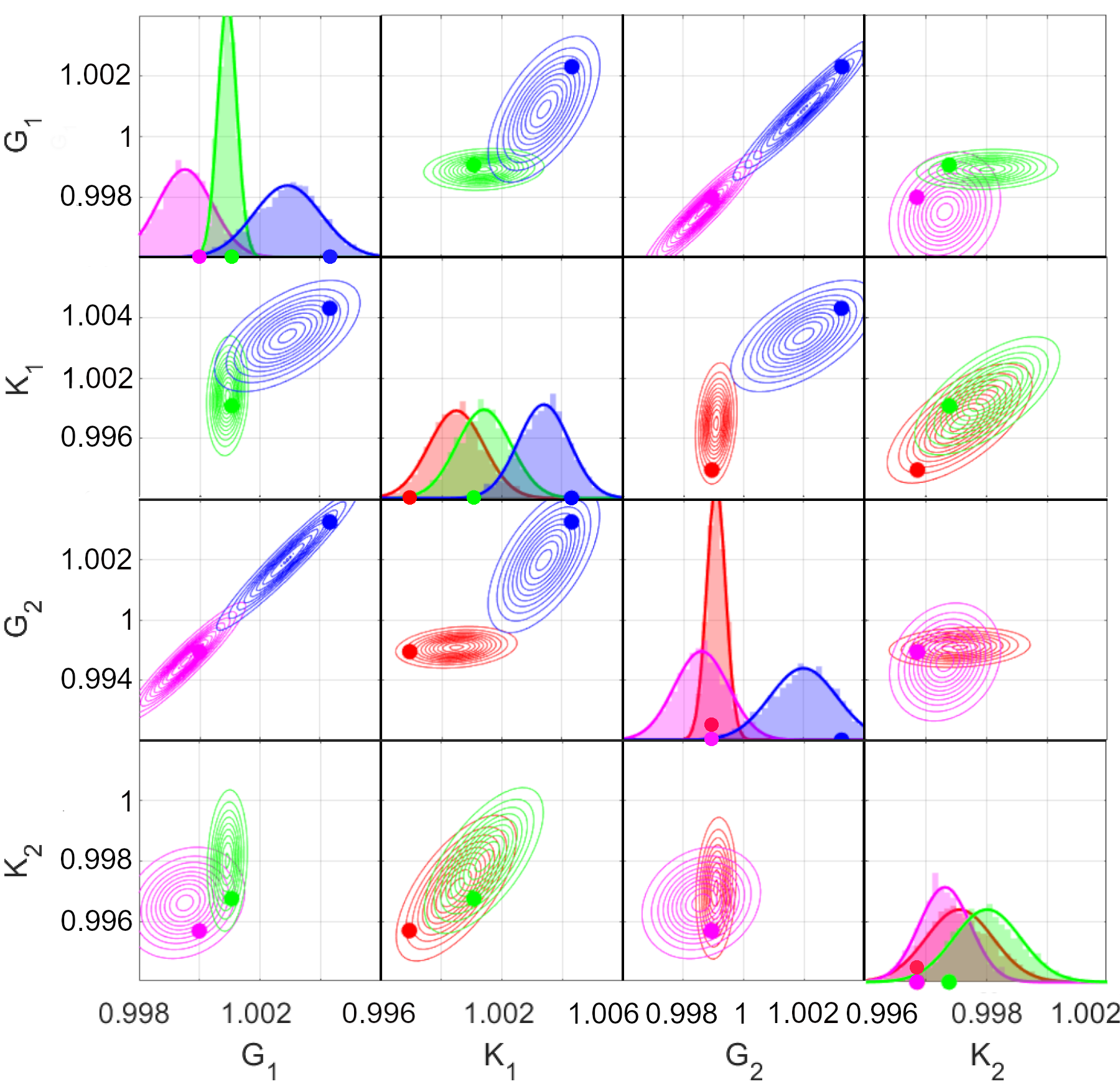}
    \label{fig:hist_ill_q=0.25_tension}}
    \caption{Posterior probability distributions for the material parameters for the non-normalized MHA and the shear test, compared to the point estimations given by BE-IDIC (denoted by dots), depending on the choice of the normalization parameter color coded according to the legend at the top.}
    \label{fig:hist_ill_tension}
\end{figure}

The most interesting feature of the non-normalized MHA is that the posterior distributions can be obtained for all the parameters by choosing any of the remaining parameters as the normalization factor, as shown in Fig.~\ref{fig:hist_ill_tension}, where the estimated posterior probability distribution resulting from the chosen normalization parameter is color-coded, as well as the points given by the deterministic BE-IDIC with the same choice for the fixed parameter. The probability density estimations for each of the parameters differ in bias and variance, based on the choice of the normalization parameter. Similarly to the normalized MHA, a relationship between the results of the deterministic method and the modes of the probability distributions can be observed, as they generally appear in the same order. Indeed, in practical applications, one might not have estimates for all the material parameters beforehand to be as free in the choice of the normalization. However, the parameter ratios can still be effectively established with the non-normalized MHA.

\subsection{Shear Test}

\begin{figure}[!htbp]
	\centering
		\subfloat[$t=2\times244$ DOFs (100\% of total).]{
		\includegraphics[trim=6.5cm 0cm 7.5cm 1cm, clip=true, width=0.4\textwidth]{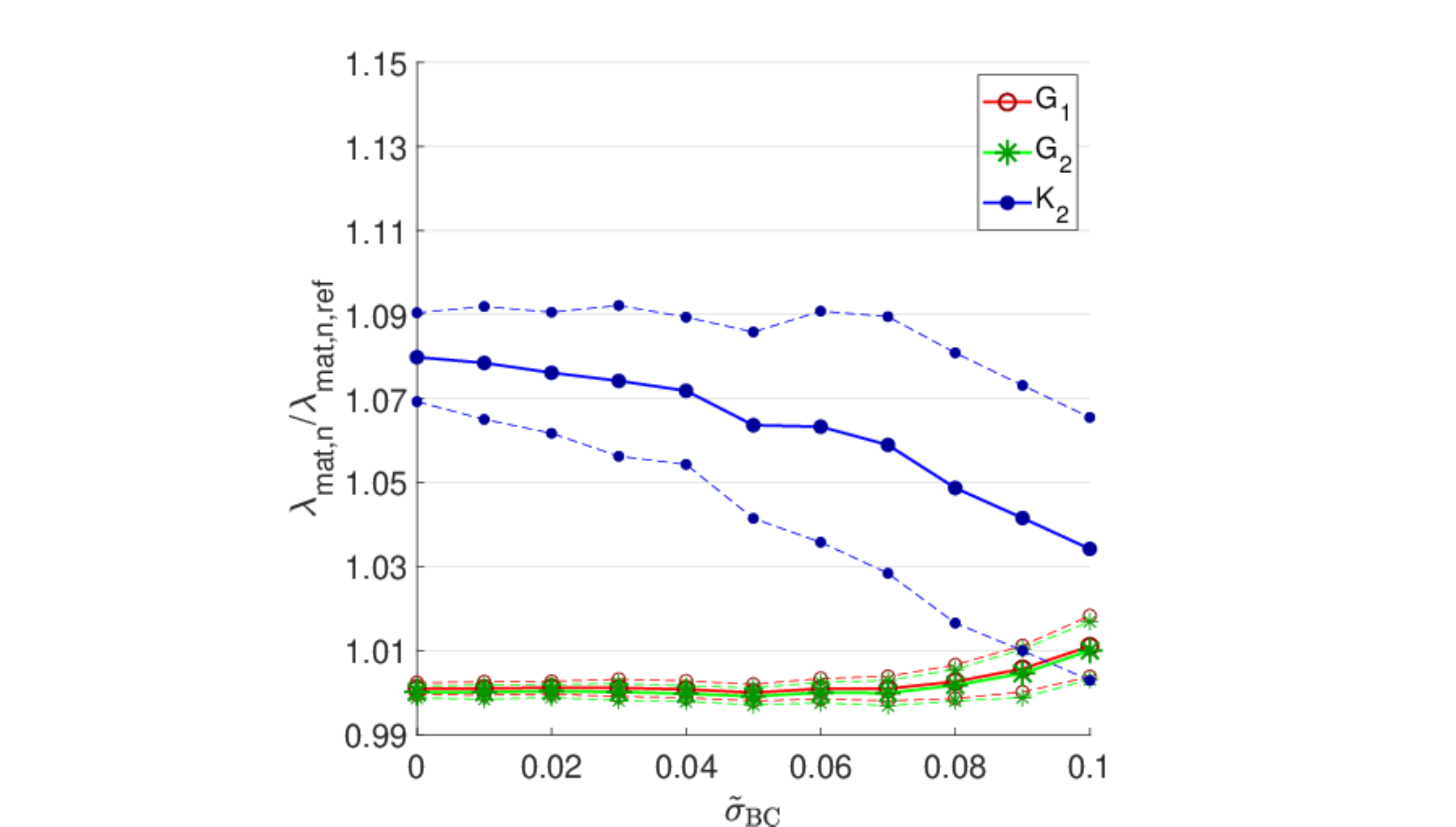}}
		\subfloat[$t=2\times61$ DOFs (25\% of total).]{
		\includegraphics[trim=6.5cm 0cm 7.5cm 1cm, clip=true, width=0.4\textwidth]{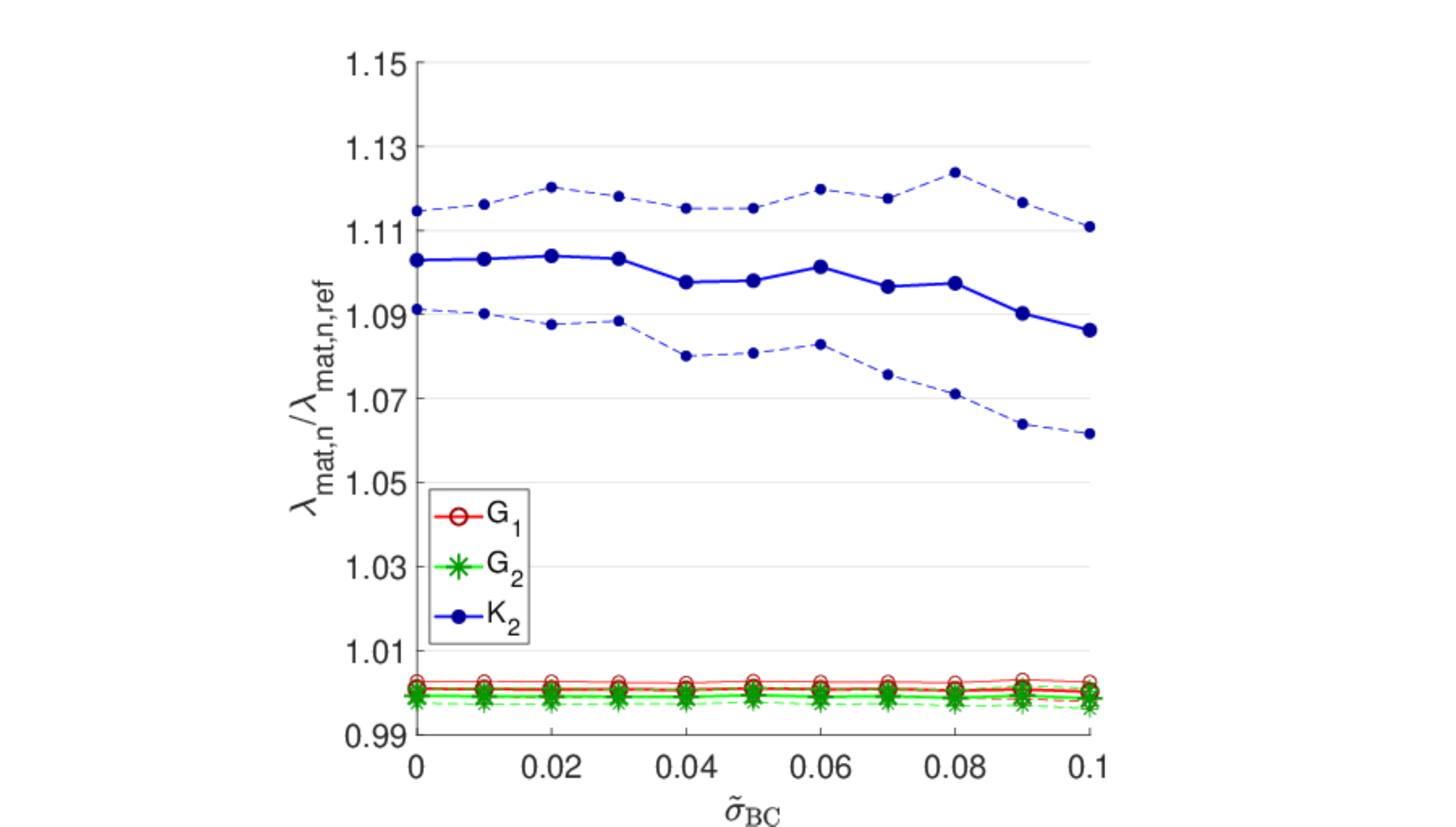}}
	\caption{Averaged modes of normalized material parameter distributions $\lambda_{\mathrm{mat},n}/\lambda_{\mathrm{mat},n,\mathrm{ref}}$ for relaxed boundary conditions, cf. Eq.~\eqref{noise_imp}, with the starting noise standard deviation $\tilde{\sigma}_{\mathrm{bc}}$ for the shear test. The means across all sampled values of all iterations (after burn-in) are plotted with solid lines and are complemented with $\pm$ standard deviations (dashed lines) with $N=24\,000$ and burn-in $N_0=22\,000$ steps. Identified parameter $K_1$ is used as a normalization factor, and is not fixed during identification.}
	\label{fig:mats_FE_shear_ill}
\end{figure}

\begin{figure}[!htbp]
\centering
	\subfloat[Normalized MHA, $\tilde{\sigma}_{\mathrm{bc}}=0$.]{
	\includegraphics[trim=5.8cm 0cm 7.5cm 0.5cm, clip=true, width=0.4\textwidth]{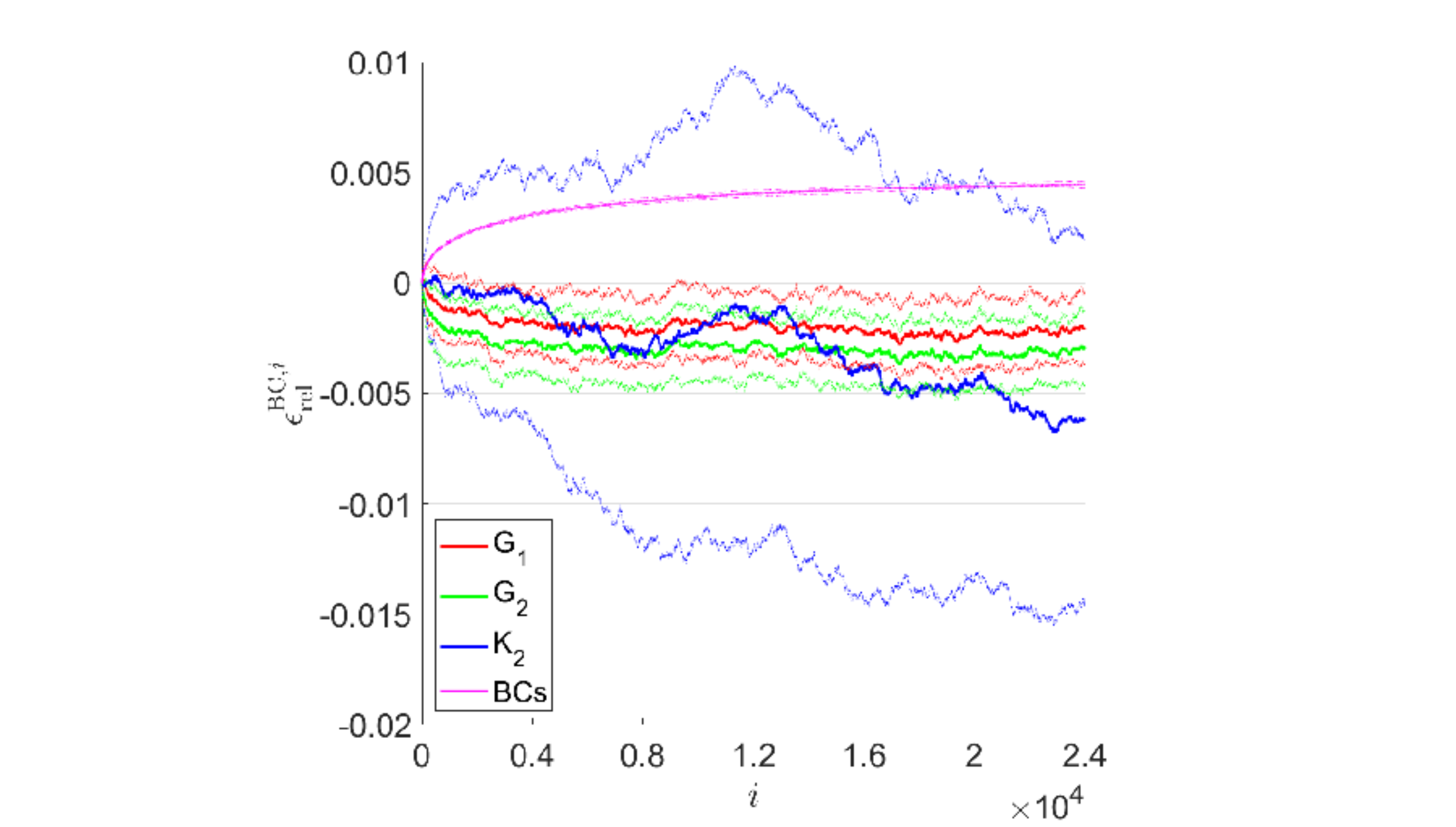}}
	\subfloat[Non-normalized MHA, $\tilde{\sigma}_{\mathrm{bc}}=0$.]{
	\includegraphics[trim=5.8cm 0cm 7.5cm 0.5cm, clip=true, width=0.4\textwidth]{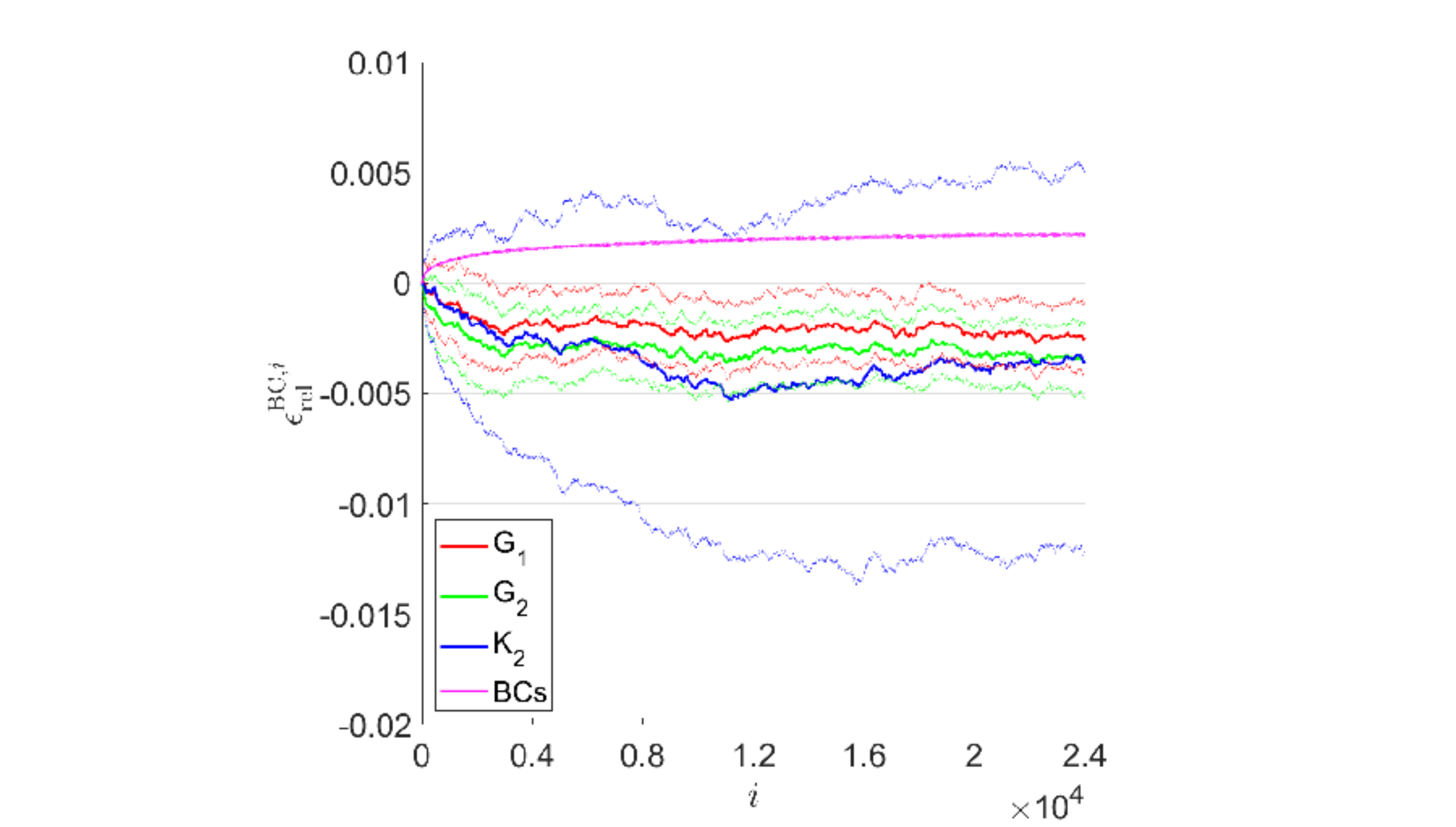}}\\
	\subfloat[Normalized MHA,  $\tilde{\sigma}_{\mathrm{bc}}=0.1$.]{
	\includegraphics[trim=5.8cm 0cm 7.5cm 0.5cm, clip=true, width=0.4\textwidth]{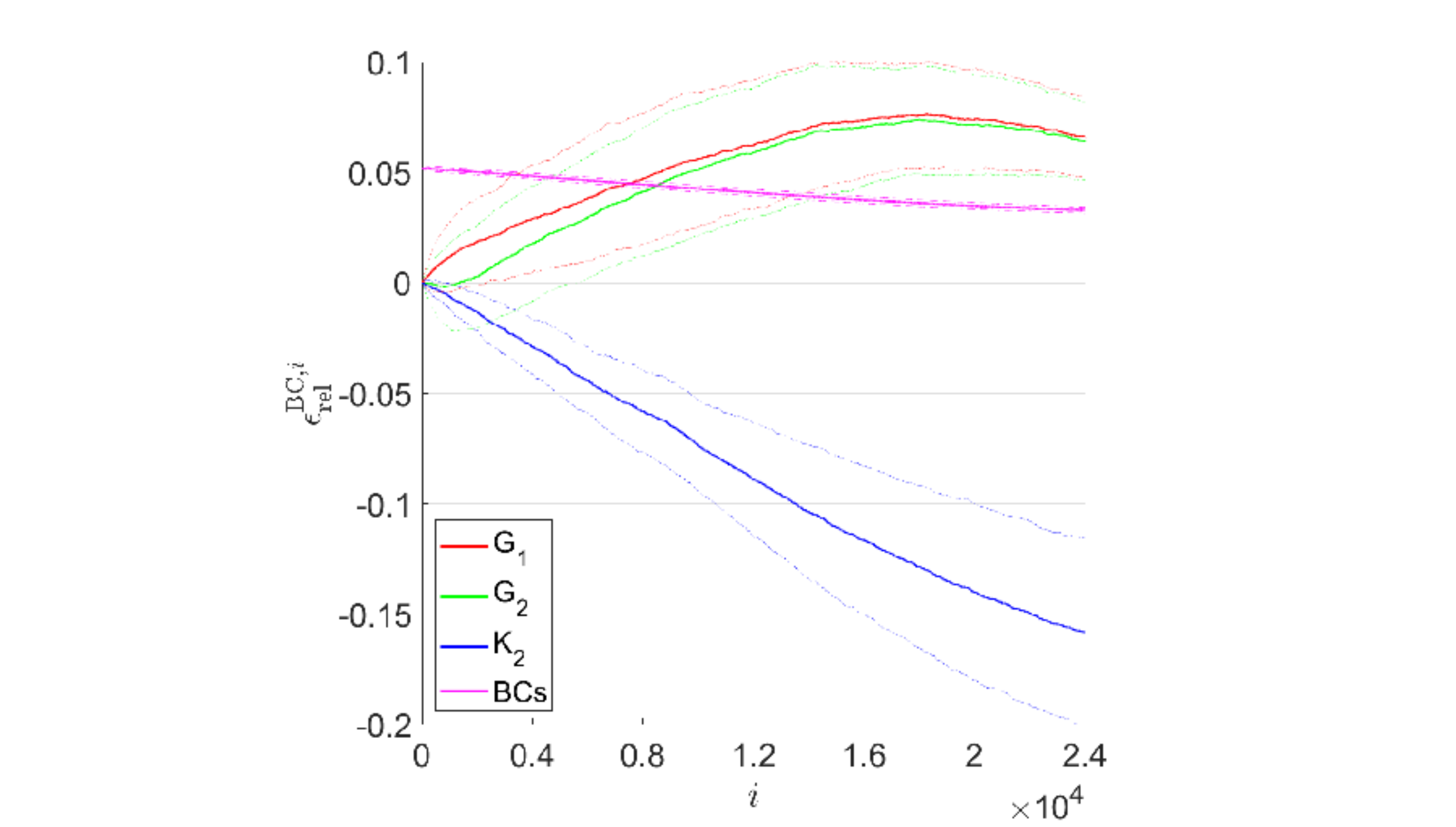}}
    \subfloat[Non-normalized MHA, $\tilde{\sigma}_{\mathrm{bc}}=0.1$.]{
	\includegraphics[trim=5.8cm 0cm 7.5cm 0.5cm, clip=true, width=0.4\textwidth]{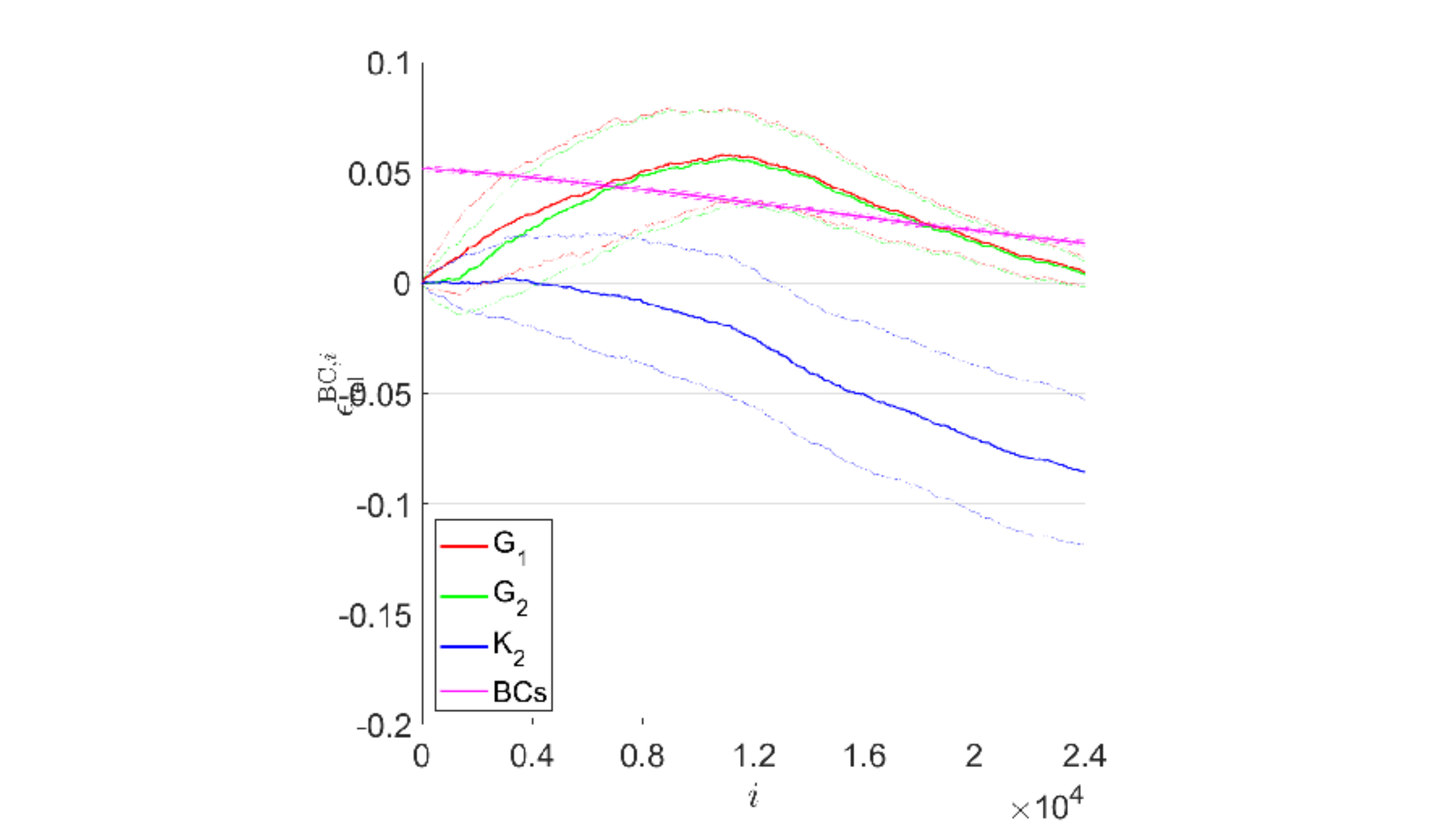}}
\caption{The evolution of the relative error in the identified material parameters and boundary conditions, cf. Eqs.~\eqref{error-total-mat} and~\eqref{error-total-kin}, as a function of MHA step $i$ for the shear test, the means across all sampled values of all iterations (after burn-in) are plotted with solid lines and are complemented with $\pm$ standard deviations (dashed lines), $t=2\times244$ (100\% of total). Identified parameter $K_1$ is used as a normalization factor, and is not fixed during identification.}
\label{fig:MH_steps_sigma_shear_comparison}
\end{figure}

\begin{figure}[!htbp]
    \centering
	\includegraphics[width=0.4\textwidth]{pics_new/legend.jpg}\\
    \subfloat[$t=2\times244$ kinematic DOFs (100\% of total).]{
	\includegraphics[trim=0cm 0cm 0cm 0cm, clip=true, width=0.49\textwidth]{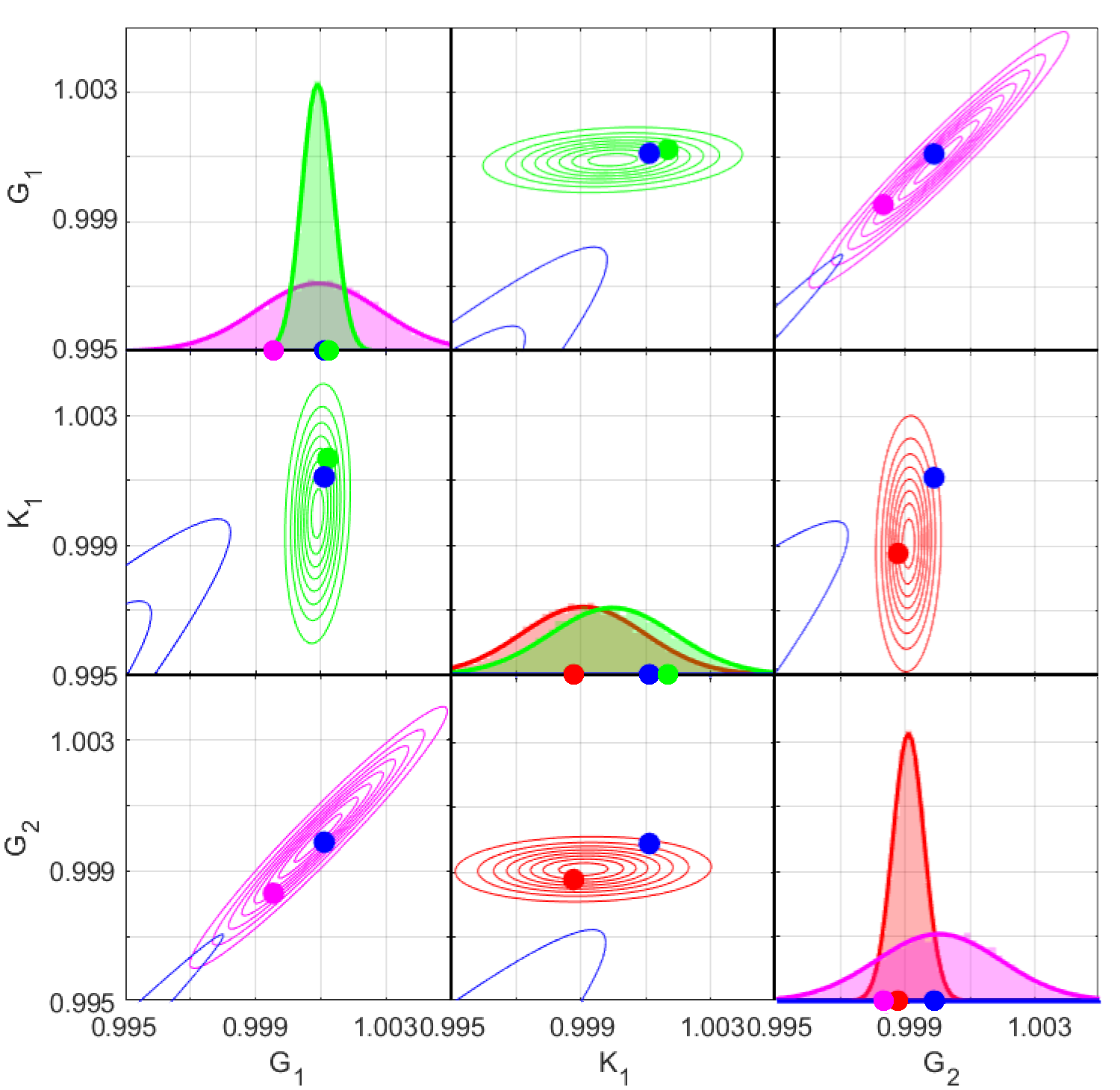}
    \label{fig:hist_ill_q=1.0_shear}}
    \subfloat[$t=2\times61$ kinematic DOFs (25\% of total).]{
	\includegraphics[trim=0cm 0cm 0cm 0cm, clip=true, width=0.49\textwidth]{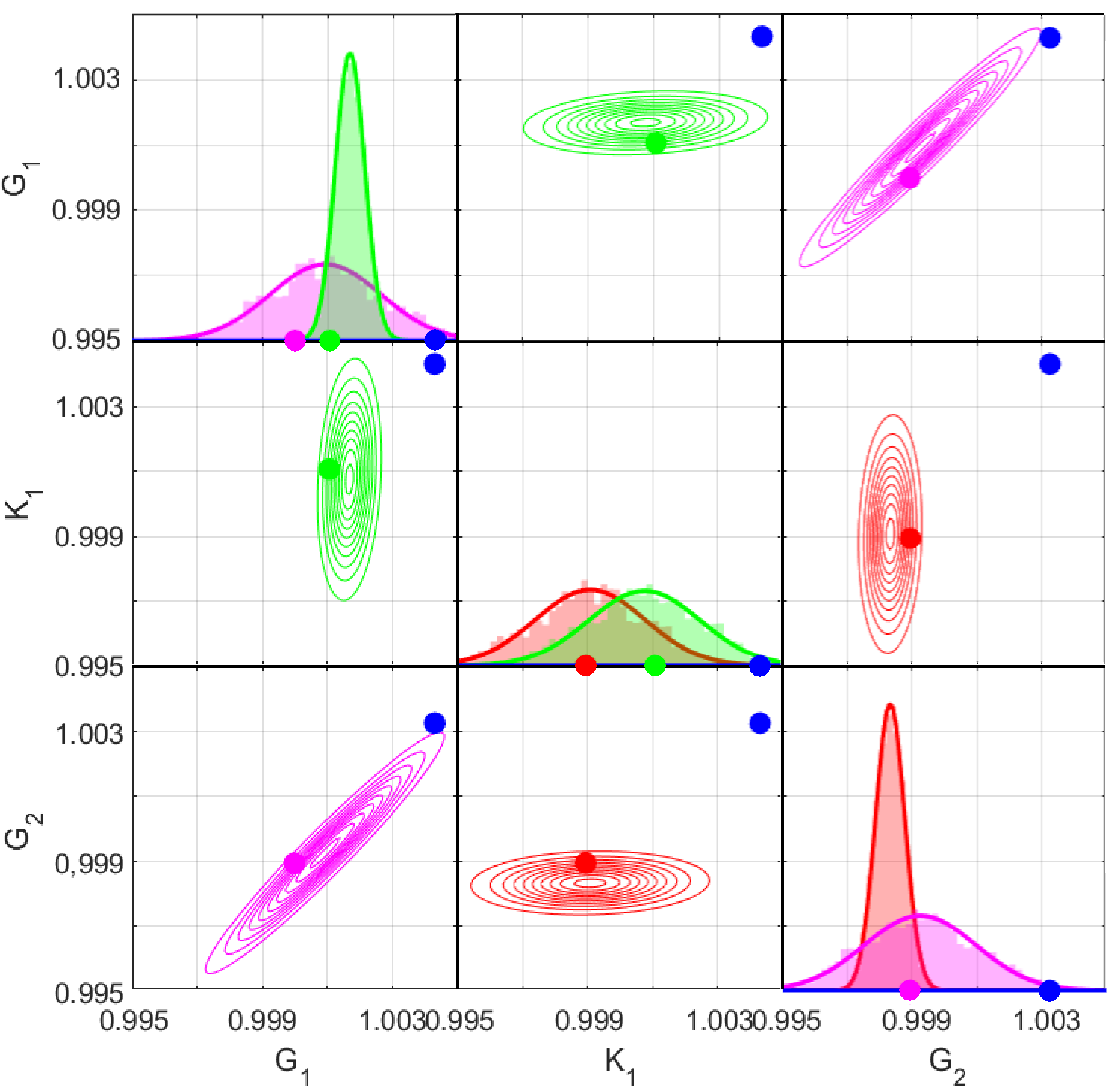}
	\label{fig:hist_ill_q=0.25_shear}}
    \caption{Posterior probability distributions for the material parameters for the non-normalized MHA and the shear test, compared to the point estimations given by BE-IDIC (denoted by dots), depending on the choice of the normalization parameter color coded according to the legend at the top.}
\end{figure}

The shear test continues to be not sensitive enough to correctly establish the parameter $K_2$ (see Fig. \ref{fig:mats_FE_shear_ill}), even though the identification accuracy in other material parameters and boundary conditions increased significantly compared to the normalized version, as shown in Fig.~\ref{fig:MH_steps_sigma_shear_comparison}. The resulting material parameter estimates have a wider confidence interval than in the tensile test, and at the same time the biases remain comparable to those of BE-IDIC, see Fig.~\ref{fig:hist_ill_q=0.25_shear}. The distributions of the parameter $K_2$ are omitted because of the low sensitivity.

\section{Summary and Conclusions}
\label{sec:summary}
In this contribution, the performance and the probabilistic robustness of the Metropolis--Hastings algorithm (MHA) is compared to the deterministic Integrated Digital Image Correlation (IDIC) method. To this end, a heterogeneous microstructural specimen with a random distribution of circular inclusions have been subjected to two virtual mechanical tests under plane strain conditions, one to primarily introduce tension, the other to introduce shear. The goal was to identify parameters of two distinct groups: materials and applied boundary conditions. The effect of errors in the applied boundary conditions was studied. First, the MHA that only identifies the material parameters with fixed boundary conditions was considered, and its sensitivity with respect to random and systematic errors in the boundary conditions was quantified and compared to the IDIC. MHA's parameter field was then expanded with two different ways of approximating the boundary conditions, and the method was compared to the Boundary-Enriched IDIC (BE-IDIC). The experiments have shown a similar behavior of MHA with fixed boundary conditions to IDIC for both the systematic as well as random error in the applied boundary conditions.

A possible way of reducing the dimensionality of the boundary condition parameters were suggested: substituting the employed discretization mesh on the boundary with a coarser one. The robustness test for the random error in the boundary conditions has shown a higher convergence rate for the price of introducing an acceptable systematic error in the material parameters.

The deterministic approach requires one of the material parameters to be fixed at the exact value, because the inverse problem with Dirichlet boundary conditions is inherently ill-posed. The benefit of the stochastic approach, on the other hand, is that this normalization, while possible, is unnecessary. It was shown that the non-normalized approach to the parameter identification, where none of the material parameters is fixed, converges faster and is more robust with respect to the initial boundary noise. As a result, this version of the algorithm was able to handle high-fidelity boundary conditions burdened with noise more efficiently than the normalized version. The material parameter estimates can then be obtained in the post-processing step from the material ratios established by MHA, which allows solving ill-posed or weakly conditioned structural optimization problems such as parameter identification in laminates \citep[e.g.,][]{CHEN2021113853}, where deterministic methods are not applicable. 

Overall, MHA with relaxed boundary conditions proved to be slightly less accurate and more computationally costly than BE-IDIC in finding the high probability region. However, both normalized and non-normalized versions of MHA can identify the material parameters with the accuracy similar to the deterministic methods in the tensile test. The irreducible errors in both the stochastic and deterministic method stem from the same source and affect the results almost equivalently. The stochastic approach was shown to have a few advantages: the relative ease of implementation (low number of hyper-parameters, usually set according to a rule-of-thumb), ability to optimize a large number of parameters that can also be dependent, as well as statistical data allowing more insight in the relationships between the parameters. The main downside is the high computational cost and lower precision than for the existing deterministic methods, which can be potentially remedied with an adaptive approach and the use of surrogate models.

\section*{CRediT Author Statement}
LG: Methodology, Software, Investigation, Writing -- Original Draft, Visualization; OR: Conceptualization, Methodology, Software, Resources, Writing -- Review \& Editing, Supervision; IP: Methodology, Writing -- Review \& Editing, Supervision; JH: Methodology, Software, Investigation, Writing -- Review \& Editing; JZ: Methodology, Writing -- Review \& Editing, Funding acquisition. 

\section*{Acknowledgements}
This work received the support from the European Regional Development Fund (Center of Advanced Applied Sciences – CAAS, CZ 02.1.01/0.0/0.0/16 019/0000778 (IP and LG)), the Czech Science Foundation (projects No. 22-35755K (LG) and No. 19-26143X (OR, JH, and JZ)), and the Student Grant Competition of CTU (projects No. SGS21/004/OHK1/1T/11 (IP) and No. SGS23/002/OHK1/1T/11 (LG)). 

\bibliography{mybibfile}

\end{document}